\catcode`\@=11

\newcount\secno
\newcount\subsecno
\newcount\subsubsecno
\newcount\prmno
\newcount\exno

\newbox\boxa
\newbox\boxb
\newdimen\dima
\newwrite\toc
\newwrite\lbl
\newread\testfile

\newif\ifnotfound
\newif\iffound
\newif\ifignore

\newif\ifpagetitre        \pagetitretrue
\newtoks\hautpagetitre    \hautpagetitre={\hfil}
\newtoks\baspagetitre     \baspagetitre={\hfil}
\newtoks\auteurcourant    \auteurcourant={\hfil}
\newtoks\titrecourant     \titrecourant={\hfil}
\newtoks\hautpagegauche   \newtoks\hautpagedroite
\hautpagegauche={\hfil\tensl\the\auteurcourant\hfil}
\hautpagedroite={\hfil\tensl\the\titrecourant\hfil}

\newtoks\baspagegauche
\baspagegauche={\hfil\tenrm\folio\hfil}
\newtoks\baspagedroite
\baspagedroite={\hfil\tenrm\folio\hfil}

\headline={\ifpagetitre\the\hautpagetitre
\else\ifodd\pageno\the\hautpagedroite
\else\the\hautpagegauche\fi\fi}

\footline={\ifpagetitre\the\baspagetitre
\global\pagetitrefalse
\else\ifodd\pageno\the\baspagedroite
\else\the\baspagegauche\fi\fi}

\input amssym.def
\input amssym

\font\capit=cmcsc10
\font\prm=cmr8

\font\Bf=cmbx12 scaled\magstep1
\font\goth=eufm10

\font\teneusm=eusm10 
\font\seveneusm=eusm7 
\font\fiveeusm=eusm5 
\newfam\eusmfam 
\textfont\eusmfam=\teneusm 
\scriptfont\eusmfam=\seveneusm
\scriptscriptfont\eusmfam=\fiveeusm

\let\cal\eusm

\input diagrams.tex
\def\hfl#1#2{\rTo^{#1}_{#2}}
\def\vfl#1#2{\dTo^{#1}_{#2}}

\def\nefl#1#2{\ruTo^{#1}_{#2}}

\def\sefl#1#2{\rdTo^{#1}_{#2}}
\def\swfl#1#2{\ldTo^{#1}_{#2}}																								

\catcode`\;=\active
\def;{\relax\ifhmode\ifdim\lastskip>\z@
\unskip\fi\kern.2em\fi\string;}

\catcode`\:=\active
\def:{\relax\ifhmode\ifdim\lastskip>\z@\unskip\fi
\penalty\@M\ \fi\string:}

\catcode`\!=\active
\def!{\relax\ifhmode\ifdim\lastskip>\z@
\unskip\fi\kern.2em\fi\string!}

\catcode`\?=\active
\def?{\relax\ifhmode\ifdim\lastskip>\z@
\unskip\fi\kern.2em\fi\string?}

\toksdef\ta=0 \toksdef\tb=2
\long\def\leftappenditem#1\to#2{\ta={\\{#1}}\tb=\expandafter{#2}%
                                \edef#2{\the\ta\the\tb}}
\long\def\rightappenditem#1\to#2{\ta={\\{#1}}\tb=\expandafter{#2}%
                                \edef#2{\the\tb\the\ta}}

\def\lop#1\to#2{\expandafter\lopoff#1\lopoff#1#2}
\long\def\lopoff\\#1#2\lopoff#3#4{\def#4{#1}\def#3{#2}}

\def\ismember#1\of#2{\foundfalse{\let\given=#1%
    \def\\##1{\def\next{##1}%
    \ifx\next\given{\global\foundtrue}\fi}#2}}

\def\namedef#1{\expandafter\def\csname #1\endcsname}
\def\nameuse#1{\csname #1\endcsname}
\def\typeout#1{\immediate\write16{#1}}
\long\def\ifundefined#1#2#3{\expandafter\ifx\csname 
  #1\endcsname\relax#2\else#3\fi}
\def\hwrite#1#2{{\let\the=0\edef\next{\write#1{#2}}\next}}
\def\lookatfile#1{\openin\testfile=\jobname.#1
    \ifeof\testfile{\immediate\openout\nameuse{#1}\jobname.#1
                    \write\nameuse{#1}{}
                    \immediate\closeout\nameuse{#1}}\fi%
    \immediate\closein\testfile}%

\let\@end\end

\def\begin#1{\ifundefined{@#1}%
                 {\def\tempa{\typeout{begin{#1} is undefined}\endgroup}}%
                 {\def\tempa{\def\@currenvir{#1}\csname @#1\endcsname}}%
             \begingroup\tempa}%
\def\end#1{\csname end#1\endcsname\endgroup%
           \ifignore\global\ignorefalse\ignorespaces\fi}

\def\@currenvir{document}
\def\openall{\openout\toc=\jobname.toc
             \openout\lbl=\jobname.lbl}
\def\closeall{\closeout\toc
              \closeout\lbl}
\def\@document{\lookatfile{lbl}
               \input\jobname.lbl
               \lookatfile{toc}
               \input\jobname.toc
               \openall}
\def\enddocument{\endgroup\closeall\@end}

\let\em\sl
\def\@em{\em}

\def\@proof{{\it D\'emonstration.}\ }

\def\@centering{\leftskip=0pt plus 1fil
                \rightskip=0pt plus 1fil
                \parfillskip=0pt
                \noindent
                \obeylines
                \let\\=\break
                }

\def\@displaymath{$$}

\def\prmlist{}
\def\bignewtheorem#1#2{\namedef{@#1}{\global\advance\prmno by 1
        \bigbreak{\capit #2 \the\secno.\the\prmno.\ ---}\it}
        \rightappenditem{#1}\to\prmlist
        \namedef{end#1}{\bigbreak\global\ignoretrue}}
\def\mednewtheorem#1#2{\namedef{@#1}{\global\advance\prmno by 1
        \bigbreak{\capit #2 \the\secno.\the\prmno.\ ---}\it}
        \rightappenditem{#1}\to\prmlist
        \namedef{end#1}{\bigbreak\global\ignoretrue}}
\def\smallnewtheorem#1#2{\namedef{@#1}{\global\advance\prmno by 1
        \bigbreak{\capit #2 \the\secno.\the\prmno.\ ---}\it}
        \rightappenditem{#1}\to\prmlist
        \namedef{end#1}{\bigbreak\global\ignoretrue}}
\def\nsmallnewtheorem#1#2{\namedef{@#1}{\global\advance\prmno by 1
        \bigbreak{\capit #2 \the\secno.\the\prmno.\ ---}}
        \rightappenditem{#1}\to\prmlist
        \namedef{end#1}{\bigbreak\global\ignoretrue}}
\def\exlist{}
\def\exercice#1#2{\namedef{@#1}{\global\advance\exno by 1
        \medbreak{\rm #2 \the\secno.\the\exno}\rm}
        \rightappenditem{#1}\to\exlist
        \namedef{end#1}{\bigbreak\global\ignoretrue}}
\exercice{ex}{Exercice}

\smallnewtheorem{defi}{D\'efinition}
\bignewtheorem{th}{Th\'eor\`eme} 
\mednewtheorem{lemme}{Lemme}
\mednewtheorem{prop}{Proposition}
\mednewtheorem{cor}{Corollaire}
\nsmallnewtheorem{rem}{Remarque}

\def\newcommand#1#2{\def#1{#2}}

\let\array\matrix
\def\matrix#1{\let\\=\cr\array{#1}}
\def\determinant#1{\let\\=\cr\left|\array{#1}\right|}
\def\@array{\def\get##1\end##2{\let\\=\cr\array{##1}\endgroup}\get}
\def\@determinant{\def\get##1\end##2{\determinant{##1}\endgroup}\get}
\def\@matrix{\def\get##1\end##2{\matrix{##1}\endgroup}\get}

\def\vspace#1{\vskip#1}
\def\hspace#1{\leavevmode\expandafter\hskip #1\relax}

\def\title#1{
             \begin{centering}\Bf{#1}\end{centering}}
\def\sectiontitle#1{%
      \hwrite\toc{\string\newtoc{{#1}\string\dotfill\the\pageno}}
                    \begin{centering}\bf{#1}\end{centering}}
\def\subsectiontitle#1{{\bf{#1}}
                       \hwrite\toc{\string\newtoc
                     {{\indent{\prm #1}\string\dotfill\the\pageno}}}}
\def\subsubsectiontitle#1{\begin{centering}{\it{#1}}\end{centering}}
\def\exercicetitle#1{{\bf{#1}}
                       \hwrite\toc{\string\newtoc
                     {{\hskip0.5cm#1\string\dotfill\the\pageno}}}}           

\def\section#1{\bigbreak\bigskip
               \def\@currenvir{section}
               \global\advance\secno by1
               \global\subsecno=0\global\subsubsecno=0\global\prmno=0
               \global\exno=0
               \sectiontitle{\number\secno.\ #1}
               \medskip}
\def\subsection#1{\bigbreak
                  \def\@currenvir{subsection}
                  \global\advance\subsecno by 1
                  \global\subsubsecno=0
                  \subsectiontitle{\number\secno.\number\subsecno.\ #1}
                  \medskip}
\def\subsubsection#1{\bigbreak
                     \def\@currenvir{subsubsection}
                     \global\advance\subsubsecno by1
                     \subsubsectiontitle{#1}             
\medbreak}

\def\exercice#1{\medbreak
               \def\@currenvir{exercice}
               \global\advance\secno by1
               \global\subsecno=0\global\subsubsecno=0\global\prmno=0
               \exercicetitle{\number\secno.\number\exno\ #1}
               \medskip}

\def\isinlabellist#1\of#2{\notfoundtrue%
   {\def\given{#1}%
    \def\\##1{\def\next{##1}%
    \lop\next\to\za\lop\next\to\zb%
    \ifx\za\given{\zb\global\notfoundfalse}\fi}#2}%
    \ifnotfound{\immediate\write16%
                 {Warning - [Page \the\pageno] No reference found}}%
                \fi}%
\def\ref#1{\ifx\labellist\empty{\immediate\write16
                 {Warning - [Page \the\pageno] No references found at all.}}
               \else{\isinlabellist{#1}\of\labellist}\fi}

\def\newlabel#1#2{\rightappenditem{\\{#1}\\{#2}}\to\labellist}
\def\labellist{}

\def\label#1{\ismember{\@currenvir}\of\prmlist%
 \iffound{\hwrite\lbl{\string\newlabel{#1}{\number\secno.\number\prmno}}}\fi%
  \def\given{section}%
  \ifx\given\@currenvir%
    {\hwrite\lbl{\string\newlabel{#1}}}\fi%
  \def\given{subsection}%
  \ifx\given\@currenvir%
    {\hwrite\lbl{\string\newlabel{#1}{\number\secno.\number\subsecno}}}\fi%
  \def\given{subsubsection}%
  \ifx\given\@currenvir%
  {\hwrite\lbl{\string%
    \newlabel{#1}{\number\secno.\number\subsecno.\number\subsubsecno}}}\fi
}

\def\newtoc#1{\rightappenditem{#1}\to\toclist}
\def\toclist{}
\def\writetoclist{{\def\\##1{\def\next{##1}\noindent\next\par}\toclist}}
\def\tableofcontents{\ifx\toclist\empty{\immediate\write16
                 {Warning: [Page \the\pageno] No Table of contents found}}
                          \else{\bigskip                               
\centerline{\capit Table des mati\`eres} 
                                \medskip
                                \writetoclist
                                \bigskip}\fi}

 \def\getlength#1{\ifx#1\end \let\next=\relax
  \else\advance\count0 by 1 \let\next=\getlength\fi\next}

\catcode`\@=12

\def\cc#1{\hfill\kern .7em#1\kern .7em\hfill}
\def\scc#1{\hfill\kern .2em#1\kern .2em\hfill}

\def\cf{{\it cf.\/}\ }
\def\ie{{\it i.e.\/}\ }

\def\up#1{\raise 1ex\hbox{\sevenrm#1}}

\def\cqfd{\unskip\kern 6pt\penalty 500
          \hbox{\vrule\vbox to 4pt{\hrule width 3pt\vfill
          \hrule}\vrule}\par}
\def\nl{\hfill\break}
\let\\\nl

\def\date{le {\the\day}
          \ifcase\month
          \or janvier
          \or f\'evrier
          \or mars
          \or avril
          \or mai
          \or juin
          \or juillet
          \or août
          \or septembre
          \or octobre
          \or novembre
          \or d\'ecembre
          \fi
          {\the\year}} 

\newbox\boxa

\def\boxit#1{\setbox\boxa=\hbox{#1}%
  \setbox\boxa=\hbox{\vrule\box\boxa\vrule}%
  \setbox\boxa=\vbox{\hrule\box\boxa\hrule}%
  \box\boxa\relax}

\let\newdiagram\diagram
\def\diagram#1{\newdiagram #1 \enddiagram}

\def\og{\leavevmode\raise.3ex\hbox{$\scriptscriptstyle\langle\!\langle$}}
\def\fg{\leavevmode\raise.3ex\hbox{$\scriptscriptstyle\,\rangle\!\rangle$}}
\mathcode`A="7041 \mathcode`B="7042 \mathcode`C="7043 \mathcode`D="7044
\mathcode`E="7045 \mathcode`F="7046 \mathcode`G="7047 \mathcode`H="7048
\mathcode`I="7049 \mathcode`J="704A \mathcode`K="704B \mathcode`L="704C
\mathcode`M="704D \mathcode`N="704E \mathcode`O="704F \mathcode`P="7050
\mathcode`Q="7051 \mathcode`R="7052 \mathcode`S="7053 \mathcode`T="7054
\mathcode`U="7055 \mathcode`V="7056 \mathcode`W="7057 \mathcode`X="7058
\mathcode`Y="7059 \mathcode`Z="705A

\font\tentr=msbm10
\font\seventr=msbm7
\font\fivetr=msbm5
\newfam\trfam
\textfont\trfam\tentr
\scriptfont\trfam\seventr
\scriptscriptfont\trfam\fivetr
\def\tr{\fam\trfam\tentr}
\def\comp{{\tr C}}          
\def\proj{{\tr P}}          
\def\entr{{\tr Z}}          

\def\mbox#1{\leavevmode\hbox{#1}}
\def\text#1{\leavevmode\hbox{#1}}
\def\Ext{\hbox{\rm Ext}}                    
\def\Tor{\hbox{\rm Tor}}                    
\def\Hom{\hbox{\rm Hom}}                    

\def\frac#1#2{{#1\over #2}}
\long\def\raisebox#1#2{\leavevmode\hbox{\raise #1\hbox{#2}}}

\def\mbox#1{\leavevmode\hbox{\rm #1}}
\def\text#1{\leavevmode\hbox{\rm #1}}
\def\ker{\hbox{\rm ker}\thinspace}
           \def\Aut{\hbox{\rm Aut}}

\def\grass{\hbox{\rm Grass}} 
	
\def\Hilb{\hbox{\rm Hilb}}
\def\Pic{\hbox{\rm Pic}}
\def\Div{\hbox{\rm Div}}

\def\Spec{\hbox{\rm Spec}}
\def\cqfd{\unskip\kern 6pt\penalty 500
          \hbox{\vrule\vbox to 4pt{\hrule width 3pt\vfill
          \hrule}\vrule}\par}

\def\coker{\hbox{\rm coker\ }}
\def\syst{\hbox{\rm Syst}}

\def\mod{\hbox{\rm\  mod\ }}

\def\ds{\displaystyle}

\def\ul{\underline}
\def\id{\hbox{\rm id}}
\def\Im{\hbox{\rm Im\thinspace}}

\def\pr{\hbox{\rm pr}}

\def\Spec{\hbox{\rm Spec$\,$}}

\def\pr{\hbox{\rm p}}

\def\Grass{\hbox{\rm Grass}}

\def\defop#1{\expandafter\def\csname #1\endcsname  
{\mathop{\sans #1}\nolimits}}                      

\font\tensans=cmss10
\font\sevensans=cmss8
\font\fivesans=cmss8 at 6.666pt
\newfam\sansfam 
 \textfont\sansfam=\tensans
 \scriptfont\sansfam=\sevensans
 \scriptscriptfont\sansfam=\fivesans
\def\sans{\fam\sansfam\tensans}

\defop{Hom}                          
\defop{End}
\defop{Mod}
\defop{Bimod}
\defop{Qcoh} 
\defop{Vect}
\defop{vect}
\defop{mod}  
\defop{Rep} 
\defop{rep}
\defop{Spec}
\defop{dim} 
\defop{Iso}
\defop{Ind}
\defop{Aut}
\defop{Max}
\defop{Irr}
\defop{trace}
\defop{det}
\defop{td}
\defop{ch}
\defop{Ext}
\defop{Syst}
\defop{pr}
\defop{Gal}
\defop{Ext}
\defop{p}
\defop{vol}
\defop{Ind}
\defop{ch}
\defop{grad}
\defop{ev}
\def\limproj{\mathop{\oalign{lim\cr\hidewidth$\longleftarrow$\hidewidth\cr}}}
\defop{rang}
\defop{tors}
\defop{Diff}
\defop{Nm}
\defop{Span}
\defop{Isom}

\catcode`@=11 
\def\wlog#1{} 

\newskip\hsssglue \hsssglue=0pt plus 1fill minus 1fill 

\newdimen\unitlength \newdimen\linethickness
\newdimen\@picheight \newdimen\@xdim \newdimen\@ydim \newdimen\@len \newdimen\@save
\newcount\@multicount \newcount\@xarg \newcount\@yarg
\newbox\@picbox \newbox\@mpbox

\font\tenln=line10     \font\tenlnw=linew10
\font\tencirc=lcircle10 \font\tencircw=lcirclew10

\def\thinlines{\let\linefont=\tenln \let\circlefont=\tencirc
  \linethickness=\fontdimen8\linefont}
\def\thicklines{\let\linefont=\tenlnw \let\circlefont=\tencircw
  \linethickness=\fontdimen8\linefont}
\thinlines

\def\beginpicture(#1,#2)(#3,#4){\@picheight=#2\unitlength
  \setbox\@picbox=\hbox to#1\unitlength\bgroup \let\line=\@line
    \kern-#3\unitlength \lower#4\unitlength\hbox\bgroup\ignorespaces}
\def\endpicture{\egroup\hss\egroup
  \ht\@picbox=\@picheight \dp\@picbox=\z@
  \leavevmode\box\@picbox}

\def\put(#1,#2)#3{\raise#2\unitlength\rlap{\kern#1\unitlength #3}\ignorespaces}

\def\multiput(#1,#2)(#3,#4)#5#6{\@multicount=#5\relax
 \@xdim=#1\unitlength \@ydim=#2\unitlength \setbox\@mpbox=\hbox{#6}%
 \loop\ifnum\@multicount>0
   \raise\@ydim\rlap{\kern\@xdim \unhcopy\@mpbox}%
   \advance\@xdim#3\unitlength \advance\@ydim#4\unitlength
   \advance\@multicount\m@ne \repeat\ignorespaces}

\def\makebox(#1,#2)#3{\setbox\@picbox=\hbox to#1\unitlength{\hss#3\hss}%
  \@ydim=\ht\@picbox \advance\@ydim-\dp\@picbox
  \ht\@picbox=#2\unitlength \dp\@picbox=\z@
  \leavevmode\lower.5\@ydim\box\@picbox}

\newif\ifneg
\def\@line(#1,#2)#3{\@xarg=#1 \@yarg=#2 \@len=#3\unitlength \leavevmode
 \ifnum\@xarg<0 \reverseline \else \negfalse \@ydim=\z@\fi
 \ifnum\@xarg=0 \@vline
 \else\ifnum\@yarg=0 \@hline \else\@sline\fi\fi
 \ifneg\kern-\@len\else\@save=\@ydim\fi}
\def\reverseline{\negtrue \kern-\@len \@xarg=-\@xarg
 \@ydim=\@len \multiply\@ydim\@yarg \divide\@ydim\@xarg \@yarg=-\@yarg}

\def\@hline{\vrule height.5\linethickness depth.5\linethickness width\@len}
\def\@vline{\kern-.5\linethickness\vrule width\linethickness
  \ifnum\@yarg<0 height\z@ depth\else depth\z@ height\fi\@len
  \kern-.5\linethickness}

\def\@sline{\setbox\@picbox=\hbox{\linefont \count@=\@xarg \multiply\count@ 8
 \ifnum\@yarg>0 \advance\count@\@yarg \advance\count@-9
 \else \advance\count@-\@yarg \advance\count@ 55 \fi \char\count@}%
 \ifnum\@yarg<0 \@picheight=-\ht\@picbox \advance\@ydim\@picheight
 \else \@picheight=\ht\@picbox \fi
 \@xdim=\wd\@picbox \@save=\@ydim
 \loop\ifdim\@xdim<\@len \raise\@ydim\copy\@picbox
  \advance\@xdim\wd\@picbox \advance\@ydim\@picheight \repeat
 \advance\@xdim-\@len \kern-\@xdim
 \multiply\@xdim\@yarg \divide\@xdim\@xarg \advance\@ydim-\@xdim
 \raise\@ydim\box\@picbox}

\def\vector(#1,#2)#3{\@line(#1,#2){#3}%
 \ifnum\@xarg=0 \@vvector \else\ifnum\@yarg=0 \@hvector \else\@svector\fi\fi}
\def\@hvector{\ifneg\rlap{\linefont\char27}\else
 \smash{\llap{\linefont\char45}}\fi} 
\def\@vvector{\ifnum\@yarg<0 \raise-\@len\rlap{\linefont\char63}%
 \else\setbox\@picbox=\rlap{\linefont\char54}\advance\@len-\ht\@picbox
 \raise\@len\box\@picbox\fi}

\def\@svector{\setbox\@picbox=\hbox to\z@{\linefont
 \ifnum\@yarg<0 \count@=55 \@yarg=-\@yarg \else\count@=-9 \fi
 \ifneg\multiply\@xarg16 \multiply\@yarg2
 \else\hss 
  \ifnum\@xarg>2 \multiply\@xarg9 \multiply\@yarg2 \advance\count@29
  \else\ifnum\@yarg>2 \multiply\@xarg16 \multiply\@yarg9 \advance\count@-20
   \else\multiply\@xarg24 \multiply\@yarg3 \fi\fi\fi
  \advance\count@\@xarg \advance\count@\@yarg \char\count@
  \ifneg\hss\fi}
 \raise\@save\box\@picbox}

\def\disk#1{\@len=#1\unitlength \count@='160 \@diskcirc}
\def\circle#1{\@len=#1\unitlength \count@='140 \@diskcirc}
\def\@diskcirc{\setbox\@picbox=\hbox{\circlefont\char\count@}\@xdim=\wd\@picbox
 \leavevmode \ifdim\@len>15.499\@xdim \@bigdc \else \@smalldc\fi}
\def\@bigdc{\ifnum\count@<'160 \@bigcirc
 \else \@len=15\@xdim \@diskcirc\fi}
\def\@smalldc{{\advance\@len-.5\@xdim
 \loop\ifdim\@xdim<\@len \advance\count@\@ne \advance\@xdim\wd\@picbox\repeat
 \hbox{\circlefont\char\count@}}}
\def\@bigcirc{{\circlefont\count@=15
 \setbox\@picbox=\hbox{\char\count@}\@xdim=\wd\@picbox
 \ifdim\@len>2.5\@xdim \@len=2.5\@xdim\fi
 \advance\@len-.125\wd\@picbox
 \loop\ifdim\@xdim<\@len \advance\count@ 4 \advance\@xdim.25\wd\@picbox\repeat
 \@ydim=.5\@xdim \advance\@ydim.5\linethickness
 \setbox\@picbox=\vbox{\hbox{\char\count@\advance\count@-3\char\count@}%
  \nointerlineskip
  \hbox{\advance\count@\m@ne\char\count@\advance\count@\m@ne\char\count@}}%
 \kern-\@ydim\lower\@ydim\box\@picbox}}

\newif\ifovaltl \newif\ifovaltr \newif\ifovalbl \newif\ifovalbr
\ovaltltrue \ovaltrtrue \ovalbltrue \ovalbrtrue
\def\oval(#1,#2){\@xdim=#1\unitlength \@ydim=#2\unitlength
 {\circlefont \setbox\@picbox=\hbox{\char0}
 \ifdim\@xdim<\wd\@picbox \@xdim=\wd\@picbox\fi
 \ifdim\@ydim<\wd\@picbox \@ydim=\wd\@picbox\fi
 \@save=\@xdim \ifdim\@ydim<\@save \@save=\@ydim \fi
 \count@=39
 \loop \setbox\@picbox=\hbox{\char\count@}\ifdim\@save<\wd\@picbox
  \advance\count@-4 \repeat
 \setbox\strutbox=\hbox{\vrule height\ht\@picbox depth\dp\@picbox width\z@
   \kern\wd\@picbox}%
 \@save=.5\wd\@picbox \advance\@save-.5\linethickness
 \setbox0=\hbox to\@xdim{\ifovaltl\char\count@\else\strut\fi
  \kern-\@save\leaders\hrule height\ifovaltl\linethickness\else\z@\fi\hfil
  \leaders\hrule height\ifovaltr\linethickness\else\z@\fi\hfil\kern\@save
  \ifovaltr\advance\count@-3\char\count@\else\strut\fi\kern-\wd\@picbox}%
  \advance\count@\m@ne
 \setbox2=\hbox to\@xdim{\ifovalbl\char\count@\else\strut\fi
  \kern-\@save\leaders\hrule height\ifovalbl\linethickness\else\z@\fi\hfil
  \leaders\hrule height\ifovalbr\linethickness\else\z@\fi\hfil\kern\@save
  \ifovalbr\advance\count@\m@ne\char\count@\else\strut\fi\kern-\wd\@picbox}%
 \@save=\@ydim \advance\@save-\wd\@picbox \divide\@save 2
 \setbox\@picbox=\vbox{\box0\nointerlineskip
  \hbox to\@xdim{\vrule height\@save width\ifovaltl\linethickness\else\z@\fi
    \hfil\ifovaltr\vrule width\linethickness\kern-\linethickness\fi}%
  \nointerlineskip
  \hbox to\@xdim{\vrule height\@save width\ifovalbl\linethickness\else\z@\fi
    \hfil\ifovalbr\vrule width\linethickness\kern-\linethickness\fi}%
  \nointerlineskip\box2}%
  \@save=.5\@ydim \advance\@save.5\linethickness \leavevmode
  \kern-.5\@xdim \kern-.5\linethickness \lower\@save\box\@picbox}}

\def\cpic#1\endcpic{\vcenter{\hbox{\beginpicture#1\endpicture}}}

\newdimen\@xi \newdimen\@xii \newdimen\@xiii \newdimen\@xiv
\newdimen\@xpt \newdimen\@xoldpt
\newdimen\@yi \newdimen\@yii \newdimen\@yiii \newdimen\@yiv
\newdimen\@ypt \newdimen\@yoldpt
\def\squine(#1,#2,#3,#4,#5,#6){\setbox\@picbox\hbox{\tencirc q}%
 \global\@xoldpt=#1\unitlength \global\@yoldpt=#4\unitlength \kern\@xoldpt
 \@xi=\@xoldpt \@xii=#2\unitlength \@xiii=#3\unitlength
 \@yi=\@yoldpt \@yii=#5\unitlength \@yiii=#6\unitlength
 \squinerec
 \@xpt=#3\unitlength \@ypt=#6\unitlength \@addpoint
 \raise\@ypt\copy\@picbox}
\newif\iffar
\def\squinerec{\farfalse \testnear\@xi\@xiii \testnear\@yi\@yiii
 \iffar \decast \fi}
\def\testnear#1#2{\@save=#1\advance\@save-#2%
 \ifdim\@save<\z@ \@save=-\@save\fi \ifdim\@save>\p@ \fartrue \fi}
\def\decast{\@xpt=\@xi \advance\@xpt\@xii \divide\@xpt2
 \advance\@xii\@xiii \divide\@xii2
 \@xiv=\@xpt \advance\@xiv\@xii \divide\@xiv2
 \@ypt=\@yi \advance\@ypt\@yii \divide\@ypt2
 \advance\@yii\@yiii \divide\@yii2
 \@yiv=\@ypt \advance\@yiv\@yii \divide\@yiv2
 \begingroup\@xii=\@xpt \@xiii=\@xiv
  \@yii=\@ypt \@yiii=\@yiv \squinerec\endgroup
 \@xpt=\@xiv \@ypt=\@yiv \@addpoint
 \@xi=\@xiv \@yi=\@yiv \squinerec}
\def\@addpoint{
 \global\advance\@xoldpt-\@xpt \wd\@picbox=-\@xoldpt
 \raise\@yoldpt\copy\@picbox \global\@xoldpt=\@xpt \global\@yoldpt=\@ypt}

\catcode`\@=12 

\begin{document}

\centerline{\capit Sur le morphisme de Barth}
\centerline{Joseph Le Potier et Alexander Tikhomirov}
\bigskip
\auteurcourant{Joseph Le Potier et Alexander Tikhomirov}
\titrecourant{Sur le morphisme de Barth}
\section {Introduction}
Soit $n$  un entier $\geq 2.$ 
Etant donn\'e un faisceau semi-stable  $F$  de rang 2 et de classes de
Chern $c_1=0, c_{2}=n$  sur le plan projectif complexe, la cohomologie 
$H^{q}(F(i))$ est  nulle pour $q=0$
  et $i\leq 0,$  ou pour $q=2$  et $i\geq -2.$  Il  r\'esulte du th\'eor\`eme de
Riemann-Roch que 
$\dim H^{1}(F(-1))=\dim H^1(F(-2))=n,$  et $\dim H^{1}(F)=n-2.$
On dit qu'une droite
$\ell\subset \proj_{2}$  est de saut pour $F $ si la
restriction $F|_{\ell}$  n'est pas isomorphe au fibr\'e trivial de rang 2
sur
$\ell$.  Cohomologiquement ceci signifie  aussi
$H^1(F(-1)|_{\ell})\not=0$  de sorte que l'ensemble des droites  de saut de
$F$  est le support du   conoyau  du morphisme canonique de faisceaux
localement libres de rang
$n$  sur le plan projectif dual
$$H^{1}(F(-2)) \otimes {\cal O}_{\proj^*_2}(-1)\hfl{}{}
H^{1}(F(-1))\otimes {\cal O}_{\proj_2^*}.$$ 
D'apr\`es le th\'eor\`eme de Grauert et M\"ulich [2], ce support est de dimension
1: autrement dit,  le d\'eterminant du morphisme ci-dessus d\'efinit une
courbe $\beta_{F}$ de degr\'e $n$ dans le plan projectif dual qu'on 
appelle la courbe des droites de saut de
$F.$

 On d\'esigne par
$M_{n}$  l'espace de modules des classes de S-\'equivalence
(\cf \S 2) de faisceaux semi-stables de rang 2  et de classes de Chern
$(0,n)$ sur $\proj_2.$ C'est une vari\'et\'e projective irr\'eductible et
normale de dimension
$4n-3$, localement factorielle [5].  Les points qui repr\'esentent des
faisceaux singuliers constituent une hypersurface
$\partial M_{n}.$ On consid\`ere d'autre part le syst\`eme lin\'eaire
complet
${\cal C}_{n}$  des courbes de degr\'e $n$ dans le plan
projectif dual: c'est un espace projectif de dimension
${1\over 2}n(n+3)$. On dispose
d'un morphisme 
$\beta:M_{n}\hfl{}{} {\cal C}_{n} $
qui associe \`a la classe d'\'equivalence d'un faisceau semi-stable $F$ la
courbe
$\beta_{F}$ de ses droites de saut: c'est le morphisme de Barth. On
sait que ce morphisme est un isomorphisme pour $n=2,$  et que pour
$n=3$, c'est un morphisme de degr\'e 3. Pour $n\geq 2
,$  on sait par
ailleurs [13] que la restriction de
$\beta$ \`a l'ouvert des faisceaux localement libres (obligatoirement
$\mu-$stables)  est  quasi-fini.  La
courbe des droites de saut d'un faisceau singulier est obligatoirement
r\'eductible; ainsi, le morphisme de Barth est un morphisme fini
au-dessus de l'ouvert des courbes irr\'eductibles. On sait d'autre part
 d'apr\`es Barth [3]
que l'image de $\beta$ rencontre l'ouvert des courbes lisses.

 Le but principal  de cet article la d\'emonstration du
th\'eor\`eme suivant: 

\begin{th} \label{main}
Pour $n\geq 4,$ le morphisme de Barth est g\'en\'eriquement injectif. 
\end{th}
 Ceci signifie qu'il existe un ouvert de Zariski non vide $U$ de $M_{n}$ 
tel que la restriction de $\beta$ \`a $U$ soit injective.
On sait par ailleurs  que l'image r\'eciproque du fibr\'e ${\cal O}(1)$  est
le fibr\'e d\'eterminant ${\cal D}$ de Donaldson.  Les nombres de Donaldson du
plan projectif sont d\'efinis par
$$q_{4n-3}= \int_{
M_{n}}c_{1}({\cal D})^{4n-3}$$
et sont  maintenant bien connus [6], [7], [18]. On a par exemple
$$q_{13}=54, q_{17}= 2540, q_{21}= 233208.$$

\begin{cor}
L'image de morphisme de Barth est une sous-vari\'et\'e ferm\'ee de ${\cal C}_{n}$
de dimension $4n-3$  et de degr\'e 
$q_{4n-3}.$  
\end{cor}
La vari\'et\'e ${\cal L}=\beta(M_{4})$  est  donc  une hypersurface
irr\'eductible  de
${\cal C}_{4}=\proj_{14}$, de degr\'e $54$,  constitu\'ees de quartiques
dites de L\"uroth:  les courbes  g\'en\'eriques   de ${\cal L}$ sont les
quartiques lisses  qui sont circonscrites \`a un pentagone complet (\cf
\S 5.4). De m\^eme, les courbes g\'en\'eriques de
$\beta(M_{5})$  sont  les quintiques lisses circonscrites \`a un hexagone
complet: ce sont les quintiques dites de Darboux; la vari\'et\'e 
$\beta(M_{5})$  est une sous-vari\'et\'e ferm\'ee de codimension 3 de ${\cal
C}_{5}=\proj_{20}$  de degr\'e
$2540.$  
\bigskip
Le th\'eor\`eme \footnote{$(^1)$}{Le th\'eor\`eme 1.1 est vrai en fait  si 
le corps de base $k$ est  alg\'ebriquement clos et de caract\'eristique
0,  ce qu'on supposera d\'esormais.
} se d\'emontre par  r\'ecurrence sur $n,
$ \`a partir du  cas  $n=4,$ qui n'est nullement \'evident. Pour
$n=4$,  il s'obtient en mettant  en \'evidence   un diviseur
${\cal D}_{2}\subset M_{n}$  satisfaisant aux conditions
suivantes, sur un ouvert non vide ${\cal D}_{2}^*$ de ${\cal D}_{2}$:

--- les fibres de
$\beta$  qui   rencontrent  
${\cal D}^*_{2}$  sont r\'eduites \`a un point. 

--- le morphisme $\beta $ est non ramifi\'e  en tout point  
de ${\cal D}^*_{2}$.

\noindent En particulier, la
restriction de $\beta$ \`a cet ouvert ${\cal D}_{2}^*$ est injective. La
construction d'un tel diviseur r\'esulte de l'\'etude de l'image r\'eciproque
$\beta^{-1}({\cal S}_{4})$ de l'hypersurface ${\cal S}_{4}\subset
{\cal C}_{4}$ des quartiques singuli\`eres;  cette hypersurface  a trois
composantes  irr\'eductibles dont l'une est le bord $\partial M_{n}$ 
et  deux autres composantes irr\'eductibles not\'ee
${\cal D}_{1}$ et ${\cal D}_{2}.$  La composante ${\cal D}_{2}$ appara\^{\i}t en
fait avec multiplicit\'e 2, et le bord $\partial M_{4}$  avec
multiplicit\'e 3 (\cf th\'eor\`eme \ref{multiplicites}).  Les points du
diviseur ${\cal D}_{2}$  sont les classes des faisceaux semi-stables $F$
de rang 2, de classes de Chern
$(0,4)$ qui ont au moins une droite de saut $\ell$ d'ordre
$\geq 2,$ c'est-\`a-dire telle que $\dim H^{1}(F(-1)|_{\ell})\geq 2.$ Une
description commode des composantes
${\cal D}_{1}$  et
${\cal D}_{2}$  s'obtient \`a partir de la description de $M_{4}$ comme
contraction  d'un espace de modules $S_{4}$  de syst\`emes coh\'erents
semi-stables
$(\Gamma,\Theta)$ o\`u $\Theta$  est faisceau coh\'erent pur de dimension 1 
de caract\'eristique d'Euler-Poincar\'e $\chi=6$, dont l'id\'eal de Fitting
d\'efinit une conique
$\hbox{\goth{c}}$ (appel\'ee support sch\'ematique de $\Theta$) et
$\Gamma$ un sous-espace vectoriel de dimension 2 de $H^0(\Theta)$ (\cf
th\'eor\`eme \ref{eclate}).
 Un point  g\'en\'eral de ${\cal D}_{2}$   provient  d'un syst\`eme
lin\'eaire  $(\Gamma,\Theta)$ o\`u   $\Theta$  est un faisceau inversible 
de degr\'e total 5,  semi-stable sur une conique singuli\`ere
$\hbox{\goth{c}}$  de
$\proj_2.$  Nous avons en fait  g\'en\'eralis\'e cet \'enonc\'e au sous-sch\'ema
$P_{n}\subset M_{n}$  des faisceaux de Poncelet, dont peut d\'ecrire une
modification en introduisant l'espace de modules $S_{n}$  des syst\`emes
coh\'erents semi-stables $(\Gamma,\Theta)$ o\`u
$\Theta$  est un faisceau coh\'erent pur de dimension un dont le support
sch\'ematique est une conique, de caract\'eristique d'Euler-Poincar\'e 
$\chi=n+2,$ et
$\Gamma$  un sous-espace vectoriel de dimension 2 de $H^0(\Theta).$ 
Les sch\'emas $S_{n}$  et $P_{n}$  sont des vari\'et\'es irr\'eductibles et
normales de dimension $2n+5$ et on dispose d'un morphisme birationnel
$S_{n}\hfl{}{} P_{n}.$  Dans le cas $n=4$, on a $P_{4}=M_{4}$ et ce
morphisme est l'\'eclatement de la sous-vari\'et\'e des faisceaux sp\'eciaux. 
D'autres propri\'et\'es  g\'eom\'etriques  de $S_{n}$ et $P_{n}$ 
sont d\'ecrites dans les propositions
\ref{S_{n}}  et
\ref{P_{n}}. 

\bigskip
Pour $n\geq 5$, on \'etudie  le comportement des fibres de  $\beta$  au
voisinage d'un point suffisamment g\'en\'eral de $B=\beta(\partial M_{n}).$
Les  points  de $\partial M_{n}$  repr\'esentent des faisceaux singuliers;
ceux qui n'ont qu'un seul point singulier s'ins\`erent dans une suite exacte
$$0\hfl{}{} F\hfl{}{} E\hfl{}{}k(x)\hfl{}{}0$$  o\`u $x$  est le point
singulier de $F$, $k(x)$ le faisceau structural du  point
$x,$ et
$E$  le bidual de
$F;$ ce faisceau localement libre $E$ est stable de classes de Chern
$c_{1}=0, c_{2}=n-1$. La courbe des droites de saut de
$F$  est alors $\beta_{F}=\beta_{E}+\check{x},$ o\`u $\check{x}$ 
est la droite du plan projectif dual    $\proj_2^*$ d\'efinie par le
pinceau de droites de $\proj_2$ passant par $x.$
 Ainsi, un point g\'en\'eral de 
$Y$ repr\'esente une courbe de degr\'e $n$ qui se d\'ecompose en une
droite  et  la courbe, de degr\'e $n-1,$ des droites de
saut de $E.$ Un r\'esultat de Str\o mme 
permet d'affirmer qu'au-dessus  d'un  point g\'en\'eral $\hbox{\goth{q}}$
de
$B$,  on ne trouve que des classes de faisceaux $F$ singuliers, et ceci 
en un seul point;  l'hypoth\`ese de r\'ecurrence permet de dire   que le
faisceau localement libre $E$ et le point singulier  sont  parfaitement
d\'etermin\'es \`a isomorphisme pr\`es par la courbe $\hbox{\goth{q}}$;  la
fibre 
$\beta^{-1}(\hbox{\goth{q}})$ est  alors  isomorphe \`a   la droite
projective
$\proj_1=\proj(E_{x})$ ainsi associ\'ee \`a  $\hbox{\goth{q}}$. 
  On consid\`ere alors  la factorisation de Stein
$$\diagram{M_{n}&\hfl{\widetilde{\beta}}{}& \widetilde{{\cal C}}_{n}\cr
&\sefl{\beta}{}&\vfl{}{\nu}\cr
&&{\cal C}_{n}
\cr}$$
Le morphisme $\widetilde{\beta}$  est   birationnel; 
au-dessus d'un point g\'en\'eral $\hbox{\goth{q}}\in Y$,  la fibre de
$\nu$  est r\'eduite \`a un point. On est ramen\'e \`a prouver  qu'au-dessus
d'un ouvert non vide de $Y$, le morphisme $\nu$  est non ramifi\'e:  c'est
la principale difficult\'e  de la d\'emonstration. Nous traduirons en fait
directement sur $\beta$ comment interpr\'eter  cette condition de non
ramification (\cf proposition \ref{cle2}).
\bigskip
Les faisceaux  de Poncelet  
ainsi que les faisceaux de Hulsbergen, jouent un grand r\^ole dans l'\'etude
pr\'esent\'ee ici.  Un  faisceau semi-stable  $F$ de rang
2 et classes de Chern $(0,n)$  sur le plan projectif est dit de Poncelet
si
$F(1)$ a  au moins deux sections lin\'eairement ind\'ependantes. Comme nous
l'avons d\'ej\`a dit,  les faisceaux de Poncelet constituent un sous-sch\'ema
ferm\'e irr\'eductible
$P_{n}$ de
$M_{n}$  de dimension $2n+5$ si $n\geq 4$(\cf proposition \ref{P_{n}}); 
pour
$n=4,$ tous les faisceaux semi-stables  sont de Poncelet. M. Toma a 
d\'emontr\'e dans [24]   que pour
$n\geq 5$ la restriction du morphisme de Barth \`a ce ferm\'e $P_{n}$  est
encore g\'en\'eriquement injective \footnote{$(^2)$}{Curieusement, la
d\'emonstration de M. Toma ne marche pas  pour $n=4$.}. L'image de
$P_{n}$  est la vari\'et\'e des courbes dites de Poncelet de degr\'e $n$; les
courbes lisses 
  qui passent par les sommets d'un polyg\^one \`a n+1  c\^ot\'es
circonscrit  \`a une conique lisse forment un ensemble dense dans le
ferm\'e $\beta(P_{n})$ des courbes de Poncelet.
Le calcul du degr\'e de
cette sous-vari\'et\'e ferm\'ee irr\'eductible  $\beta(P_{n})$ de dimension
$2n+5$ de ${\cal C}_{n}$ se ram\`ene encore \`a un calcul de nombres
d'intersection sur la vari\'et\'e
$M_{n}.$  Par exemple, pour $n=5,$  il existe un faisceau universel
${\bf F}$ sur
$M_{5}\times
\proj_{2}$, et le degr\'e de $\beta(P_{5})$  est donn\'e par 
$$\int_{M_{5}}c_{2}(\pr_{1!}({\bf F}(1)))c_{1}({\cal D})^{15}.$$

Les faisceaux dits de Hulsbergen sont les faisceaux semi-stables $F$ de
rang 2 et de classe de Chern $(0,n)$  tels que $F(1)$  aient au moins une
section non nulle. Si $n\leq 5$, c'est le cas de tous les faisceaux
semi-stables; pour $n\geq 6$ ils constituent (\cf [17], th\'eor\`eme 4.7) un 
ferm\'e irr\'eductible
$H_{n}$ de
$M_{n}$  de dimension
$3n+2.$  L'image de $H_{n}$ par le morphisme de Barth $\beta$ est
encore un ferm\'e irr\'eductible de dimension $3n+2$ qui contient comme
ensemble dense les points repr\'esentant 
 les courbes lisses  de degr\'e  $n$ passant par les sommets d'un
polyg\^one  quelconque \`a $n+1$ c\^ot\'es; ces courbes sont appel\'ees  courbes de
Darboux de degr\'e
$n$.  Une question int\'eressante qui \`a notre connaissance
n'est pas r\'esolue   est de d\'eterminer si la restriction de $\beta$ au
ferm\'e
$H_{n}$  est encore g\'en\'eriquement injective. Pour $n=5,$ on a
$H_{5}=M_{5}$  et  ceci est donc vrai. Si c'\'etait encore  le cas pour 
$n\geq 6$, le degr\'e du ferm\'e des courbes de Darboux pourrait encore se
calculer comme nombre d'intersection sur l'espace de modules
$M_{n}.$ 
\bigskip

Voici le plan de cet article:
les sections 2 et 3 sont consacr\'ees \`a des rappels sur la notion de
stabilit\'e et semi-stabilit\'e des faisceaux alg\'ebriques coh\'erents et des
syst\`emes coh\'erents. Dans  la section 4, nous \'etudions les faisceaux de
Poncelet  et les singularit\'es des courbes de Poncelet de degr\'e
$n$ associ\'ees  \`a ces faisceaux de Poncelet dans la section 5: l'\'enonc\'e
obtenu \'etend aux syst\`emes coh\'erents consid\'er\'es dont le support
sch\'ematique peut \^etre une conique singuli\`ere  un r\'esultat bien connu
(\cf Maruyama [20], Trautmann [25], Vall\`es [27]) pour les syst\`emes
lin\'eaires sur les coniques lisses. 
 Ceci permet  l'\'etude dans la section 6 du diviseur 
de $M_{4}$ des faisceaux semi-stables dont la courbe des droites de saut
associ\'ee est singuli\`ere. Les sections 7  et 8   sont la cl\'e de la
d\'emonstration dans le cas $n=4$: on associe \`a une quartique singuli\`ere
g\'en\'erale
$\hbox{\goth{q}}$ l'unique  conique $\hbox{\goth{p}}$ qui passe  par les
6 points de contact des tangentes issues du point singulier; \'ecrire que
cette conique est singuli\`ere  revient
 \`a \'ecrire que la quartique $\hbox{\goth{q}}$ est la courbe des droites de
saut d'un fibr\'e de Poncelet associ\'e \`a un pinceau de diviseurs de degr\'e
5  sur une conique singuli\`ere; de plus, cette conique singuli\`ere et le
pinceau sont d\'etermin\'es de fa\c{c}on unique \`a partir de la quartique.
 Cette pr\'esentation
simplifie beaucoup l'esquisse que nous avions  pr\'esent\'ee dans  [16]. 
De plus le point fondamental de la d\'emonstration, \`a savoir le calcul
du rang de
$\beta$  au point g\'en\'erique du diviseur ${\cal D}_{2}$  n'y \'etait pas
abord\'e.  C'est l'\'etude des \'equations  du ferm\'e $\beta({\cal D}_{2})$
qui permet de calculer ce rang
 dans la section 9:  on montre qu'au-dessus
d'une telle quartique de L\"uroth, le morphisme de Barth est non ramifi\'e;
ceci s'obtient par un calcul de d\'eriv\'ee, en d\'eformant la conique et le
pinceau. Dans les sections 10 \`a 14 on d\'emontre, par r\'ecurrence sur $n,$ 
le th\'eor\`eme principal.

\section{L'espace de modules $M_{X}(c)$}
On rappelle ici les notions de stabilit\'e et semi-stabilit\'e [22]. 
Soit $(X,{\cal
O}{_X}(1))$ une vari\'et\'e projective lisse
 polaris\'ee de dimension $n.$ On d\'esigne par $K(X)$ le groupe de
Grothendieck des classes de ${\cal O}{_X}$-modules coh\'erents ; il est
\'equip\'e de la forme quadratique enti\`ere $q$ d\'efinie par $ q(u) =
\chi(u^2) $, dont la forme polaire associ\'ee est not\'ee $ \langle \ ,\
\rangle$. Si
$F$ est un faisceau alg\'ebrique coh\'erent sur $X$, on d\'esigne par $P_F$
le polyn\^ome de Hilbert de $F$. Si $Y$ est une section hyperplane, on
d\'esigne par  $h$ la classe  de
${\cal O}_Y$ , de sorte que $\langle F,h^i \rangle$ est la
caract\'eristique d'Euler-Poincar\'e de la restriction de $F$ \`a une
intersection compl\`ete g\'en\'erale de $i$ sections hyperplanes.  Le
polyn\^ome de Hilbert de $F$ est donn\'e par la formule $$ P_{F}(m) =
\langle F,(1-h)^{-m}\rangle =
\sum_{n\geq i\geq 0} C_{m+i-1}^i \langle F,h^i\rangle $$ 
Ces d\'efinitions  s'\'etendent aux classes de Grothendieck $c \in K(X).$ Le
noyau de la forme quadratique $q$  est constitu\'e des classes $c \in K(X)$ 
dont le rang et les classes de Chern, vues dans l'anneau d'\'equivalence
num\'erique  sont  nuls.  En particulier, le polyn\^ome de Hilbert  d'une
telle classe est  identiquement  nul.  On d\'esigne par
$K_{\rm num}(X)$  l'alg\`ebre quotient
$$K_{\rm num}(X)= K(X)/\ker q.$$
Sur $K_{\rm num}(X)$  la forme quadratique a encore un sens; de m\^eme, le
polyn\^ome de Hilbert
$ P_{c}$ d'une  classe
$c\in K_{\rm num}(X)$ a un sens. 
Si la vari\'et\'e  $X$  est une surface, 
et si $c\in K_{\rm num}(X)$ a pour  rang $r,$ classe de Chern $c_1$  et
caract\'eristique d'Euler-Poincar\'e
$\chi$, la forme quadratique et le polyn\^ome de Hilbert sont donn\'es par les
formules suivantes:
$$\eqalign{q(c)&= 2r\chi + c_{1}^2-r^2\chi({\cal O}_{X})\cr
P_{c}(m)&={1\over 2} \langle c,h^2 \rangle {m(m+1)}+ \langle c,h 
\rangle m +\chi.\cr}$$
Quand $X$  est le plan projectif, la classe $h$  est la classe du
faisceau structural d'une droite; on a
$K(X)=
\entr^3$ et la forme  quadratique  ci-dessus est unimodulaire.  Alors
$K_{\rm num}(X)= K(X).$

\begin{defi}  Soit $F$  un faisceau alg\'ebrique coh\'erent de dimension
$d$ \footnote{$(^3)$}{\rm La dimension d'un faisceau
alg\'ebrique coh\'erent est la dimension de son support.} sur
$X.$ On appelle multiplicit\'e de
$F$ le nombre entier
$r$ d\'efini par  $r=r(F) =
\langle F,h^d\rangle .$
\end{defi}

 Si le faisceau $ F$ est sans torsion, on a $r=
{\rm rang}(F).\hbox {\rm degr\'e }(X)$.
 Si $F$ est de dimension $n-1,$ la multiplicit\'e $r$ est
aussi le degr\'e de la vari\'et\'e de Fitting de $F$. 
Une classe  de Grothendieck $c$  de  $K(X)$ est dite de dimension $\leq d$ si
elle appartient \`a l'id\'eal
$F_{d}K(X)$ engendr\'e par les  classes des  faisceaux coh\'erents de dimension
$\leq d.$  Si $c$  est  de dimension $d,$  la multiplicit\'e de $c$  est le
nombre $\langle c,h^d\rangle.$
Une classe de Grothendieck $c\in K(X)$ est dite $\geq 0$ si  le
polyn\^ome de Hilbert
$$m\mapsto\langle c,(1-h)^{-m} \rangle$$  est $\geq 0;$  s'il
n'est pas nul, ceci signifie qu'il existe un entier $j\geq 0$  tel que
$\langle c, h^{i}\rangle =0$  pour
$i>j$ et $\langle c, h^{j}\rangle >0.$  Les classes $c$  de polyn\^ome de
Hilbert nul  constituent un  sous-groupe  $N$ de $K(X).$ La relation
ci-dessus d\'efinit  une relation d'ordre sur  le groupe ab\'elien quotient
$K(X)/N.$ On \'ecrit $c \equiv 0$  si $c$  est de polyn\^ome de Hilbert nul.
\begin{defi} (Puret\'e)
Un ${{\cal O}_X}$-module coh\'erent $F$ de dimension d est dit
pur s'il n'a pas de sous-module coh\'erent non nul de
dimension $<$ d.

\end{defi}
\begin{defi} (Semi-stabilit\'e)
Soit $F$ un ${{\cal O}_X}$-module coh\'erent de dimension d, de
multiplicit\'e $ r$ et classe de Grothendieck $c.$ On dit que $F$ est {\it
semi-stable} si

{\rm (i)} il est  pur de dimension $d$

{\rm (ii)}  pour tout sous-module coh\'erent non nul $F'\subset F$ , de 
classe de Grothendieck $c'$ et de multiplicit\'e $r'$, on a 
$${c' \over r'} \leq {c\over r}$$
\end{defi}
On d\'efinit de m\^eme la notion de {\it stabilit\'e } en demandant
l'in\'egalit\'e  stricte lorsque $F'$ est distinct de $F$.

Soit $c\in K_{\rm num}(X)$  une classe de
Grothendieck de   dimension $d$ et de multiplicit\'e $r$. On consid\`ere la
cat\'egorie
des faisceaux semi-stables $F$ de dimension $d,$ de
multiplicit\'e $r$  tels que 
$rc(F)\equiv r(F)c. $ C'est une cat\'egorie ab\'elienne,
artinienne et noeth\'erienne. Elle a donc des filtrations de
Jordan-H\"older. Deux faisceaux coh\'erents $F$ et $G$  sont dit
S-\'equivalents s'ils ont m\^eme classe de Grothendieck et si leurs
gradu\'es  de Jordan-H\"older sont isomorphes.
\begin{th} \label{module} Soit $c\in K_{\rm num}(X)$. 

{\rm (i)} Il existe 
un espace de modules grossier $M_X (c) $ pour les
faisceaux semi-stables F de dimension d, de classe de Grothendieck
$c.$

{\rm (ii)} L'espace de modules $M_{X}(c)$ est un sch\'ema projectif  dont
les points sont les classes de S-\'equivalence de faisceaux semi-stables. 
\end{th}
Ce r\'esultat, \'enonc\'e sous une forme voisine, est d\'emontr\'e avec cette
g\'en\'eralit\'e par C. Simpson [22]; il  \'etend l'\'enonc\'e bien connu
de Gieseker et Maruyama [19].

\section {Syst\`emes coh\'erents}
 Soit $(X,{{\cal O}_X}(1))$ une vari\'et\'e alg\'ebrique projective et
lisse, de dimension $n,$ polaris\'ee.
\begin{defi}  On appelle syst\`eme coh\'erent de dimension $d$ sur $X$ un
couple
$\Lambda=(\Gamma,F)$ form\'e d'un faisceau alg\'ebrique coh\'erent $F$ de
dimension $d$ sur $X$
 et d'un sous-espace vectoriel $\Gamma\subset H^0(F)$.
\end{defi}

Un morphisme de syst\`emes coh\'erents 
$\Lambda=(\Gamma,F)\hfl{}{}\Lambda'=(\Gamma',F')$  est la donn\'ee d'un
morphisme de faisceaux alg\'ebriques coh\'erents $f:F\hfl{}{} F'$   tel que
$f(\Gamma)\subset
\Gamma'$.  Les syst\`emes coh\'erents constituent une cat\'egorie additive. On
peut la plonger dans une cat\'egorie ab\'elienne en introduisant la notion
de syst\`eme alg\'ebrique [9]: un syst\`eme alg\'ebrique est un triplet
$\Lambda=(\Gamma,i,F)$  o\`u $F$  est un faisceau alg\'ebrique, $\Gamma$  un
espace vectoriel et $i:\Gamma\hfl{}{}H^0(F)$  une application lin\'eaire.
Un morphisme de syst\`emes alg\'ebriques $$\Lambda=(\Gamma,i,F) \hfl{}{}
\Lambda'=(\Gamma',i',F')$$  est un couple $(f,g)$  o\`u $f:\Gamma\hfl{}{}
\Gamma'$  est une application lin\'eaire et $g:F\hfl{}{}F'$  un morphisme
de faisceaux alg\'ebriques tels que le diagramme
$$\diagram{
\Gamma &\hfl{f}{}&\Gamma'\cr
\vfl{i}{}&&\vfl{}{i'}\cr
H^0(F)&\hfl{g}{}&H^0(F')\cr
}$$  soit commutatif.  Les syst\`emes alg\'ebriques
constituent une cat\'egorie ab\'elienne qui a suffisamment d'objets injectifs. 
En particulier, on peut d\'efinir, \'etant donn\'es deux syst\`emes coh\'erents
$\Lambda=(\Gamma,F)$ et
$\Lambda'=(\Gamma',F')$ les espaces vectoriels
$\Ext^{q}(\Lambda,\Lambda').$ Ce sont des espaces vectoriels de
dimension finie, et  on a une suite exacte longue
$$\eqalign{0\hfl{}{}
\Hom(\Lambda,\Lambda')\hfl{}{}
\Hom(F,F')&\hfl{}{}\Hom(\Gamma,H^0(F')/\Gamma')\hfl{}{}\cr
\Ext^{1}(\Lambda,\Lambda')&\hfl{}{}\Ext^{1}(F,F')\hfl{}{}
\Hom(\Gamma,H^{1}(F'))\hfl{}{}\dots\cr}$$
On associe \`a un  syst\`eme coh\'erent $\Lambda=(\Gamma,F)$ la classe de
Grothendieck  dans
$K(X)$
$$c(\Lambda)=\dim \Gamma + c(F)$$

\begin{defi}
 Un  syst\`eme coh\'erent $\Lambda= (\Gamma,F)$ de dimension $d$ est dit
semi-stable  si $F$  est pur de dimension $d$  et si pour tout
sous-module $F'\subset F$ non nul de multiplicit\'e $r'$, on a, en posant
$\Gamma'=\Gamma\cap H^0(F)$  et $\Lambda'=(\Gamma',F')$
$${c(\Lambda')\over {r'}}\leq {c(\Lambda)\over {r}}$$
\end{defi}
Ceci signifie donc que $\ds{{\dim \Gamma'}\over {r'}}\leq \ds{{\dim
\Gamma}\over {r}}$ et en cas d'\'egalit\'e $\ds{c(F')\over {r'}}\leq \ds{c(F)\over
{r}}$.

 Quand on fixe le polyn\^ome de Hilbert de $F$   on obtient une famille de
syst\`emes coh\'erents semi-stables
$(\Gamma,F)$
 qui est  limit\'ee. On peut \'evidemment d\'efinir la notion de
filtration de Jordan-H\"older pour de tels syt\`emes coh\'erents: ceci tient
encore  au  fait qu'une classe  $c\in K(X)$ de dimension $d$ et
de multiplicit\'e $r$ \'etant donn\'ee la cat\'egorie ${\cal S}(c)$ des
syst\`emes coh\'erents
$\Lambda=(\Gamma,F)$ de multiplicit\'e $r$ tels que $rc(\Lambda)\equiv r(F)c$
est une cat\'egorie ab\'elienne,  noeth\'erienne et artinienne.  Deux syst\`emes
coh\'erents  de dimension $d$ et multiplicit\'e $r$ 
 sont  S-\'equivalents  s'ils ont m\^eme classe de Grothendieck et  si leur  
gradu\'e de Jordan-H\"older sont isomorphes.  Soient
$k$ un entier, et $c \in K_{\rm num}(X)$; il existe alors  pour les
syst\`emes coh\'erents  
$\Lambda=(\Gamma,F)$ tels que $\dim\ \Gamma=k, c(F)=c$  un
espace de modules grossier
${\Syst}_{X}(c,k)$; c'est un sch\'ema projectif dont les
points sont les classes de 
S-\'equivalence de  syst\`emes coh\'erents semi-stables [15], [9].\bigskip
\goodbreak

\medskip
\section{Faisceaux de Poncelet}

Soit $n$  un entier $\geq 4.$ On consid\`ere  l'espace de modules $M_n$ des
classes de $S-$\'equivalence de faisceaux semi-stables de rang 2, de
classes de Chern $(0,n)$ sur le plan projectif.  C'est une vari\'et\'e
irr\'eductible de dimension
$4n-3$; elle est lisse si $n$  est impair; quand $n$ est pair,  elle est 
normale  et localement factorielle; son lieu singulier
est de dimension $2n$ et constitu\'e des classes d'\'equivalence de faisceaux
strictement semi-stables. Dans cet espace de modules, les classes qui se
repr\'esentent par un faisceau singulier constituent une hypersurface not\'ee
$\partial M_{n}.$

Soit $F$  un faisceau semi-stable de rang 2  et de classes de Chern
$(0,n)$. On a $H^{2}(F(j))=0$ pour $j\geq -3.$
 La suite des  nombres de Hodge
$j\mapsto \dim H^1(F(j))$ d'un tel faisceau semi-stable est
strictement d\'ecroissante  pour $j\geq 0$ jusqu'\`a ce qu'elle s'annule (\cf
[10], lemme 1.1);  pour
$1\leq j\leq n-2$  on a donc $\dim H^1(F(j))\leq  n-2-j,$ et pour $j\geq n-2$ on a
$H^{1}(F(j))=0$.  D'apr\`es la formule de Riemann-Roch, on a en outre
$$\dim H^0(F(j))-\dim H^{1}(F(j))=\chi(j):=(j+1)(j+2)-n.$$  En
particulier,
$\dim H^0(F(1))\leq 3.$

\begin{defi} Un faisceau semi-stable $F$ de rang 2, de classes de Chern
(0,n)  est appel\'e faisceau de Poncelet  s'il satisfait \`a l'une des
conditions \'equivalentes 

$\dim H^0(F(1))\geq 2$

$\dim H^{1}(F(1))\geq n-4.$

\end{defi}
Les classes d'\'equivalence de faisceaux de Poncelet constituent un 
 ferm\'e de $M_{n}$  sous-jacent 
\`a un sous-sch\'ema ferm\'e $P_{n}$  d\'efini de la mani\`ere suivante: 
quand
$n$  est impair, il existe un faisceau universel $\bf F$ sur $M_{n}\times
\proj_2$  param\'etr\'e par
$M_{n}$   et le sous-sch\'ema $P_{n}$  est d\'efini par l'id\'eal de Fitting
d'indice $n-5$ (\cf Lang [11] p. 483) du faisceau  $R^{1}\pr_{1*}({\bf
F}(1)).$  Quand $n$  est pair, il n'existe plus de faisceau universel.
Si $F$  est semi-stable de classe de Grothendieck $2-nh^2,$
le faisceau $F(\ell)$  est engendr\'e par ses sections et
$H^{q}(F(\ell))=0$  pour $q\geq 1 $ et $\ell\geq n-1.$ Consid\'erons   pour
$\ell\geq n-1$
 un espace vectoriel $H$  de dimension $\chi(\ell)$,  le
fibr\'e vectoriel
$E=H\otimes{\cal O}(-\ell) $  et le sch\'ema de Hilbert
$R= \Hilb(E,2-nh^2) $ des faisceaux coh\'erents quotients de $E$  de
classe de Grothendieck $2-nh^2.$ L'espace de modules
$M_{n}$ est le quotient de l'ouvert $R^{ss}$ des points semi-stables 
par l'action du groupe r\'eductif $SL(H).$ Il existe alors un faisceau
universel
${\bf F}$  param\'etr\'e par
$R^{ss}$  et le sch\'ema
$P_{n}$  est le quotient de Mumford  du sous-sch\'ema ferm\'e d\'efini par
l'id\'eal de Fitting d'indice $n-5$  du faisceau $R^{1}\pr_{1*}({\bf
F}(1)).$ 
\bigskip
\subsection{Construction des faisceaux de Poncelet}

Consid\'erons un syst\`eme coh\'erent semi-stable $\Lambda=(\Gamma,\Theta)$ 
de classe de Grothendieck $c=2h+nh^2$  et  tel que $\dim \Gamma=2.$ Le
faisceau coh\'erent  $\Theta$ est un faisceau pur de dimension 1, de
classe  de Grothendieck
$c=2h+4h^2$:   sa  multiplicit\'e est 2 et sa
caract\'eristique d'Euler-Poincar\'e est $\chi=n+2$; le support
sch\'ematique d'un tel faisceau  est une conique $\hbox{\goth{c}}$. Le
faisceau $\Theta$  peut alors \^etre vu comme un ${\cal
O}_{\hbox{\goth{c}}}$-module; les points au voisinage  desquels ce
faisceau n'est pas un ${\cal
O}_{\hbox{\goth{c}}}$-module libre sont dits points singuliers de
$\Theta.$ Un tel syst\`eme coh\'erent d\'efinit  un \'el\'ement de
$\Ext^{1}(\check{\Theta},{\cal O})\otimes
\Gamma^*=\Ext^{1}(\check{\Theta},\Gamma^*\otimes {\cal
O})$ o\`u  le faisceau
$\check{\Theta}$  est   le dual, d\'efini par
$$\check{\Theta}=\ul{\Ext}^{1}(\Theta,{\cal O})= \Theta^* \otimes_{{\cal
O}_{\hbox{\goth{c}}}} {\cal O}_{\hbox{\goth{c}}}(2).$$  On en  d\'eduit une
extension
$$0\hfl{}{} \Gamma^* \otimes {\cal O}\hfl{}{} F(1)\hfl{}{}
\check{\Theta}\hfl{}{}0 \eqno (4.1)\quad$$

\begin{defi}  Un faisceau semi-stable de rang 2 et de classes de Chern
$(0,n)$  est dit sp\'ecial  si $\dim H^0(F(1))=3.$
\end{defi}

\begin{th}  \label{modification}  On pose $S_{n}=\Syst(2h+nh^2,2).$
\item{\rm (i)}  Le faisceau  $F$ construit ci-dessus est un faisceau de
Poncelet de classes de Chern $(0,n)$.

\item{\rm (ii)}  Le morphisme  canonique
$$\pi:S_{n}\hfl{}{} P_{n}$$
 qui associe \`a la classe du syst\`eme coh\'erent semi-stable
$(\Gamma,\Theta)$  la classe du faisceau $F$  est surjectif; c'est  un
isomorphisme au-dessus de l'ouvert de $P_{n}$ des classes de faisceaux de
Poncelet non sp\'eciaux.
\end{th}
On obtient ainsi un diagramme  de vari\'et\'es projectives
$$ \diagram{
S_{n}&\hfl{\sigma}{}& {\cal C}_{2}=\proj_5\cr
\vfl{\pi}{}&& \cr P_{n} \cr}
$$ dans lequel le morphisme $\pi$  est le morphisme d\'ecrit ci-dessus; 
le morphisme 
$\sigma$   associe \`a un point 
de $\Syst(2h+nh^2,2)$ repr\'esent\'e par un syst\`eme coh\'erent
$(\Gamma,\Theta)$ la conique $\hbox{\goth{c}}$ d\'efinie par le support
sch\'ematique  de
$\Theta$.  On  v\'erifiera aux paragraphes 4.3  et 4.4 que $S_{n}$ sont
des sch\'emas int\`egres (et m\^eme normaux)  de dimension
$2n+5.$ Il  r\'esulte de l'\'enonc\'e ci-dessus que le morphisme $\pi$  est 
birationnel.  

\begin{proof} D'apr\`es le
corollaire 4.22 de [15], le syst\`eme coh\'erent
$(\Gamma,F(1))$  est semi-stable. Pour en d\'eduire que $F$  est
semi-stable, la d\'emonstration est identique \`a celle qui est donn\'ee dans
le lemme 6.6  de [15]. La suite exacte (4.1)  montre que le faisceau
$F$  est bien un faisceau de Poncelet.
  R\'eciproquement, soient  $F$  un faisceau de Poncelet, et $\Gamma$  un
sous-espace vectoriel de dimension 2 de  $H^0(F(1)).$  Alors la
d\'emonstration du lemme 6.6  de [15]  montre le syst\`eme coh\'erent
$(\Gamma,F(1))$  est semi-stable.
Soit $\check{\Theta}$  le conoyau du morphisme d'\'evaluation $\Gamma
\otimes {\cal O}\hfl{\ev}{} F(1),    $  et $\Theta=
\ul{\Ext}^1({\check{\Theta}},{\cal O})$.    Alors $\Theta$ est un
faisceau pur de dimension 1, de classe de Grothendieck 
$2h+nh^2 $ et l'inclusion canonique
$\Gamma^*\hookrightarrow \Ext^1(\check{\Theta},{\cal O})=H^{0}(\Theta)$
obtenue en appliquant le foncteur ${\Hom}(-,{\cal O})$ \`a la suite
exacte 
$$0\hfl{}{} \Gamma\otimes {\cal
O}\hfl{\ev}{}F(1)\hfl{}{}\check{\Theta}\hfl{}{}0$$ d\'efinit  un syst\`eme
coh\'erent semi-stable 
$(\Gamma^*,\Theta)$ de dimension 1  d'apr\`es le corollaire 4.22
de [15], et le faisceau de Poncelet associ\'e est  
isomorphe   \`a
$F.$  Pour obtenir un  isomorphisme au-dessus de l'ouvert des faisceaux
non sp\'eciaux, il suffit de travailler en familles et d'utiliser les
propri\'et\'es de modules grossiers.
\end{proof}

\bigskip
\subsection{Description des faisceaux sp\'eciaux.}
On d\'esigne par  $Q$  le fibr\'e quotient canonique de rang 2 sur $\proj_2.$
\begin{prop}  Soit $F$ un faisceau semi-stable  de rang 2  et de classes
de Chern $(0,n)$  sur le plan projectif.  Les 
conditions suivantes  sont \'equivalentes:

{\rm (i)} le faisceau $F$  est sp\'ecial;

{\rm (ii)} $\dim H^{1}(F(1))= n-3;$

{\rm (iii)}  pour tout $1\leq i\leq n-2$, on a $\dim H^{1}(F(i))=n-2-i;$

{\rm (iv)} le faisceau $F$  a une droite de saut d'ordre $n-1;$

{\rm (v)}  
il existe  une droite $\ell\subset \proj_2$ et une suite exacte 
$$0\hfl{}{}Q^*\hfl{}{}F\hfl{}{}{\cal O}_{\ell}(1-n)\hfl{}{}0$$

Dans ces conditions la droite de saut d'ordre $n-1$  est unique.
\end{prop}
\begin{proof}  Les assertions (i)  et (ii)  sont \'evidemment \'equivalentes
d'apr\`es le th\'eor\`eme de Riemann-Roch, et (iii)  implique trivialement (ii).

Montrons que (i) $\Rightarrow $ (iii). Soit un $F$  un faisceau
sp\'ecial; consid\'erons un syst\`eme coh\'erent semi-stable
$(\Gamma,\Theta)$  o\`u
$\Theta$  est de classe de Grothendieck $2h+nh^2$ de support une conique
$\hbox{\goth{c}}$ et
$\Gamma$  de dimension 2, et tel  que le faisceau semi-stable
associ\'e soit isomorphe \`a
$F.$  Le faisceau  $\check{\Theta},$ qui a pour
caract\'eristique d'Euler-Poincar\'e $\chi(\check{\Theta})=4-n$ est alors
instable car sinon, le faisceau structural ${\cal O}_{\hbox{\goth{c}}}$ 
\'etant stable  de caract\'eristique d'Euler-Poincar\'e $1,$ par
semi-stabilit\'e
$H^{0}(\check{\Theta})= 0$,  et on 
aurait alors $\dim H^{0}(F(1))=2,$ ce qui contredit le fait que $F$ est
sp\'ecial. Il en r\'esulte que 
$\Theta$   est un faisceau instable de dimension 1.
 Le support
de $\Theta$  est donc obligatoirement une conique $\hbox{\goth{c}}$ 
singuli\`ere; plus pr\'ecis\'ement,  la filtration de Harder-Narasimhan de
$\Theta$  montre qu'il   existe une suite exacte  
$$0\hfl{}{} \Theta'\hfl{}{}\Theta \hfl{}{} \Theta''\hfl{}{} 0$$
 o\`u $\Theta'$ et $\Theta''$ sont des faisceaux
inversibles 
sur des droites $\ell'$ et $\ell''$ , de degr\'es
respectifs $a'$  et $a''$   tels que $a'>a''.$  
On a bien s\^ur $a'+a''=n.$
 Le fait
que
$(\Gamma,\Theta)$  soit semi-stable  implique en outre que $a'\leq n-1$  ou
ce qui revient au m\^eme $a''\geq 1.$  Par dualit\'e on a alors une suite exacte
$$0\hfl{}{}\check{\Theta}''\hfl{}{}\check{\Theta}\hfl{}{}
\check{\Theta}'\hfl{}{}0$$
Puisque $\check{\Theta}''={\cal O}_{\ell''}(-a''+1)$, et que 
$\check{\Theta}'$  est  un faisceau inversible  de degr\'e $-a'+1<0,$  le
fait que $F$  soit sp\'ecial est \'equivalent \`a $a''=1$  ou encore \`a 
$a'=n-1$.  Si c'est le cas on a une surjection 
$F(1)\hfl{}{} {\cal O}_{\ell''}(-n+2)$ ce qui signifie que $\ell''$  est
une droite de saut d'ordre $n-1$; ceci d\'emontre (iii). 

 Montrons que (iii)$\Rightarrow $(iv).  Dire que $\ell$  est une droite de
saut d'ordre $n-1$  signifie que  l'on a un morphisme non nul
$F\hfl{}{}{\cal O}_{\ell}(1-n):$    ceci implique   
$H^1(F(n-3))\not=0$ . 
D'apr\`es le comportement de la suite $i\mapsto \dim H^{1}(F(i))$  d\'ecrite
en pr\'eambule  ceci implique que
$\dim H^1(F(i))=n-2-i$  pour $1\leq i\leq n-2,$ ce qui d\'emontre (iv). 

 Montrons
que (i) \'equivaut \`a (v).  On a vu que si (i)  est vrai, il existe une
droite $\ell$  et un  morphisme surjectif $F\hfl{}{} {\cal
O}_{\ell}(1-n).$ Le noyau $K$ d'un tel morphisme est un faisceau sans
torsion de rang 2 de classes de Chern $(-1,1).$  Montrons qu'il est
semi-stable:  soit $L$  un sous-faisceau maximal de rang 1  de $K.$
Puisque $L$  est aussi un sous-faisceau de $F$, on a $c_{1}(L)\leq 0.$  Il
s'agit de montrer en fait que $c_{1}(L)\leq -1.$  Mais au cas o\`u
$c_{1}(L)=0,$ on aurait \`a cause de la semi-stabilit\'e de $K$
$$c_{2}(L)\geq {n\over 2}$$ et
d'autre part le quotient $K/L$  serait sans torsion de rang 1  et
aurait pour classe de Chern $${{1-h+h^2}\over
{1+c_{2}(L)h^2}}=1-h+(1-c_{2}(L)h^2).$$
Ainsi ${n\over 2}\leq c_{2}(L)\leq 1$ et $n\leq 2$ ce qui est contraire \`a l'hypoth\`ese.
Ceci d\'emontre donc qu'un tel sous-faisceau n'existe pas. Puisque $Q^*$ 
est le seul faisceau semi-stable de classes de Chern $(-1,1)$ le
faisceau 
$K$  est isomorphe \`a $Q^*.$  Ceci d\'emontre (v).  Evidemment (v) entra\^{\i}ne
l'existence d'une droite de saut d'ordre $n-1.$
Qu'une telle droite de saut soit unique, r\'esulte du 
fait que $\Hom (Q^*,{\cal O}_{\ell}(1-n))=0$  pour $n>2.$
\end{proof}

Un faisceau  sp\'ecial $F$ est obligatoirement stable; il s'\'ecrit 
de mani\`ere unique comme extension
 $$0\hfl{}{}Q^*\hfl{}{} F\hfl{}{}{\cal O}_{\ell}(1-n)\hfl{}{}0.$$ 
ils constituent une sous-vari\'et\'e
$\Sigma_{n}$ lisse  et irr\'eductible  de dimension $2n+2$  qui ne
rencontre pas le lieu singulier de
$M_{n}.$

\subsection{Le sch\'ema $S_{n}=\Syst(2h+nh^2,2)$}

\begin{prop}  \label{S_{n}}  Soit $n$  un entier $\geq 4.$

{\rm (i)}  Le sch\'ema $S_{n}$    est 
irr\'eductible et normal de dimensions $2n+5.$  

{\rm (ii)} L'ouvert des points
stables de $S_{n}$   est localement intersection compl\`ete.  En
particulier, si
$n$  est impair, $S_{n}$ est  localement intersection compl\`ete.

\end{prop}

\begin{proof}   La d\'emonstration n\'ecessite l'introduction  du sch\'ema
de Hilbert qui permet la construction de $S_{n}$, \`a cause de la pr\'esence
de points semi-stables  quand $n$  est pair.  Elle  peut \^etre
simplifi\'ee  
  pour $n$  impair car dans ce cas tous  les points
semi-stables sont stables;  en outre les faisceaux stables $\Theta$   
sont  obligatoirement localement libres sur leur support 
[14].  
 Si
$(\Gamma,\Theta)$  est un syst\`eme coh\'erent semi-stable repr\'esentant un
point de
$S_{n}$,
$\Theta$  est engendr\'e par ses sections et
$H^{1}(\Theta)=0$. On  consid\`ere un espace vectoriel $H$  de dimension
$N=n+2. $ Soient $c=2h+nh^2,$  et $\Hilb(H\otimes {\cal O},c)$  le
sch\'ema de Hilbert des faisceaux quotient de classe de Grothendieck $c.$
D\'esignons par  $R_{n}^{ss}$ le sous-sch\'ema ouvert de
$$R_{n}:=\Grass(2,H)\times \Hilb(H\otimes {\cal O},c)$$
des couples $(\Gamma,\Theta)$ satisfaisant
aux conditions suivantes   

(a) $H^{1}(\Theta)=0;$

(b) l'application lin\'eaire $H\hfl{}{} H^{0}(\Theta)$  est inversible;

(c) le couple $(\Gamma,\Theta)$   d\'efinit un syst\`eme coh\'erent
semi-stable.

L'espace de modules $S_{n}$  est un bon quotient de $R_{n}^{ss}$  par
l'action naturelle du groupe $PGL(H),$  et l'action de $PGL(H)$  est
libre sur l'ouvert $R_{n}^{s}$ des points stables. D\'esignons pour
$(\Gamma,\Theta)\in R_{n}^{ss}$ par 
$K$  le noyau du morphisme
$H\otimes {\cal O}\hfl{}{}\Theta$.  Au voisinage d'un tel point, l'espace
de modules $R_{n}$  est  isomorphe (\cf [9]) au sch\'ema des z\'eros d'un
morphisme d'un germe de sch\'ema lisse d'espace tangent de Zariski 
$\Hom(\Gamma,H^0(\Theta)/\Gamma)) \times \Hom(K,\Theta)$  \`a valeurs dans
$\Ext^1(K,\Theta).$  Il  r\'esulte  du lemme de Krull que toutes les
composantes irr\'eductibles sont de dimension  $\geq 2n +\dim
\Hom(K,\Theta)-\dim\Ext^1(K,\Theta)= 2n+\chi(K,\Theta)=2n+4+N^2.$
Consid\'erons la projection canonique 
$$\rho: R_{n}^{ss}\hfl{}{} S_{n}\hfl{\sigma}{} {\cal C}_{2}=\proj_5$$ 
Au-dessus de l'ouvert des coniques lisses, ce morphisme est lisse de
dimension  relative $2n+N^2-1$, et la fibre au-dessus d'un point est
isomorphe \`a 
$\Grass(2,H)\times PGL(H)$  par cons\'equent l'image r\'eciproque $U$ de cet
ouvert  est irr\'eductible de dimension $2n+4+N^2.$
 On va d\'emontrer  que $R_{n}^{ss}$  est
un sch\'ema irr\'eductible localement intersection compl\`ete de dimension
$2n+4+N^2$, donc de Cohen-Macaulay. Il suffit  de v\'erifier  que le
compl\'ementaire de cet ouvert
$U$ est de dimension
$<2n+4+N^2.$ Ceci r\'esulte du  lemme suivant:

\begin{lemme} \label{fibres} Les   fibres du morphisme
$\rho: R_{n}^{ss}\hfl{}{}
\proj_{5}$  sont de dimension $\leq 2n+N^2-1$.
\end{lemme}

En fait, cet \'enonc\'e entra\^{\i}nant que $R_{n}^{ss}$  est localement
intersection compl\`ete et donc de Cohen-Macaulay, il r\'esulte   que ces
fibres sont de Cohen-Macaulay  de dimension $2n+N^2-1:$  en particulier,
les fibres  du morphisme $\rho$ sont pures de dimension $2n+N^2-1.$
\medskip
\begin{proof}
Supposons d'abord que la conique $\hbox{\goth{c}}$ soit  r\'eduite:
au-dessus de $\hbox{\goth{c}}$, il n'y a alors, \`a isomorphisme pr\`es,
qu'un  nombre fini
$\Theta_{1}, \dots ,\Theta_{\ell}$ 
de faisceau $\Theta,$ satisfaisant \`a la propri\'et\'e suivante: pour tout
quotient $\Theta\hfl{}{}\Theta''$ de multiplicit\'e 1, $\chi(\Theta'')\geq 1.$
 Chacun de ces $\Theta_{i}$  est engendr\'e par ses sections et 
on a $H^{1}(\Theta_{i})=0$,  de plus,  le groupe des automorphismes de
$\Theta_{i}$ est r\'eduit \`a
$k^*$, sauf si $\Theta_{i}$  n'est pas localement libre sur
$\hbox{\goth{c}}.$ On a donc  un morphisme  canonique 
$$f_{i}:W_{i}=\grass(2,H) \times \Isom(H,H^0(\Theta_{i}))\hfl{}{}
R_{n};$$ d\'esignons par $W_{i}^{ss}$  l'image r\'eciproque de l'ouvert
$R_{n}^{ss}$ des points semi-stables.  Les fibres de 
$f_{i}$  sont les trajectoires sous  l'action libre du groupe 
$\Aut(\Theta_{i})$ ; par suite l'adh\'erence de l'image de $W_{i}^{ss}$
est au plus de dimension
$2n+N^2-1.$  Ces ferm\'es recouvrent la fibre, ce qui d\'emontre l'\'enonc\'e
dans le cas o\`u $\hbox{\goth{c}}$  est r\'eduite.

Supposons maintenant que la conique $\hbox{\goth{c}}$  soit une droite
double $\ell^2.$ Le faisceau $\Theta$  est ou localement libre
auquel cas  son groupe d'automorphismes est de dimension 1, ou 
isomorphe \`a une  extension 
$$0\hfl{}{}{\cal O}_{\ell}(a)\hfl{}{}\Theta\hfl{}{}{\cal
O}_{\ell}(b)\hfl{}{}0$$
o\`u $a$ et $b$  sont des entiers tels que $a+b=n, a\geq b\geq 1.$  Une telle
extension peut \^etre scind\'ee, auquel cas le groupe des automorphismes
de $\Theta$ est de dimension $a-b+3$  ou 4 suivant que $a>b$ ou $a=b$,
ou non scind\'ee, auquel cas le groupe des automorphismes est de
dimension
$a-b+2.$ Les  extensions  non scind\'ees
  sont param\'etr\'ees  par l'espace projectif 
$P_{a,b}=\proj(\Ext^1({\cal O}_{\ell}(b),{\cal O}_{\ell}(a)))$, de
dimension $a-b+1$  et  on dispose d'une extension universelle
$\Theta_{a,b}$ sur
$P_{a,b}
\times \proj_2$. Consid\'erons le sch\'ema 
$W_{a,b}=\grass(2,H)\times\ul{\Isom}(H,\pr_{1*}(\Theta_{a,b}))$     et le
morphisme canonique
$f_{a,b}:W_{a,b}\hfl{}{} R_{n}.$ 
D\'esignons
par $W_{a,b}^{ss}$  l'image r\'eciproque de $R_{n}^{ss}$. On a donc
deux morphismes 
$$\diagram
{ W_{a,b}^{ss}&\hfl{f_{a,b}}{}& R_{n}^{ss}\cr
\vfl{\nu}{}&&\cr
P_{a,b}&&\cr
}$$ 
Consid\'erons un point $\lambda \in P_{a,b}$  auquel est associ\'ee une
extension
$\Theta_{\lambda}$  et la fibre $\nu^{-1}(\lambda).$  La fibre de la
restriction de 
$f_{a,b}$  \`a $\nu^{-1}(\lambda)$  est la trajectoire par 
l'action  du groupe
$\Aut(\Theta_\lambda)$ ;  cette action \'etant libre,  on voit que le
ferm\'e $\overline{f_{a,b}(\nu^{-1}(\lambda))}$  est de dimension
$\leq 2n -a+b-2+N^2.$  Par suite $\dim
\overline{f_{a,b}(W_{a,b})}\leq 2n+N^2-1.$ Les extensions scind\'ees se
traitent  comme ci-dessus, et leur contribution fournit des ferm\'es de
dimension $<2n+N^2-1$.  Comme tous ces ferm\'es recouvrent la
fibre de $\rho$  au-dessus du point d\'efini par la
conique $\hbox{\goth{c}}$  ceci ach\`eve la d\'emonstration du lemme
\ref{fibres}. 

 Pour voir que $R_{n}^{ss}$   est une vari\'et\'e normale, il reste  d'apr\`es
le crit\`ere de Serre \`a s'assurer  que son ensemble singulier est de
codimension
$\geq 2.$

\begin{lemme}  
{\rm (i)} Si $n=4$  le sch\'ema $R_{n}^{ss}$  est une vari\'et\'e lisse.

{\rm (ii)} Si $n>4,$
le ferm\'e des points  singulier   de $R_{n}^{ss}$   est de codimension
$\geq 3.$ 
\end{lemme}

\begin{proof} Il suffit d'apr\`es  la description locale de v\'erifier que le
ferm\'e de $R_{n}^{ss}$ des  points   $\Lambda=(\Gamma,\Theta)$  tels
que
$\Ext^2(\Theta,\Theta)\not=0$    est vide si $n=4$ ou de
codimension
$\geq 3$ si $n>5$. 
 On a \'evidemment $\Ext^2(\Theta,\Theta)=0$   si $\Theta$  est localement
libre,  ou si le support de $\Theta$  est d\'eg\'en\'er\'e en deux droites
distinctes. Si $\Theta$  est port\'e par une droite double $\ell^2$, et non
localement libre, il s'\'ecrit comme extension
$$0\hfl{}{}{\cal O}_{\ell}(a)\hfl{}{}\Theta\hfl{}{}{\cal
O}_{\ell}(b)\hfl{}{}0$$
avec $a\geq b-1,$  et alors $\Ext^2(\Theta,\Theta) \simeq \Ext^2({\cal
O}_{\ell}(a),{\cal O}_{\ell}(b))$  et par suite $$\dim
\Ext^2(\Theta,\Theta)=\left\{\matrix{0 \  {\rm si}\  a-b\leq 2\cr
a-b-2 \  {\rm si}\ a\geq b+3}\right.$$  Puisque $a\leq n-1$  et $a+b=n$, ce ferm\'e
est vide dans le cas $n=4;$  dans le cas $n\geq 5,$  l'\'enonc\'e r\'esulte du
lemme \ref{fibres}.
\end{proof}
Il en r\'esulte que $R_{n}^{ss}$  et par suite $S_{n}$  est une vari\'et\'e
normale; puisque sur l'ouvert des points stables l'action de $PGL(H)$
est libre, on voit que $S_{n}$  est un sch\'ema irr\'eductible et normal de 
dimension $2n+5.$ Ceci d\'emontre l'assertion (i).  Pour d\'emontrer
(ii),  il suffit de remarquer qu'au voisinage d'un point stable
$\Lambda=(\Gamma,\Theta)$  de
$S_{n}$  le sch\'ema $S_{n}$  est localement isomorphe (\cf Min He, [9])
au sch\'ema des z\'eros d'un morphisme  d\'efini sur un germe de sch\'ema lisse
d'espace tangent de Zariski $\Ext^1(\Lambda,\Lambda)$  et \`a valeurs
dans
$ \Ext^{2}(\Lambda,\Lambda). $ Or  $\sum_{i}(-1)^{i}\dim
\Ext^{i}(\Lambda,\Lambda)=-4-2n$  et puisque $\Lambda$  est stable, $\dim
\Hom(\Lambda,\Lambda)=1.$  Donc  
$$\dim \Ext^1(\Lambda,\Lambda)-\dim \Ext^2(\Lambda,\Lambda)=2n+5$$
Par suite,
$S_{n}$  est localement intersection
compl\`ete au voisinage de $\Lambda.$ Ceci ach\`eve la d\'emonstration de
la proposition \ref{S_{n}}.

\bigskip
\subsection{Le sch\'ema $P_{n}$ des faisceaux de Poncelet}

\begin{prop} \label{P_{n}}

{\rm (i)} Le sch\'ema $P_{n}$  est irr\'eductible et normal de
dimension
$2n+5.$ 

{\rm (ii)} L'ouvert des points stables de $P_{n}$  est de Cohen-Macaulay.

\end{prop}

\begin{proof}
Il r\'esulte de l'\'enonc\'e ci-dessus que $P_{n}$  est irr\'eductible de
dimension $2n+5.$
Soit $\ell$  un entier $\geq n-1,$  et $K$  un espace vectoriel de dimension
$\chi(\ell).$
 Le sch\'ema $P_{n}$  est obtenu par passage au quotient sous le groupe
r\'eductif $SL(K)$ d'un  sous-sch\'ema d\'eterminantiel $P'_{n}$ d'un  
ouvert lisse ${\cal H}^{ss}$ de
${\cal H}=
\Hilb(E,2-nh^2)
$,  o\`u $E={K\otimes{\cal O}(-\ell)}.$  Soit ${\bf F}$ le
faisceau universel  sur ${\cal H}^{ss}$. Il existe  des fibr\'es
vectoriels $A$  et $B$ sur ${\cal H}^{ss}$  de rangs respectifs $a$  et
$b=a+n-6$ et un morphisme de fibr\'es vectoriels  dont le noyau  est
$\pr_{1*}({\bf F}(1))$  et le conoyau $R^{1}\pr_{1*}({\bf F}(1)).$ Le
sch\'ema 
$P'_{n}$  est d\'efini par l'id\'eal  des mineurs $a-1\times a-1$ de ce
morphisme ; il est  de dimension $2n+4 +\chi(\ell)^{2}.$  
Consid\'erons en effet le morphisme canonique $R_{n}^{ss}\hfl{}{}M_{n}$  qui
associe au couple $(\Gamma,\Theta)$  la classe du faisceau $F$  d\'efini par
l'extension
$$0\hfl{}{}\Gamma^*\otimes
{\cal O}\hfl{}{}F(1)\hfl{}{}\check{\Theta}\hfl{}{}0$$ Cette construction
d\'efinit en fait  une famille  ${\bf F}'$ de
 faisceaux
$F
$  param\'etr\'ee par
$R_{n}^{ss}$  de sorte que si on introduit  le fibr\'e des rep\`eres de
$R'=\Isom(K\otimes {\cal O}_{R_{n}^{ss}},\pr_{1*}({\bf F}'))$  on
obtient un diagramme commutatif
$$\diagram{R'&\hfl{q}{}& R_{n}^{ss}\cr
\vfl{p}{}&&\vfl{}{}\cr
{\cal H}^{ss}&\hfl{r}{}&M_{n}\cr
}
$$
Il r\'esulte du th\'eor\`eme pr\'ec\'edent que $R'$  est un  sch\'ema
de Cohen-Macaulay irr\'eductible et normal de dimension $2n+4 +N^2+
\chi(\ell)^{2}.$ Les morphismes verticaux sont
$GL(H)-$\'equivariants; le morphisme $p$  a pour image $P'_{n}$. Au-dessus
de  l'ouvert
$U'\subset {\cal H}^{ss}$ des faisceaux quotients non sp\'eciaux le
morphisme
$p$
 est  un fibr\'e principal de groupe structural $GL(H)$ et de
base le sous-sch\'ema 
$P'_{n}\cap U'.$  Donc $\dim P'_{n} = 2n+4 + \chi(\ell)^{2}$ 
est irr\'eductible de codimension $2(n-4).$ Localement c'est l'image
r\'eciproque par un morphisme de sch\'emas lisses d'un sous-sch\'ema de
Cohen-Macaulay; puisqu'il est de la bonne dimension,   il est donc de
Cohen-Macaulay.  L'ensemble singulier de $R'$  est de codimension $\geq 3$ 
et invariant par
$GL(H),$ et la trace dans $U'$  de son image est l'ensemble singulier de
$P'_{n}\cap U'.$ Donc  cet ensemble singulier est de codimension $\geq 3.$
Il en r\'esulte que $P'_{n}$ est une vari\'et\'e normale d'apr\`es le crit\`ere de
Serre.  Sur l'ouvert des points stables de $P_{n}'$, qui est non vide, 
l'action  de
$PGL(K)$ est libre; le quotient
$P_{n}$  est donc une vari\'et\'e irr\'eductible et normale de dimension $2n+5.$
Ceci d\'emontre (i).

Soit $x \in P_{n}$  un point stable,  et ${\cal O}_{x}$  l'anneau local
de $M_{n}$  en $x.$  Au-dessus de  $\Spec (\hat{\cal O}_{x})\times
\proj_2$, il existe une famille universelle, ce qui dispense de se placer
sur le sch\'ema de Hilbert.  Le m\^eme argument que ci-dessus montre qu'alors
que le compl\'et\'e  de l'anneau local de $P_{n}$  en $x$  est de
Cohen-Macaulay.  Puisque profondeur et dimension sont conserv\'ees par
passage au compl\'et\'e, on voit que $P_{n}$  est lui-m\^eme  de
Cohen-Macaulay en $x.$  D'o\`u (ii). 
\end{proof}  
Consid\'erons le morphisme $\pi:S_{n}\hfl{}{}P_{n}.$  L'image r\'eciproque
du sous-sch\'ema $\Sigma_{n}$ des faisceaux sp\'eciaux est un diviseur de
Weil;  il  correspond aux syst\`emes coh\'erents $(\Gamma,\Theta)$ stables
tels que le faisceau
$\Theta$   sous-jacent  soit  extension 
$$0\hfl{}{}{\cal O}_{\ell'}(n-1)\hfl{}{} \Theta\hfl{}{} {\cal
O}_{\ell''}(1)\hfl{}{} 0$$ 
o\`u $\ell'$  et $\ell''$  sont deux droites de $\proj_2.$
 Un tel syst\`eme coh\'erent est stable; on obtient ainsi un  ferm\'e de
$S_{n}$ d\'efini par  les syst\`emes coh\'erents satisfaisant \`a la condition
d\'eterminantielle
$H^1(\Theta(-3))\not=0.$  C'est un  diviseur de Weil; nous ignorons si ce
diviseur est un diviseur de Cartier si $n\geq 5$.
  C'est vrai pour $n=4.$  Dans ce cas, on a $P_{4}=M_{4}$,
le sous-sch\'ema $\Sigma_{4}$  est une sous-vari\'et\'e lisse de
codimension 3 qui \'evite les points singuliers de $M_{4}$; le r\'esultat
suivant est d\'emontr\'e dans {\rm [15]}: 

\begin{th} \label{eclate}Le morphisme $\pi:S_{4}\hfl{}{} M_{4}$  est
l'\'eclat\'e de
$M_{4}$  le long de $\Sigma_{4}.$

\end{th}

Les seuls points singuliers de $S_{4}$  sont les points  strictement
semi-stables.  Ceci ne reste pas vrai pour $S_{n},$ pour $n\geq 5.$

 \section{Courbe des droites de saut}
 Soient 
$(\Gamma,\Theta)$  un syst\`eme coh\'erent semi-stable  de classe de
Grothendieck
$c=2h+nh^2$ et  tel que
$\dim
\Gamma=2,$  et
$F$  le faisceau  de Poncelet de classes de Chern $(0,n)$ 
associ\'e  dans la correspondance  d\'ecrite au \S 4.1 ci-dessus. 
La courbe $\beta_{F}$ des droites de saut de $F$, dite courbe de
Poncelet, peut en fait  se d\'ecrire directement en termes de
 syst\`emes coh\'erents;  on introduit pour ceci le diagramme
standard associ\'e \`a la vari\'et\'e d'incidence des couples (droites, points):
$$\diagram{ D&\hfl{\pr_{2}}{}&\proj_{2}\cr
\vfl{\pr_{1}}{}&&\cr
\proj_{2}^{*}\cr}$$ 
La courbe des droites de
saut $\hbox{\goth{q}}=\beta_{F}$ de
$F$  est d\'efinie par l'id\'eal de Fitting du faisceau pur de dimension un 
 sur le plan projectif dual $R^{1}\pr_{1*}(\pr_{2}^*(F(-1)))$.
\subsection{Le faisceau ${\cal G}=\pr_{1*}(\pr_{2}^*(\Theta))$ }
\begin{lemme} \label{G}{\rm (i) } Sur le plan projectif dual, l'image
directe $R^{1}\pr_{1*}(\pr_{2}^*(\Theta))$ est nulle, et le faisceau
${\cal G}=\pr_{1*}(\pr_{2}^*(\Theta))$  
est un faisceau alg\'ebrique coh\'erent sans
torsion de rang 2, de classes de Chern
$(n,\ds{1\over 2}n(n+1))$. 

{\rm (ii)} Ce faisceau
est non singulier  en tout point  $\check{\ell}$ d\'efinissant une  droite
$\ell$  non   contenue dans le support de
$\Theta.$  

{\rm (iii)} Soit ${\check{\ell}} $  un point de $\proj_2^*, $ d\'efinissant
une droite $\ell$  de $\proj_2.$ Le morphisme canonique
$${\cal G}_{\check{\ell}}\otimes_{{\cal O}_{\check{\ell}}}k
\hfl{}{}H^0(\Theta|_{\ell})$$  est un isomorphisme; de plus, on a un
isomorphisme  canonique
$$\Tor_{1}({\cal G}_{\check{\ell}},k)\simeq
H^{0}(\ul{\Tor}_{1}(\Theta,{\cal O}_{\ell})).$$
\end{lemme}

\begin{proof}
De la condition de semi-stabilit\'e de $(\Gamma,\Theta )$, on d\'eduit
$H^1(\Theta (i))=0 $ pour $i\geq -2$, ce qui exprime le fait que le faisceau
$\Theta $ n'est pas trop instable.  Le faisceau $\Theta$  \'etant pur, il
en est de m\^eme de son image r\'eciproque sur $\proj_2^* \times \proj_2$;
ceci implique la condition de transversalit\'e
$\Tor_{i}^{{\cal O}_{{\proj_2}^* \times {\proj_2}}}({\cal
O}_{D},\Theta)=0$
pour $i>0.$
  De la r\'esolution de ${\cal O}_{D}$ sur $\proj_2^*\times
\proj_2$ 
$$0\hfl{}{}{\cal O}(-1,-1)\hfl{}{} {\cal O}\hfl{}{} {\cal O}_{D}\hfl{}{}
0
\eqno(5.1)\qquad$$ et du fait que $H^1(\Theta (-1))=H^{1}(\Theta)=0$
on d\'eduit que  $R^{1}\pr_{1*}(\pr_2^*(\Theta ))=0$; de plus, on a la
pr\'esentation de l'image directe ${\cal G}=\pr_{1*}(\pr_2^*(\Theta ))$
$$0\hfl{}{} K\otimes{\cal O}_{\proj_2^*}(-1)\ \hfl{\alpha}{}\  H\otimes
{\cal O}_{\proj_{2}^*}\hfl{}{} {\cal G}\hfl{}{} 0 \eqno{(5.2)}\qquad$$
o\`u  $K=H^0(\Theta (-1))$
 et $ H=H^0(\Theta)$; d'apr\`es la formule de
Riemann-Roch, ces espaces vectoriels sont de dimensions
respectives  $n$ et
$n+2$. Il en r\'esulte que le faisceau ${\cal G}$  est un faisceau de rang
2 et de profondeur
$\geq 1$  en tout point $\ell$, et de classes de Chern $(n,\ds{1\over
2}n(n+1))$.  Les assertions (ii) et (iii) sont  \'evidentes si le support
de
$\Theta$ est lisse, car alors $\pr_2^*(\Theta)$  est fini et plat sur
$\proj_{2}^*$ et elles r\'esultent du th\'eor\`eme de changement de base. Ce
n'est plus le  cas si le faisceau 
$\Theta$  a son support singulier; pour obtenir  le fait  que ${\cal G}$
est sans torsion,  il suffit d'apr\`es le crit\`ere de Serre  de remarquer 
que l'ensemble singulier de
${\cal G}$  est fini.  Par d\'evissage, il suffit en fait de v\'erifier  
le lemme suivant,  qui est beaucoup plus pr\'ecis, et conduit en m\^eme
temps \`a l'assertion (iii):

\begin{lemme}\label{droites} Soit $\ell$  une droite de $\proj_2$, et
$i$  un entier
$\geq 0$.
 D\'esignons par $\check{\ell}$  le point de $\proj_2^*$ 
d\'efini par $\ell,$ et par $Q=Q_{\proj_2^*}$  le fibr\'e quotient universel
de rang 2 sur le plan projectif dual; soit $t$  une section de $Q$ 
dont le sch\'ema des z\'eros est le point $\check{\ell}.$ 

{\rm (i)} L'isomorphisme canonique  $H^{0}({\cal O}_{\proj_2}(1))\simeq
H^{0}(Q_{\proj_2^*})$  induit  un isomorphisme de faisceaux
$$\pr_{1*}(\pr_{2}^*({\cal
O}_{\ell}(i)))\simeq \hbox{\goth{m}}_{\check{\ell}}^i(i)$$

{\rm (ii)}   Le diagramme  associ\'e \`a cet isomorphisme
$$\diagram{
H^{0}({\cal O}(i)) &\hfl{}{}&
H^{0}({\cal O}_{\ell}(i)) \cr
\vfl{\wr}{}&&\vfl{\wr}{}\cr
H^{0}(S^{i}Q) &\hfl{t^i}{}&
H^{0}(\hbox{\goth{m}}_{\check{\ell}}^{i}(i))\cr}$$  dans lequel la	
seconde fl\`eche horizontale est d\'efinie 
  par le morphisme $Q\hfl{}{} {\cal
O}(1)=\wedge^2Q$  associ\'e \`a $t$  est commutatif.

{\rm (iii)}  On a 
$R^{1}\pr_{1*}(\pr^*_{2}({\cal O}_{\ell}(i)))=0.$
\end{lemme}

\begin{proof} On a $S^{i}Q= \pr_{1*}(\pr_2^*({\cal O}(i))$. Supposons
choisies des coordonn\'ees
$(x,y,t)$  dans
$V^*=H^0({\cal O}_{\proj_2}(1))$ de sorte que l'\'equation de ${\ell}$ 
soit $t=0 .$ La r\'esolution standard de ${\cal O}_{\ell}$  conduit
imm\'ediatement \`a l'assertion (iii), et  fournit en outre, par application
du foncteur $\pr_{1*}\pr_2^*$ une   r\'esolution
$$0\hfl{}{}S^{i-1}Q\hfl{\tau}{} S^iQ\hfl{}{}
\pr_{1*}(\pr_{2}^*({\cal O}_{\ell}(i))) \hfl{}{}0$$
 o\`u le
morphisme $\tau$  est la multiplication par  $t,$ consid\'er\'ee 
 comme section de $ V^*=H^0(Q).$ On a un morphisme $S^{i}Q\hfl{}{}
\hbox{\goth{m}}^{i}_{\check{\ell}}(i)$, qui s'annule sur l'image de
$\tau;$ il induit donc un morphisme $$\pr_{1*}(\pr_{2}^*({\cal
O}_{\ell}(i)))\hfl{}{}\hbox{\goth{m}}_{\check{\ell}}^i(i)$$ 
qui  est \'evidemment  un isomorphisme en dehors
du point 
$\check{\ell}$. Le premier membre est un faisceau  sans torsion; par
suite,  ce morphisme est donc injectif.  Les deux membres ayant les m\^emes
classes de Chern, c'est un isomorphisme. Ceci d\'emontre l'assertion (i).
Quant \`a l'assertion (ii),  il suffit de constater que le diagramme  de
faisceaux sur
$\proj_2^*$
$$\diagram{
\pr_{1*}\pr_{2}^*({\cal O}(i)) &\hfl{}{}&
\pr_{1*}\pr_2^*({\cal O}_{\ell}(i)) \cr
\vfl{\wr}{}&&\vfl{\wr}{}\cr
S^iQ &\hfl{t^i}{} &\hbox{\goth{m}}^{i}_{\check{\ell}}(i)\cr}$$ est
commutatif.  Ceci r\'esulte trivialement de la construction de l'isomorphisme
de droite.
\end{proof} \goodbreak
{\it Fin de la d\'emonstration du lemme \ref{G}}

Le lemme \ref{droites} ci-dessus entra\^{\i}ne bien s\^ur l'assertion (ii) du
lemme \ref{G}. Pour v\'erifier l'assertion (iii), 
on utilise  l'isomorphisme  donn\'e dans la cat\'egorie
d\'eriv\'ee   des  espaces vectoriels par la formule de projection
$$R{\pr_{1*}}(\pr_{2}^*(\Theta))_{{\cal O}_{\check{\ell}}} \otimes_{{\cal
O}_{\check{\ell}}}^{L} k
\simeq R\pr_{1*}(\pr^*_{2}(\Theta )_{{\cal
O}_{\check{\ell}}}\otimes^{L}{\cal O}_{\ell}))$$ Compte-tenu du fait que
la projection $\pr_2$  est un morphisme plat, le second membre est aussi
$R\pr_{1*}(\pr^*_{2}(\Theta
\otimes^{L}{\cal O}_{\ell}))$. La cohomologie du premier membre est
 l'aboutissement d'une suite spectrale 
$$'E_{2}^{p,q}= \Tor^{{\cal
O}_{\check{\ell}}}_{-p}(R^q{\pr_{1*}}(\pr_{2}^*(\Theta)),k))$$
  L'assertion (iii) est \'evidente si le support de $\Theta$  est lisse.
  Si ce support 
est une conique singuli\`ere, le fait que  le syst\`eme coh\'erent
$(\Gamma,\Theta)$ soit  semi-stable implique qu'il existe une suite
exacte
$0\hfl{}{} \Theta'\hfl{}{}\Theta\hfl{}{}\Theta''\hfl{}{}0$ o\`u $\Theta'$ 
et $\Theta''$  sont des faisceaux inversibles de degr\'es positifs sur des
droites $\ell'$ et $\ell''$ de $\proj_2.$ Il r\'esulte du lemme \ref{droites}
ci-dessus que
$'E_{2}^{p,q}$  est nul sauf si $q=0$  et $p=0$  ou $-1$.  Consid\'erons la
deuxi\`eme suite spectrale, de m\^eme aboutissement, 
$$''E_{2}^{p,q}= R^{p}{\pr_{1*}}(\pr_{2}^*(\ul{\Tor}_{-q}(\Theta,{\cal
O}_{\ell})))$$ Du fait que $\ul{\Tor}_1({\cal
O}_{\ell},{\cal O}_{\ell})= {\cal O}_{\ell}(-1) $ on d\'eduit encore du
lemme \ref{droites} que
$''E_{2}^{p,q}=0$  sauf si $p=0$  et $q=0$  ou $q=-1.$ Ceci conduit  
au r\'esultat attendu. \cqfd
\bigskip

\subsection{Description de  la courbe  $\hbox{\goth{q}}=\beta_{F}$}

 Consid\'erons le morphisme canonique sur le plan projectif dual
$$\varphi:\Gamma\otimes {\cal O}_{\proj_{2}^*}\hfl{}{}
\pr_{1*}(\pr_{2}^*(\Theta))={\cal G}$$
  \begin{prop} Le morphisme $\varphi$ est  g\'en\'eriquement
injectif, et son conoyau a m\^eme id\'eal de Fitting que le faisceau
$R^{1}\pr_{1*}(\pr_2^*(F(-1)))$.

\end{prop}
\begin{proof} Le noyau de l'homomorphisme canonique 
$$\Gamma \otimes {\cal O}_{\proj_2}\hfl{}{} \Theta$$  est isomorphe \`a
$F^*(-1),$ et son conoyau \`a $\ul{\Ext}^1(F,{\cal O}(-1)).$ En
consid\'erant une droite g\'en\'erique, on est ramen\'e pour v\'erifier la
premi\`ere assertion, \`a  constater que
$H^{0}(F^*(-1)|_{\ell})=0,$  ou ce qui revient au m\^eme par dualit\'e que
$H^{1}(F(-1)|_{\ell})=0.$  Il suffit de prendre pour 
$\ell$  une droite qui n'est pas de saut pour $F.$
Consid\'erons maintenant les projections $\pr_{i}$ figurant dans le
diagramme  standard.
 Les faisceaux $R^1\pr_{1*}(\pr_2^*(F(-1))$  et
$\ul{\Ext}^1_{\pr_{1}}(\pr_{2}^*(F),{\cal O}_{D}(0,-1))$  sont deux
faisceaux coh\'erents de dimension 1  sur $\proj_2^*$  de m\^eme id\'eal de
Fitting (\cf  [10], lemme 3.4): ceci peut se voir par
exemple en utilisant une pr\'esentation de
$F$  comme faisceau de cohomologie d'une monade, d'o\`u on d\'eduit  en fait un
isomorphisme
$$\ul{\Ext}^1_{\pr_{1}}(\pr_{2}^*(F),{\cal O}_{D}(0,-1))\simeq
\ul{\Ext}^1(R^1\pr_{1*}(\pr_2^*(F(-1))),{\cal
O}_{\proj_2^*}(-1))\eqno{(5.3)}\qquad$$  Soit
$\check{\Theta}= \ul{\Ext}^1(\Theta,{\cal O}).$ En appliquant le foncteur
$\ul{\Hom}_{\pr_{1}}(\pr_2^*(-),{\cal O}_{D})$ \`a la suite exacte 
$$0\hfl{}{} \Gamma^*\otimes {\cal O}\hfl{}{} F(1) \hfl{}{} \check{\Theta
}\hfl{}{} 0
\eqno{(5.4)}\qquad$$ on obtient la suite exacte sur $\proj_2^*$
$$0\hfl{}{} \Gamma\otimes {\cal O}_{\proj_2^*}\hfl{v}{}
\ul{\Ext}^1_{\pr_{1}}(\pr_{2}^*(\check{\Theta}) ,{\cal
O}_{D})\hfl{}{}\ul{\Ext}^1_{\pr_{1}}(
\pr_2^*(F(1)),{\cal O}_{D})\hfl{}{}0$$ En consid\'erant la r\'esolution (5.1)
de la vari\'et\'e d'incidence $D$, on obtient  exactement pour le faisceau
$\ul{\Ext}^1_{\pr_{1}}(\pr_{2}^*(\check{\Theta}) ,{\cal O}_{D})$, la
pr\'esentation (5.2) , ce qui montre qu'il est isomorphe \`a ${\cal G}$, le
morphisme $v$ s'identifiant au morphisme induit par le morphisme
d'\'evaluation, c'est-\`a-dire \`a $\varphi$. 
\end{proof}

 \begin{cor}\label{droitesdesaut} La courbe des droites de saut de $F$ 
est donn\'ee par l'id\'eal de Fitting du conoyau du morphisme 
$\varphi:\Gamma\otimes {\cal O}_{\proj^*_2}\hfl{}{} {\cal G}.$
\end{cor}

Il r\'esulte de l'isomorphisme (5.3) que  le conoyau  de $\varphi$ est 
pur de dimension 1, et que son support sch\'ematique est une quartique  qui
co\"{\i}ncide avec la quartique
$\hbox{\goth{q}}=\beta_{F}$ des droites de saut de $F.$  Il r\'esulte du
lemme
\ref{G} qu'une 
droite $\ell$ d\'efinit un point de $\beta_{F}$  si et seulement
si l'application lin\'eaire
$\varphi_{\ell}:s\mapsto s|_{\ell}$
$$\Gamma \hfl{}{} H^0(\Theta|_{\ell})$$ n'est pas inversible. 

\begin{cor}  Si  le support de $\Theta$  est une conique singuli\`ere
$\hbox{\goth{c}}$, la courbe $\hbox{\goth{q}}=\beta_{F}$  contient  les
points de
$\proj_{2}^*$ correspondant aux droites $\ell$ contenues dans  la
conique et telles que $\chi(\Theta|_{\ell})\geq 3$. 
\end{cor}

\begin{cor}\label{singuliers-evidents}
Si  le support de $\Theta$  est une conique singuli\`ere
$\hbox{\goth{c}}$, une droite  $\ell$ contenue dans  la
conique et telle que $\chi(\Theta|_{\ell})\geq 4$ d\'efinit un point
singulier $\check{\ell}$ de la courbe $\hbox{\goth{q}}=\beta_{F}$. 
\end{cor}

\subsection{Etude des points singuliers}

Il est fondamental pour la suite de pouvoir d\'ecider 
qu'une  droite de saut $\ell$ de $F$  d\'efinit ou non un point singulier
 de la quartique de L\"uroth associ\'ee.  C'est facile lorsque $F$  provient
d'un syst\`eme lin\'eaire
$\Lambda=(\Gamma,\Theta)$ de $\Syst(2h+nh^2,2)$ sans point de base  sur
une conique
$\hbox{\goth{c}}$  \'eventuellement singuli\`ere. Commen\c{c}ons par \'etendre
cette notion aux syst\`emes coh\'erents $(\Gamma,\Theta)$  tels que
$\Theta$  ne soit plus obligatoirement localement libre sur son support.

\begin{defi} On dit que le syst\`eme coh\'erent $(\Gamma,\Theta)$  est sans
point de base si le morphisme d'\'evaluation  ${\rm ev}:\Gamma \otimes
{\cal O}_{\proj_2}\hfl{}{} \Theta$  est surjectif.
\end{defi} 
Si le syst\`eme coh\'erent  $(\Gamma,\Theta)$ est sans point de base, le
faisceau
$F$ associ\'e est localement libre. Si le syst\`eme coh\'erent
$(\Gamma,\Theta)$ est stable et a un  point de base 
$a,$   on peut trouver un
 sous-faisceau 
$\Theta'\subset \Theta$ de multiplicit\'e 2   , de caract\'eristique
d'Euler-Poincar\'e $\chi(\Theta')=n-1$    tel que $\Gamma\subset
H^0(\Theta').$  Le syst\`eme coh\'erent
$(\Gamma,\Theta')$  est encore semi-stable, et d\'efinit une courbe
$\hbox{\goth{q}}'$  de degr\'e
$n-1.$  Il r\'esulte  du corollaire \ref{droitesdesaut}
que la courbe
$\hbox{\goth{q}}=\beta_{F}$  est d\'ecompos\'ee  en la courbe
$\hbox{\goth{q}}'$  et une droite
$$\hbox{\goth{q}}=\hbox{\goth{q}}'+ \check{a}$$ o\`u $\check{a}$  est
la droite du plan projectif dual d\'efinie par le point $a.$  En
particulier, la courbe $\hbox{\goth{q}}$  est singuli\`ere. C'est
aussi le cas si on part d'un syst\`eme coh\'erent $(\Gamma,\Theta)$
strictement semi-stable:  ceci impose que $n$ est pair, et la courbe
$\hbox{\goth{q}}$  est alors l'union de deux courbes de degr\'e $m={n\over
2}$, chacune d'elles \'etant  l'union de $m$  droites concourantes.
 L'\'enonc\'e suivant
g\'en\'eralise l'\'enonc\'e bien connu pour les courbes de Poncelet associ\'ees \`a
des syst\`emes lin\'eaires sur des courbes lisses (Maruyama [20], Trautmann
[24], Vall\`es [26]).  Il est important pour la suite d'examiner les
points  singuliers des courbes de Poncelet associ\'ees \`a des syst\`emes
coh\'erent semi-stables dont le support sch\'ematique est une conique
singuli\`ere.

\begin{prop} \label{points-singuliers} Soit $z \in
H^{0}({\cal O}(1))$ non nul et $\ell $  la droite d'\'equation
 $z=0.$  On suppose que

 {\rm (i)} la	 droite $\ell$  n'est pas associ\'ee \`a $\Theta;$ 
 
  {\rm (ii)} le syst\`eme coh\'erent $(\Gamma,\Theta)$  est semi-stable
et sans point de base. 

La droite 
$\ell$ est d\'efinit un point singulier  $\check{\ell}$ de la courbe
$\hbox{\goth{q}}$ des droites de saut $\hbox{\goth{q}}=\beta_{F}$ si et
seulement si 
$$z^{2}H^{0}(\Theta (-2))\cap \Gamma \not=0.$$
\end{prop}

\begin{proof} Puisque la droite $\ell$  n'est pas associ\'ee \`a $\Theta$, le
morphisme de multiplication par $z$ $$\Theta(-1)\hfl{z}{} \Theta$$
est injectif, et  il en est de m\^eme  du morphisme induit sur les
espaces de sections globales. On a alors $\dim H^{0}(\Theta|_{\ell})=2$ 
et $\ell$  est un point non singulier pour ${\cal G}.$ Consid\'erons  la
filtration d\'ecroissante
$\Gamma_{i}=\Gamma
\cap
z^{i} H^{0}(\Theta(-i))$. Puisque $(\Gamma,\Theta)$  est sans point de
base, on a $ \dim zH^{0}(\Theta(-1)) \cap \Gamma \leq 1.$  
Si cette intersection est nulle, le morphisme  de restriction
$\varphi_{\ell}:\Gamma\hfl{}{} H^{0}(\Theta|_{\ell})$ est
injectif et donc c'est un isomorphisme.  Donc $\ell$ n'appartient pas \`a
la courbe
$\hbox{\goth{q}}$ des droites de saut.

Supposons $zH^{0}(\Theta(-1))\cap \Gamma\not=0.$  Alors $\ell$  d\'efinit
un point de $\hbox{\goth{q}}.$
 Il existe  une section non nulle $s \in
\Gamma$   de la forme
$s=zt$, o\`u $t $  est une section de $\Theta(-1).$  L'espace tangent
de Zariski \`a $\hbox{\goth{q}}$ en $\check{\ell}$  est le noyau de
l'application lin\'eaire tangente de Petri $d_{\check{\ell}}\varphi$ de
$\varphi$ au point
$\check{\ell}$
$$T_{\check{\ell}}(\proj_{2}^*)=H^{0}({\cal O}_{\ell}(1))\hfl{}{}
\Hom(\ker \varphi_{\ell},\coker \varphi_{\ell}).$$ Le noyau de
$\varphi_{\ell}$  est  le sous-espace vectoriel
engendr\'e par $ s=tz$;  l'application lin\'eaire
tangente  de Petri
$d_{\check{\ell}}\varphi$ est d\'efinie pour $\tau \in H^{0}({\cal
O}_{\ell}(1))$ par 
$$d_{\check{\ell}}\varphi(\tau)(s)= p(\tau t
|_{\ell}) $$  o\`u $p$  est la projection canonique
$H^0(\Theta|_{\ell})\hfl{}{}\coker \varphi_{\ell}.$ 

Supposons $z^2H^{0}(\Theta(-2))\cap \Gamma\not=0.$
Il s'agit de v\'erifier que $\check{\ell}$ est un point singulier de
$\hbox{\goth{q}}.$ 
La section $s$  est alors divisible par $z^2$, et $t$  est
divisible par $z.$ Cette diff\'erentielle est alors nulle; et par
cons\'equent le point $\check{\ell}$ est un point singulier de
$\hbox{\goth{q}}=\beta_{F}$. R\'eciproquement, si cette diff\'erentielle est
nulle, 
$\tau t|_{\ell}$  appartient \`a $\Im \varphi_{\ell}$  pour tout
$\tau.$ Nous allons en d\'eduire que $t|_{\ell}$  est nulle.

Soit $\hbox{\goth{c}}$  le support de $\Theta.$ Remarquons tout d'abord
que  le morphisme de restriction
$$H^{0}({\cal O}_{\ell}(1))
 \hfl{}{} H^{0}({\cal O}_{\hbox{\goth{c}}\cap \ell}(1))$$  est
un isomorphisme. Si $\ell$  ne rencontre aucun des points
singuliers de $\Theta,$ l'espace vectoriel des sections
$H^{0}(\Theta|_{\ell})$  est un module libre de rang 1 
sur l'alg\`ebre 
$H^{0}({\cal O}_{\hbox{\goth{c}}\cap \ell})$. L'image de
$\varphi_{\ell}$  est un sous-espace vectoriel de dimension 1, qui,
d'apr\`es l'hypoth\`ese (ii) n'est contenu
 dans aucun
des sous-modules $K_{x}$ de $H^{0}(\Theta|_{\ell})$ des \'el\'ements
s'annulant en un point $x$ de
$\hbox{\goth{c}}\cap
\ell.$ Le stabilisateur de $\Im \varphi_{\ell}$  dans le groupe  des
\'el\'ements inversibles $H^{0}({\cal O}_{\hbox{\goth{c}}\cap \ell})^*$ est 
r\'eduit au sous-groupe $k^*$  des homoth\'eties. Donc, son orbite n'est pas
r\'eduite \`a un point. Il en r\'esulte que
$t|_{\ell}=0$  et par suite la section 
$s$ appartient \`a $\Gamma \cap z^2H^0(\Theta(-2)) $.
Si $\ell$  rencontre un point singulier $a$ de $\Theta,$  le module
$\Theta|_{\ell}$  est de longueur 2, de support le point $a$ et non
isomorphe \`a
${\cal O}_{\hbox{\goth{c}}\cap \ell}.$  Donc $\Theta|_{\ell
}$   est
isomorphe \`a
$k(a)\oplus k(a).$ Le morphisme  de restriction 
$H^0(\Theta|_{\ell})\hfl{}{}\Theta
\otimes k(a)$ est alors un isomorphisme. Puisque $(\Gamma,\Theta)$  est
sans point de base, l'application  lin\'eaire
$\varphi_{\ell}:\Gamma\hfl{}{} H^0(\Theta|_{\ell})$  est un isomorphisme
et donc $\ell$  n'est pas de saut.
\end{proof}

\subsection{Exemple} 

 On voit en particulier  que si le support de
$\Theta$  est une conique lisse $\hbox{\goth{c}}$, le faisceau $\Theta$ 
est un  faisceau inversible de degr\'e $n+1$  sur $\hbox{\goth{c}}$; si
$\Gamma$  contient une section qui a $n+1$ z\'eros distincts, la courbe de
saut 
$\hbox{\goth{q}}=\beta_{F}$  passe par les ${n(n+1)\over2}$ sommets du
 polyg\^one  de $\proj_2^*$ \`a $n+1$ c\^ot\'es   d\'etermin\'e par ces $n+1$
points. Ce polyg\^one est circoncrit \`a la conique duale de
$\hbox{\goth{c}}$.  
De plus, l'interpr\'etation ci-dessus  permet d'\'etudier les singularit\'es 
de la courbe $\beta_{F}$: ceci entra\^{\i}ne   que si
$(\Gamma,\Theta)$ est un syst\`eme lin\'eaire de degr\'e $n+1$ sans point de
base sur une conique lisse, et si
$\Gamma$ ne contient pas de sections ayant 2 z\'eros doubles, la courbe
$\hbox{\goth{q}}=\beta_{F}$  est lisse (\cf 
 Maruyama [20],  Trautmann [25], Vall\`es [27]). 

R\'eciproquement,  \'etant donn\'e un
 polyg\^one  \`a $n+1$ c\^ot\'es
 dans le plan projectif dual,  circonscrit \`a une conique lisse
$\hbox{\goth{c}}^*$,  ce polyg\^one  d\'efinit  $n+1$ points distincts sur
une conique lisse
$\hbox{\goth{c}}$ de
$\proj_{2}$;
 ces points  sont les z\'eros  d'une section $s$  d'un fibr\'e inversible
$\Theta$ de degr\'e $n+1$  sur $\hbox{\goth{c}}.$ Les syst\`emes lin\'eaires
$\Gamma\subset H^0(\Theta)$  qui contiennent la section $s$ 
constituent dans la grassmannienne $\grass(2,H^0(\Theta))$  des
sous-espaces de dimension 2   de $H^0(\Theta)$ un espace projectif de
dimension $n$; d'apr\`es le corollaire \ref{droitesdesaut}, en restriction
\`a cet espace projectif, le morphisme
$\beta$   est une homographie; cette homographie  identifie cet espace
projectif avec l'espace projectif des courbes de degr\'e $n$ passant par
les  sommets du pentagone donn\'e (\cf Barth [3], Lemma 2). Par suite,
toute  courbe 
$\hbox{\goth{q}}$ de
$\proj_{2}^{*}$ passant par les ${n(n+1)\over 2}$ sommets d'un tel
polyg\^one d\'efinit un point de l'image  de $\beta.$ Une telle courbe
$\hbox{\goth{q}}$ quand elle est lisse, est appel\'ee classiquement une
courbe de Poncelet associ\'ee \`a la conique
$\hbox{\goth{c}}^*.$

\section{Faisceaux dont la courbe des droites de saut est singuli\`ere.}
On se limite dans les sections 6-9  au cas $n=4:$  la g\'en\'eralisation aux
faisceaux de Poncelet de classe de Chern quelconque n'est pas
 \'evidente; nous ignorons par exemple  si $S_{n}$  est
localement factorielle. Dans le cas $n=4$, on a vu l'espace de modules
$S_{4}=\Syst(2h+4h^2,2)$  s'identifie \`a l'\'eclat\'e de $M_{4}$  le long
d'une sous-vari\'et\'e lisse
$\Sigma_{4}$ de codimension 3  qui \'evite les points singuliers. Tout
comme $M_{4}$  c'est une vari\'et\'e irr\'eductible, normale, localement
factorielle. L'image ${\cal L}$ du
morphisme de Barth $\beta:M_{4}\hfl{}{}{\cal C}_{4}$ est une
hypersurface dont les points repr\'esentent les courbes de Poncelet de
degr\'e 4. Ces quartiques sont appel\'ees quartiques de L\"uroth; 
l'hypersurface
${\cal L}$ sera appel\'ee  hypersurface de L\"uroth.

On pose dans la
suite
$c=2h+4h^2.$ 
On note $\gamma= \beta\circ \pi: S_{4}\hfl{}{}{\cal C}_{4}$ le
morphisme qui associe \`a un syst\`eme coh\'erent la courbe des droites de
saut. 

\subsection{Le groupe de Picard de $S_{4}=\Syst(c,2)$}

On sait d\'ecrire le groupe de Picard de $M_{4};$  c'est un groupe ab\'elien
libre \`a 2 g\'en\'erateurs ([5],[22]).  Il r\'esulte du th\'eor\`eme \ref{module}
que  le groupe de Picard de $S_{4}$  est un groupe ab\'elien libre \`a 3
g\'en\'erateurs. Le diviseur exceptionnel $\hbox{\goth{e}}$ est exactement
le diviseur des syst\`emes coh\'erents $\Lambda=(\Gamma,\Theta)$  tels que
$\Theta$  soit instable. On se propose  dans ce paragraphe de d\'ecrire une
base explicite de
$\Pic(S_{4}).$

On consid\`ere le groupe ab\'elien libre 
$\widetilde{K}(\proj_2)=K(\proj_{2})
\oplus
\entr,$  muni de la forme quadratique unimodulaire $(u,m)\mapsto
\chi(u^2)+m^2$ .  Consid\'erons une famille $\Lambda=(\Gamma,\Theta) $ de
syst\`emes coh\'erents de
$S_{4}$ param\'etr\'ee par une vari\'et\'e alg\'ebrique lisse $S:$  ceci
signifie 
 que $\Theta $  est une famille  $S-$plate de  faisceaux
coh\'erents de classe de Grothendieck $c$ et que $\Gamma$  est un
sous-faisceau localement libre de $\pr_{1*}(\Theta)$  induisant
au-dessus de chaque point un syst\`eme coh\'erent.
On
consid\`ere le morphisme de groupes ab\'eliens
$$\lambda_{\Lambda}:\widetilde{K}(\proj_2)\hfl{}{}
\Pic(S)$$ d\'efini par $(u,m)\mapsto \det \pr_{1!}(\Theta.u) 
\otimes (\det
\Gamma)^{\otimes m},$  o\`u $\Theta.u$  d\'esigne par abus le produit des
classes
$\Theta.\pr_2^*(u)$  dans $K(S\times \proj_2).$  Ce produit a un sens car
$\Theta$ a une r\'esolution  localement libre de longueur 1. Pour tout
faisceau inversible
$A
$ sur
$S$, consid\'erons la famille $\Lambda \otimes \pr_{1}^*(A)=(\Gamma
\otimes
\pr_{1}^*(A),
\pr_{1}^*(A)\otimes\Theta)$.  Il est alors clair que 
$$\lambda_{\Lambda \otimes \pr_{1}^*(A)}(u,m)= \lambda_{\Lambda}(u,m)
\otimes A^{\langle c,u \rangle + 2m}$$

Comme on l'a d\'ej\`a vu, l'espace de modules
$S_{4}$  se construit en quotientant un ouvert convenable ${\cal
G}$ d'une grassmannienne relative au-dessus d'un sch\'ema Quot  par
l'action d'un groupe r\'eductif.  Cet ouvert param\`etre une famille
universelle de syst\`emes coh\'erents et par application du crit\`ere de
descente de Kempf
 on peut alors construire  un homomorphisme  de groupes [16]
$$\lambda_{S_{4}}: (c,2)^{\bot}\hfl{}{}\Pic(S_{4})$$ 
caract\'eris\'e par la propri\'et\'e universelle suivante: pour toute famille
$\Lambda=(\Gamma,\Theta)$  param\'etr\'ee par $S$, on a 
$$\lambda_{\Lambda}(u,m)=f_{\Lambda}^*
(\lambda_{S_{4}}(u,m))$$
o\`u $f_{\Lambda}:S\hfl{}{}S_{4}$  est le morphisme
modulaire associ\'e \`a la famille.

\begin{prop} \label{61} 
Le morphisme canonique $$\lambda_{S_{4}}:(c,2)^{\bot}\hfl{}{}
\Pic(S_{4})$$ est un isomorphisme.

\end{prop}

\begin{proof}
D'apr\`es le th\'eor\`eme \ref{eclate}, le groupe de Picard de l'espace de
modules
$S_{4}$  est \'evidemment 
somme directe du groupe de Picard de
$M_{4}$  et du groupe engendr\'e par le fibr\'e inversible  $\hbox{\goth{e}}$
engendr\'e par le diviseur exceptionnel.  L'\'enonc\'e sera donc une
cons\'equence de la description du groupe de Picard de $M_{4}$ et
l'interpr\'etation du diviseur exceptionnel. Consid\'erons  une famille 
$\Lambda=(\Gamma,\Theta)$ de  syst\`emes coh\'erents semi-stables 
param\'etr\'ee par
$S$  et la famille $F$ de faisceaux semi-stables de classe de
Grothendieck
$2-4h^2$  associ\'ee, d\'efinie par l'extension canonique  sur $S \times
\proj_2$
$$0 \hfl{}{} \Gamma^{*}\otimes {\cal O}_{\proj_2}\hfl{}{} F(1) \hfl{}{}
{{\check{\Theta}}}\hfl{}{} 0$$ o\`u ${\check{\Theta}}=
\ul{\Ext}^1(\Theta,{\cal O}).$

\begin{lemme} \label{exceptionnel} Soit $\hbox{\goth{e}}$  le fibr\'e
inversible associ\'e au diviseur exceptionnel. Dans $\Pic(S_{4}),$ on a 
$$\hbox{\goth{e}}= -\lambda_{S_{4}}((1-h)^{3},0)$$
\end{lemme}
\begin{proof} Le support du diviseur exceptionnel est  
le ferm\'e   des syst\`emes coh\'erents semi-stables $(\Gamma,\Theta)$  tels que
$\Theta$  soit instable.  Un tel faisceau  a pour support une conique
singuli\`ere  et la filtration de Harder-Narasimhan de $\Theta$ 
s'\'ecrit
$$0\hfl{}{}{\cal O}_{\ell'}(3)\hfl{}{}\Theta \hfl{}{}{\cal
O}_{\ell''}(1)\hfl{}{}0$$
o\`u $\ell'$  et $\ell''$  sont deux droites  de $\proj_2.$
Donnons une description sch\'ematique de ce diviseur exceptionnel.
On peut trouver une famille $\Lambda$ param\'etr\'ee par une
vari\'et\'e lisse et connexe $S$,  telle que le morphisme modulaire
$f_{\Lambda}$ soit surjectif et telle que  le morphisme 
$$f_{\Lambda}^*:\Pic(S_{4})\hfl{}{} \Pic (S)$$
soit injectif. 
L'image r\'eciproque dans $S$ du diviseur exceptionnel  est alors 
d\'efini par l'id\'eal de Fitting du faisceau $R^{1}\pr_{1*}(\Theta(-3)).$ Ce
faisceau coh\'erent est de dimension homologique 1 et l'image
directe
$\pr_{1*}(\Theta (-3))$ est nulle.  En effet, soit $C\hfl{}{}
S$  le support sch\'ematique de $\Theta$: c'est une conique relative, et
$\Theta(-3)$  est un ${\cal O}_{C}-$module qu'on \'ecrit comme conoyau
d'un morphisme  
$0\hfl{}{}A\hfl{}{} B\hfl{}{} \Theta(-3) \hfl{}{} 0$
o\`u $B={\cal O}_{C}(-m)^{N}$  est un faisceau localement libre  sur $C$ 
suffisamment n\'egatif.  Alors par image directe on obtient une suite
exacte
$$0\hfl{}{} \pr_{1*}(\Theta(-3))\hfl{}{} R^{1}\pr_{1*}(A)\hfl{}{}
R^{1}\pr_{1*}(B)\hfl{}{} R^{1}\pr_{1*}(\Theta(-3))\hfl{}{}0$$
Les faisceaux $R^{1}\pr_{1*}(A)$  et $R^{1}\pr_{1*}(B)$  sont
localement libres.  Ceci prouve que l'image directe
$\pr_{1*}(\Theta(-3))$  est sans torsion; comme ce faisceau est
g\'en\'eriquement nul, il est identiquement nul. 
Ceci entra\^{\i}ne que  le sous-sch\'ema d\'efini par   l'id\'eal de Fitting de
$R^{1}\pr_{1*}(\Theta(-3))$ est
 le sch\'ema des z\'eros d'une section non identiquement nulle du fibr\'e
inversible
$-\det \pr_{1!}(\Theta(-3))=-\lambda_{\Lambda}((1-h)^{3},0).$ L'\'enonc\'e en
r\'esulte. 
\end{proof}

 Posons pour tout $u \in K(\proj_2), \check{u}=u^* \otimes
\omega_{\proj_2},$ o\`u $u\mapsto u^{*}$  d\'esigne l'homomorphisme de
conjugaison, qui envoie la classe d'un faisceau localement libre   sur
son dual \footnote{$(^4)$}{Cette conjugaison est un automorphisme
d'alg\`ebre qui envoie $h$  sur $h^*=-h-h^2$}.   Consid\'erons
 l'homomorphisme $\lambda_{\Theta}: K(\proj_{2})\hfl{}{} \Pic(S)$
associ\'e \`a la famille $\Theta$ (\cf [16], paragraphe 2.1), d\'efini par
$\lambda_{\Theta}(v)= \det \pr_{1!}(\Theta.u)$.
 On a par dualit\'e de Serre  pour tout $v \in K(\proj_2)$
$$\pr_{1!}(\check{\Theta}(v))=
-\pr_{1!}(\Theta(\check{v}))^*$$    par suite
$\lambda_{\check{\Theta}}(v)= \lambda_{\Theta}(\check{v}).$ Consid\'erons
d'autre part l'homomorphisme
$\lambda_{F}:K(\proj_2)\hfl{}{}
\Pic(S)$  associ\'e \`a la famille $F$, d\'efini par $\lambda_{F}(u)=
\det
\pr_{1!}(F(u))$ (\cf [16], paragraphe 2.1).  On  a
$$\lambda_{F}(u)=  (\det
\Gamma )^{\otimes-\chi(u(-1))}\otimes \lambda_{\Theta}(u^*(-2))
=\lambda_{\Lambda}(u^*(-2),-\chi(u(-1)))  $$
On est amen\'e \`a introduire l'homomorphisme de
groupes ab\'eliens $K(\proj_{2})\hfl{}{}
\widetilde{K}(\proj_2)$  d\'efini par $\xi:u\mapsto
(u^*(-2),-\chi(u(-1)).$

\begin{lemme} {\rm (i)} On a $\langle \xi(u),(c,2)\rangle =-\langle u,
2-4h^2 \rangle $

{\rm (ii)}  L'image de $\xi$  est le sous-groupe de
$\widetilde{K}(\proj_2)$  des couples $(v,m)$  tels que $\chi(v)+m=0.$

\end{lemme}
\begin{proof} C'est calcul facile reposant sur la dualit\'e de Serre, le
fait que dans $K(\proj_2)$  on ait $[{\cal O}_{\proj_2}(-1)]=1-h,$  et
la d\'efinition de $\xi.$ \end{proof}

Il en r\'esulte que l'image par $\xi$ de $(2-4h^2)^{\bot}$   est contenue
dans $(c,2)^{\bot}.$  D'apr\`es le
th\'eor\`eme
\ref{eclate},  on a un scindage 
$$\Pic(S_{4}) = \Pic(M_{4}) \oplus \entr \hbox{\goth{e}}$$  o\`u
$\hbox{\goth{e}}$  est la  classe du diviseur exceptionnel.
On a alors le diagramme commutatif  de suites
exactes:
$$\diagram{
0&\hfl{}{} &(2-4h^2)^{\bot}&\hfl{}{}&(c,2)^{\bot}& \hfl{\rho}{}& \entr&
\hfl{}{} 0 \cr
&&\vfl{\lambda_{M_{4}}}{}&&\vfl{\lambda_{S_{4}}}{}&&\vfl{}{}&&\cr
0&\hfl{}{}& \Pic(M_{4})& \hfl{}{} &\Pic(S_{4})&\hfl{}{}&
\entr&\hfl{}{} 0
\cr}$$
o\`u  $\lambda_{M_{4}}$ l'homomorphisme d\'efini dans [16]  et $\rho$  est 
induit par  l'homomorphisme $(v,m)\mapsto \chi(v)+m $ 
Cet homomorphisme est
surjectif, par exemple parce que
$(h^2,0)$  appartient \`a
$(c,2)^{\bot}$. Il est en fait plus habile de remarquer que 
$((1-h)^3,0)$  appartient aussi \`a
$(c,2)^{\bot}$  et que l'on a $\chi((1-h)^3)=1:$ ainsi, $((1-h)^3,0)$
r\'ealise un scindage de la premi\`ere suite exacte. La premi\`ere fl\`eche
verticale est un isomorphisme d'apr\`es le th\'eor\`eme de Dr\'ezet. 
D'apr\`es le lemme \ref{exceptionnel}, l'image de
$((1-h)^{3},0)$ par $\lambda_{S_{4}}$  est $-\hbox{\goth{e}}.$  Par 
suite l'homomorphisme $\lambda_{S_{4}}:(c,2)^{\bot}\hfl{}{}
\Pic(S_{4})$  est un isomorphisme. Ceci ach\`eve la d\'emonstration
de la proposition \ref{61}.\end{proof}

L'orthogonal   $(1-2h^2)^{\bot}$ dans $K(\proj_2)$ est engendr\'e par
$[{\cal O}_{\ell}(-1)]= h-h^2$  et $[{\cal
O}(1)]-h^2=(1-h)^{-1}-h^2=1+h.$ Remarquons que $$\xi(h-h^2)=(-h,1)
\  ;\
 \xi(1+h)=([{\cal O}(-3)]-h^2,0).$$ 
\begin{defi}
On pose
$$\eqalign{\hbox{\goth{d}}&=\lambda_{M_{4}}(-h+h^2)=
\lambda_{S_{4}}(h,-1)\cr
\hbox{\goth{k}}&=\lambda_{S_{4}}(h^2,0)\cr}$$ 
\end{defi}
Le fibr\'e inversible $\hbox{\goth{d}}$ ainsi d\'efini \`a isomorphisme
pr\`es  est l'image r\'eciproque du fibr\'e d\'eterminant de Donaldson ${\cal
D}$  sur $M_{4}.$

\begin{cor}  {\rm (i)} Le groupe de Picard $\Pic(S_{4})$ est un groupe
ab\'elien libre de rang 3, de base $\hbox{\goth{d}}, \hbox{\goth{k}},
\hbox{\goth{e}}.$  

{\rm (ii)} L'image du groupe de Picard de $M_{4}$ est \'egale  au
sous-groupe de $\Pic(S_{4})$ engendr\'e par les
fibr\'es d\'eterminant $\hbox{\goth{d}}$  et 
$\hbox{\goth{a}}=-\lambda_{M_{4}}(1+h)=\hbox{\goth{k}}+\hbox{\goth{e}}$.
\end{cor}

La classe $\hbox{\goth{k}}$  a une interpr\'etation g\'eom\'etrique
remarquable:
\begin{prop}  Consid\'erons le morphisme canonique
$\sigma:S_{4}\hfl{}{}{\cal C}_{2}= \proj_5.  $  La
classe
$\hbox{\goth{k}}$  est l'image r\'eciproque de  la classe  ${\cal
O}_{\proj_5}(1)$  par $\sigma.$
\end{prop}
\begin{proof}  On utilise la famille universelle
$\Lambda=(\Gamma,\Theta)$ de syst\`emes coh\'erents introduite dans la
d\'emonstration du lemme
\ref{exceptionnel}.  
  Alors une r\'esolution  localement libre
$$0\hfl{}{}A\hfl{}{}B\hfl{}{}\Theta\hfl{}{}0$$
fournit une section du fibr\'e inversible sur $S\times \proj_2$ d\'efini par 
$\det B \otimes (\det A)^{-1}= \lambda_{\Lambda}(h^2,0)\boxtimes {\cal
O}_{\proj_2}(2).$  Cette section peut \^etre vue comme un morphisme 
surjectif de fibr\'es sur $S $ 
$$H^0({\cal O}_{\proj_2}(2))^* \otimes {\cal O}_{S}\hfl{}{}
\lambda_{\Lambda}(h^2,0) .$$  On en d\'eduit par la propri\'et\'e universelle
de l'espace projectif un morphisme $\tau :S\hfl{}{} \proj_5$
 tel que $\tau^*({\cal O}(1))=\lambda_{\Lambda}(h^2,0).  $  Ce morphisme
se factorise \`a travers $\sigma$ suivant le diagramme
$$\diagram{
S&\hfl{f_{\Lambda}}{} &S_{4}\cr
&\sefl{}{\tau}&\vfl{}{\sigma}\cr
&&\proj_5
}$$
Puisque d'apr\`es le choix de la famille universelle  le morphisme induit
par $f_{\Lambda}$ sur les groupes de Picard est injectif, il en r\'esulte
que
$\sigma^*({\cal O}(1))=\hbox{\goth{k}}.$

\subsection{Quatre diviseurs de $S_{4}$}

Outre le diviseur exceptionnel $\hbox{\goth{e}}$  de $S_{4}$  on va
d\'ecrire trois autres diviseurs de $S_{4}$. 
\medskip
{\it Le bord $\partial S_{4}.$}
Le bord
$\partial S_{4}$  de
$S_{4}$  est constitu\'e des syst\`emes coh\'erents \`a point de base. 
On obtient ainsi  un diviseur irr\'eductible de $S_{4}$  qui n'est
autre que  l'image r\'eciproque  du bord $\partial M_4.$

\medskip
{\it Le diviseur $D_{1}.$} On consid\`ere  le sous-ensemble localement
ferm\'e des classes de syst\`emes coh\'erents stables $(\Gamma,\Theta)$ tels
que  le support $\hbox{\goth{c}}$ de  $\Theta$  soit lisse, et tels qu'il
existe une section $s$ de $\Gamma$  ayant deux z\'eros doubles.
Autrement dit, il existe une droite $\ell$ d'\'equation $z=0$  telle que
$\Gamma
\cap z^2H^0(\Theta(-2))\not=0.$
  On d\'esigne par $D_{1}$ son adh\'erence: c'est donc un diviseur
de
$S_{4}.$ D'apr\`es la proposition \ref{points-singuliers}, la quartique
de L\"uroth associ\'ee \`a un point de ce diviseur est singuli\`ere au point de
$\proj_{2}^*$  correspondant \`a la droite $\ell.$  Une telle quartique est
appel\'ee quartique de  L\"uroth singuli\`eres de type I.

\medskip
{\it Le diviseur $D_{2}.$}
Consid\'erons dans
$S_{4}$ l'image r\'eciproque du diviseur $\Delta$ des coniques
singuli\`eres  par le morphisme
$\sigma:S_{4}\hfl{}{}M_{4}.$ 

\begin{lemme} 
Le diviseur $\sigma^{-1}(\Delta)$ est r\'eduit.
\end{lemme}
 \begin{proof} En effet, il
r\'esulte de la suite exacte longue  d\'ecrite dans la
section 3  que le morphisme
$\sigma$  est une submersion sur l'ouvert $U$ des syst\`emes coh\'erents
$\Lambda=(\Gamma,\Theta)$ stables et tels que le faisceau $\Theta$
soit non singulier:
dans cet ouvert, l'image r\'eciproque de
$\Delta$  est lisse car  la conique associ\'ee en un point de ce
diviseur  est obligatoirement  r\'eduite, et donc $\Delta$ est lisse au
point d\'efini par cette conique. Il r\'esulte du calcul effectu\'e dans la
d\'emonstration du lemme \ref{fibres}  que le compl\'ementaire de
$U$ dans $S_{4}$ est de codimension
$\geq 2$; ceci entra\^{\i}ne que $\sigma^{-1}(\Delta)$  est r\'eduit. 
\end{proof}

Ce  diviseur  a deux
composantes irr\'eductibles: le diviseur
exceptionnel $\hbox{\goth{e}}$, qui ne rencontre pas le ferm\'e des
points strictement semi-stables. L'autre composante
$D_{2}$ est  constitu\'ee de  l'adh\'erence de l'ensemble  
 des classes d'isomorphisme de  syst\`emes coh\'erents  stables 
$\Lambda=(\Gamma,\Theta)$ tels que le support sch\'ematique de $\Theta$
soit une conique singuli\`ere, et que $\Theta$  soit semi-stable: la
description de ce ferm\'e $D_{2}$ sera d\'etaill\'ee dans la section 8; nous
verrons notamment qu'il est bien irr\'eductible.
\bigskip

\begin{prop} \label{calcul} Dans le groupe de Picard $\Pic(S_{4})$ on a,
en notation  additive, $$\eqalign{\partial S_{4}&= 5
\hbox{\goth{d}}-2\hbox{\goth{k}}-2\hbox{\goth{e}}\cr
D_{1}&= 12 \hbox{\goth{d}} \mod \hbox{\goth{e}}\cr
D_{2}&= 3\hbox{\goth{k}} - \hbox{\goth{e}}\cr}$$

\end{prop}

La premi\`ere assertion r\'esulte de la  formule qui donne le bord dans
le groupe de Picard de $M_{4}$ (\cf par exemple Le Potier  [17])  et du
fait que ce bord ne contient pas la sous-vari\'et\'e $\Sigma_{4}$  des
diviseurs sp\'eciaux.  La derni\`ere \'egalit\'e r\'esulte de la d\'efinition du
diviseur 
$D_{2}.$ La seconde assertion est plus difficile et r\'esulte
essentiellement de la formule de Porteous.  Comme elle n'est pas
indispensable pour la d\'emonstration du th\'eor\`eme principal, nous
n'aborderons pas ce calcul ici. Elle ne sera utilis\'ee que dans le
corollaire \ref{degreII}.
\bigskip

\subsection{Le diviseur $\beta^{-1}({\cal S}_{4})$}

On d\'esigne par ${\cal D}_{1}$ et ${\cal D}_{2}$ les images des
diviseurs
$D_{1}$ et $D_{2}$ de $S_{4}$ par le morphisme $\pi:S_{4}\hfl{}{}
M_{4}$.  Puisque  ces diviseurs $D_{1}$  et $D_{2}$ ne sont pas contenus dans
le diviseur exceptionnel, ${\cal D}_{1}$ et ${\cal D}_{2}$ sont des diviseurs
de
$M_{4}.$  Le diviseur ${\cal D}_{2}$ a l'interpr\'etation suivante, en
termes de droites de saut:
\begin{lemme}\label{diviseur2}
 Le diviseur  ${\cal D}_{2}$ de $M_{4}$ est  le  diviseur repr\'esentant les
faisceaux semi-stables qui ont au moins une droite de saut d'ordre $\geq 2.$
\end{lemme}
Ce lemme  r\'esulte de la
description d'un point g\'en\'eral de $D_{2}$ et sa d\'emonstration est
report\'ee \`a la fin de la section 8.

\begin{th} \label{multiplicites} Soit $\beta:M_{4}\hfl{}{} {\cal C}_{4}$  le
morphisme de Barth.  L'image r\'eciproque du diviseur ${\cal S}_{4}$  des
quartiques singuli\`eres de ${\cal C}_{4}$ est le diviseur
$$\beta^{-1}({\cal S}_{4})={\cal D}_{1}+2{\cal D}_{2}+3 \partial
M_{4}\eqno(6.1)\quad$$
\end{th}

\begin{proof}
Soit $\gamma= \beta \circ \pi:S_{4}\hfl{}{} {\cal C}_{4}$ le
morphisme de Barth. Il revient  au m\^eme de v\'erifier que 
$$\gamma^{-1}({\cal S}_{4})=D_{1}+2D_{2}+3\partial S_{4} \mod
\hbox{\goth{e}}. \eqno(6.2)\quad$$ D'apr\`es la proposition
\ref{points-singuliers}, on a
$\gamma(D_{1})\subset {\cal S}_4$ et le
corollaire \ref{singuliers-evidents} permet de voir
 que
$\gamma(D_{2})\subset {\cal S}_{4}$  (voir aussi section 8).
L'image du bord $\partial S_{4}$  est contenu dans le ferm\'e des
quartiques r\'eductibles, et l'image du diviseur exceptionnel est
constitu\'e de quartiques ayant au moins un point triple puisque c'est
l'image par
$\beta$ de la vari\'et\'e $\Sigma$ des faisceaux sp\'eciaux. L'image d'un
syst\`eme coh\'erent qui n'appartient \`a l'un de ces  diviseurs est une
quartique lisse d'apr\`es la proposition \ref{points-singuliers}.  On a donc
ensemblistement
$$\gamma^{-1}({\cal S}_{4})= D_{1}\cup D_{2} \cup \partial S_{4} \cup
\hbox{\goth{e}}$$
D'apr\`es la proposition \ref{calcul}, dans le groupe de Picard
$\Pic(S_{4})$  on a
$D_{1}=12
\hbox{\goth{d}}
\mod \hbox{\goth{e}}$, 
$D_{2}=3\hbox{\goth{k}} \mod \hbox{\goth{e}}$  et $\partial S_{4}=
5\hbox{\goth{d}}-2\hbox{\goth{k}} \mod \hbox{\goth{e}}$. Compte-tenu du
fait que l'hypersurface ${\cal S}_{4}$ est de degr\'e 27, on a
$\gamma^{-1}({\cal S}_{4}) = 27 \hbox{\goth{d}}$,  d'o\`u il r\'esulte
la formule attendue.
\end{proof}
La formule (6.1) r\'epond en partie \`a une question qui  nous a \'et\'e pos\'ee
il y a une dizaine d'ann\'ees par C. Peskine.
 Elle
montre que $\beta$  est transverse \`a ${\cal S}_{4}$  sur un ouvert de
${\cal D}_{1}:$  pour obtenir le
th\'eor\`eme \ref{main} il suffirait  donc de prouver que  $\beta|_{{\cal
D}_{1}}$  est g\'en\'eriquement injective.  Nous ne savons malheureusement
pas traiter cette question.  Par contre, l'\'etude de $\beta|_{{\cal
D}_{2}}$  est plus facile; la formule ci-dessus dit  que $\beta$  
n'est pas transverse \`a
${\cal S}_{4}$  sur un ouvert ${\cal D}_{2},$ ce qui complique  la
d\'emonstration du th\'eor\`eme \ref{main}  dans le cas $n=4$. 

\section{G\'eom\'etrie des quartiques singuli\`eres}
On pose  $V= H^{0}({\cal O}_{\proj_2}(1))^*$  de sorte  que le plan
projectif consid\'er\'e est l'espace des droites vectorielles de $V.$ 
On d\'esigne par ${\cal C}_{n}=\proj(H^0({\cal
O}_{{\proj}_2}(n)))$  le syst\`eme lin\'eaire  des courbes planes de degr\'e
$n$, et par ${\cal S}_{n} \subset {\cal C}_{n} $  l'hypersurface des
courbes singuli\`eres.  Si $\hbox{\goth{q}}$  est une courbe de degr\'e
$n,$  repr\'esent\'ee  par un polyn\^ome homog\`ene $f\in S^{n}V^*$ la courbe
polaire
$\hbox{\goth{q'}}(a)$ d'un point
$a
\in
\proj_2$  est la  courbe de degr\'e
$n-1$ d\'efinie par $x\mapsto \partial_{x}f.a'$, o\`u $\delta_{x}f$ 
d\'esigne la diff\'erentielle de $f$ au point $x$ et $a'$  est un  vecteur
non nul de
$V$ repr\'esentant le point
$a$ Cette courbe polaire recoupe la courbe
$\hbox{\goth{q}}$  aux points de contact
$x$
 des tangentes \`a $\hbox{\goth{q}}$ issues de $a.$  On s'int\'eresse ici
aux quartiques singuli\`eres.
\bigskip
\subsection{Conique associ\'ee \`a une quartique singuli\`ere.}

On se propose de montrer que dans le cas d'une quadrique singuli\`ere
suffisamment g\'en\'erale, les points de contacts des tangentes issues du
point singulier appartiennent \`a une conique $\hbox{\goth{p}}$  bien
d\'etermin\'ee.

\begin{prop}
Soit $\hbox{\goth{q}}$  une quartique plane singuli\`ere en un
seul point double
$O$, et $\hbox{\goth{q}}'$  la cubique polaire de $O$. On consid\`ere le nombre
d'intersection $(\hbox{\goth{q},\hbox{\goth{q}}}')_{O}$.

{\rm (i)}  On a $(\hbox{\goth{q}},\hbox{\goth{q}}')_{O}\geq 6,$ 
l'\'egalit\'e n'ayant lieu que si $\hbox{\goth{q}}'$ int\`egre.

{\rm (ii)}  Si $(\hbox{\goth{q}},\hbox{\goth{q}}')_{O}=6$ 
l'intersection $\hbox{\goth{q}}\cap \hbox{\goth{q}}'$ est la r\'eunion
  d'un  sous-sch\'ema de longueur $6$ concentr\'e en
$O,$ contenu dans le c\^one tangent en $O$ \`a $\hbox{\goth{q}}$ et 
d'un sous-sch\'ema de longueur 6  contenu dans une conique
$\hbox{\goth{p}}$ ne passant pas par $O.$

\end{prop}
\begin{proof}
L'hypoth\`ese assure que la quartique $\hbox{\goth{q}}$ est int\`egre, et
donc le nombre d'intersection 
$(\hbox{\goth{q}},\hbox{\goth{q}}')_{O}$  a un sens. On peut supposer
que le point singulier
$O$ est l'origine dans le plan affine $A$ d\'efini par le 
compl\'ementaire  de la droite $\Delta$ d'\'equation
$t=0,$  o\`u $t \in H^{0}({\cal O}_{\proj_2}(1)),$ qu'on prend pour
droite de l'infini. On d\'esigne par ${\hbox{\goth{m}}}_{O}$   l'id\'eal
du point $O.$ La quartique
$\hbox{\goth{q}}$  est le sch\'ema des z\'eros d'une section $f \in
H^0({\cal O}_{\proj_2}(4))$  s'\'ecrivant
$$f=t^{2}f_{2}+tf_{3}+f_{4}\eqno(7.1)\quad$$  o\`u $f_{i}\in
H^{0}(\hbox{\goth{m}}^{i}_{O}(i))\simeq H^{0}({\cal O}_{\Delta}(i))$ 
est un polyn\^ome homog\`ene de degr\'e
$i$, 
 de sorte que  la cubique $\hbox{\goth{q}}'$  polaire  de $O$
est le sch\'ema des z\'eros de la forme cubique
$$g=2tf_{2}+f_{3}\eqno (7.2)\quad$$
Il r\'esulte des formules (7.1) et (7.2)  ci-dessus que $\hbox{\goth{q}}'$ 
et $\hbox{\goth{q}}$ ont m\^eme c\^one $C$ \`a l'origine,  d\'efini par  la forme
quadratique
$f_{2}.$ Le nombre d'intersection \`a l'origine
$(\hbox{\goth{q}},\hbox{\goth{q}}')_{O}$  est donn\'e par
$$(\hbox{\goth{q}},\hbox{\goth{q}}')_{O}=(tf_{3}+2f_{4},2tf_{2}+f_{3})_{O}$$
Il en r\'esulte que $(\hbox{\goth{q}},\hbox{\goth{q}}')_{O}\geq 6$  avec
\'egalit\'e si et seulement si les formes binaires $f_{2}$  et $f_{3}$  n'ont
pas de z\'ero commun sur $\Delta.$ Il revient au m\^eme de demander que le
c\^one
$C$  et la cubique $\hbox{\goth{q}}'$  n'ont pas de composante commune,
ou encore que la  cubique  $\hbox{\goth{q}}'$ est int\`egre, ce qui
d\'emontre (i) .

 Supposons que $(\hbox{\goth{q}},\hbox{\goth{q}}')_{O}=6$.
Sur la  droite de l'infini $\Delta$
les formes binaires  
$f_{2}
$ et $f_{3}$  n'ont pas de z\'ero en commun. 
Le complexe de Koszul associ\'e \`a $(f_3,f_2)$ montre qu'il
existe  une forme  lin\'eaire $\phi$ et une forme quadratique
$\psi$ uniques
 telles que  
$f_{4}=\phi f_{3} + \psi f_{2}.$
On a alors $$\eqalignno{f&= (t^{2}+ \psi)f_{2}
+(t+\phi)f_{3}&\cr
&= (t+\phi)g-(t^2+2t\phi-\psi )f_{2}& (7.3)\quad  \cr
}$$
La conique $\hbox{\goth{p}}$ d'\'equation $t^2+2t\phi-\psi =0$  ne passe
pas par l'origine. Le cycle  intersection  de la quartique
$\hbox{\goth{q}}$
  et de la cubique $\hbox{\goth{q}}'$  est  donn\'e par
$$(\hbox{\goth{q}},\hbox{\goth{q}}')=(g,f_{2})_{O}O +
\sum_{m}(\hbox{\goth{q}}',\hbox{\goth{p}})_{m}m
$$
ce qui d\'emontre l'assertion (ii).

\subsection{Le morphisme $\hbox{\goth{q}}\mapsto
(\hbox{\goth{q}}',\hbox{\goth{p}})$}

On d\'esigne par ${\cal S}_{4}$  l'hypersurface de ${\cal C}_{4}$  des
quartiques singuli\`eres, et par ${\cal S}^{*}_{4}$  l'ouvert de ${\cal
S}_{4}$  des courbes  singuli\`eres $\hbox{\goth{q}}$ ayant au plus un
point singulier ordinaire
$O$,  satisfaisant \`a la condition $(\hbox{\goth{q}},\hbox{\goth{q}}')_{O}=6,$
autrement dit telles que la polaire $\hbox{\goth{q}}'$  soit int\`egre.

On d\'esigne de m\^eme par ${\cal S}_{3}$  l'hypersurface de ${\cal
C}_{3}$  des cubiques singuli\`eres, et par ${\cal S}_{3}^*$ l'ouvert de
${\cal S}_{3}$ des cubiques qui ont exactement un point double ordinaire.

On dispose d'un morphisme 
${\cal S}_{4}^*\hfl{}{} \proj_2$
qui associe \`a la quartique $\hbox{\goth{q}}$  son  unique  point
singulier $O$. Consid\'erons en effet l'hypersurface  ${\cal S}_{4}$ des
quartiques singuli\`eres.
 C'est l'image de la 
sous-vari\'et\'e $\widetilde{{\cal S}}_{4}$ lisse  de codimension
$3$  de ${\cal C}_{4}\times \proj_2$  des couples $(\hbox{\goth{q}}, m)
\in {\cal C}_{4}\times
\proj_2$  tels que 
${\rm j}_{m}f=0,$
o\`u $f$  d\'esigne l'\'equation de $\hbox{\goth{q}},$ et ${\rm j}_{m}f$
le jet d'ordre un de $f$ au point $m$,  par  la projection
$(\hbox{\goth{q}},m)\mapsto \hbox{\goth{q}}$. Cette projection  a alors une
section au-dessus de l'ouvert ${\cal S}_{4}^*.$ La compos\'ee de cette
section avec la projection $\widetilde{{\cal S}}_{4} \hfl{}{} \proj_2$ 
fournit le morphisme annonc\'e.  La m\^eme remarque vaut pour l'ouvert des
cubiques int\`egres \`a point double ordinaire: on dispose d'un morphisme ${\cal
S}_{3}^*\hfl{}{} \proj_2.$

\begin{prop}\label{produit}
Le morphisme $${\cal S}^{*}_{4}\hfl{}{} {\cal S}^*_{3}\times {\cal
C}_{2}$$  qui associe \`a
$\hbox{\goth{q}}$  le couple $(\hbox{\goth{q}}',\hbox{\goth{p}})$ est
un isomorphisme  sur un ouvert $U$ de la vari\'et\'e des couples form\'es
d'une cubique irr\'eductible singuli\`ere en un point double ordinaire
$O$  et d'une conique ne passant pas par $O.$

\end{prop}

\begin{proof}
Les deux membres sont des ouverts de fibr\'es localement triviaux en
espaces projectifs de  dimension 11  au-dessus de 
$\proj_2$, et le morphisme ci-dessus est compatible avec cette
projection, de sorte qu'on peut fixer le point singulier $O$, et
choisir  la droite de l'infini d'\'equation $t=0$ ne passant pas par $O.$  
 Toute  cubique $\hbox{\goth{q}}'$  n'ayant qu'un point
double ordinaire  est le sch\'ema des z\'eros d'une section de la forme
$g=2tf_{2}+f_{3}$ o\`u $f_{i}\in
H^{0}(\hbox{\goth{m}}^{i}_{O}(i))=H^{0}({\cal O}_{\Delta}(i))$, et o\`u
$f_{2}$  est une forme quadratique non d\'eg\'en\'er\'ee.
Toute conique
$\hbox{\goth{p}}$ du plan ne passant pas par l'origine $O$ a pour
\'equation 
$$t^2+2\phi t -\psi=0$$ 
o\`u $\phi \in H^{0}(\hbox{\goth{m}}_{O}(1)) $  et $\psi \in
H^{0}(\hbox{\goth{m}}^{2}_{O}(2))$, de sorte que le morphisme inverse 
 est donn\'e par $(\hbox{\goth{q}}',\hbox{\goth{p}})\mapsto 
\hbox{\goth{q}}$  o\`u la quartique $\hbox{\goth{q}}$ a pour \'equation
$f=0$, o\`u $f=t^2f_{2}+tf_{3}+f_{4},$ et o\`u $f_{4}$  se calcule par la
formule
$$f_{4}=\psi f_{2}+ \phi f_3.$$
En g\'en\'eral, si $\phi$  et $\psi$  sont quelconques, 
 la quartique $\hbox{\goth{q}}$ ainsi construite  peut \^etre singuli\`ere en
l'un des points  d'intersection de $\hbox{\goth{q}}'$  et $\hbox{\goth{p}}$,
par exemple si $f_{4}=0.$  Pour \^etre s\^ur  que $O$  est le seul point
singulier,  il faut  \'ecrire que $\hbox{\goth{q}}$  est lisse en ces points,
ce qui est une condition ouverte.
\end{proof}

\subsection{Calcul de la d\'eriv\'ee du morphisme $\hbox{\goth{q}}\mapsto
\hbox{\goth{p}}$}

 Soit
$t$  une forme lin\'eaire sur
$V,$ fix\'ee une fois pour toute, consid\'er\'ee comme \'equation de la droite
de l'infini
$\Delta$ et
$\xi$  un vecteur tel que $\langle t,\xi \rangle =1.$ Le plan 
projectif contient le plan affine $A\subset V$  d\'efini par
l'\'equation
$\langle t,x \rangle=1. $ On \'etudie les
quartiques singuli\`eres au voisinage d'un point $O\in A$  fix\'e.

   Tout vecteur de
$V$  s'\'ecrit
$x= \xi t + u$  avec $u \in W= \ker t,$
de sorte que l'\'equation $f$  de $\hbox{\goth{q}}$  s'\'ecrit 
$$f(\xi t + u)= \sum_{i}t^{4-i}f_{i}(\xi)(u) $$  o\`u $i!f_{i}(\xi)=
\partial^{i}_{\xi}f$  est la diff\'erentielle i-i\`eme de $f$
 en $\xi;$ c'est un polyn\^ome homog\`ene de degr\'e
$i$  sur
$W$ et $\xi \mapsto f_{i}(\xi)$  est homog\`ene de degr\'e $4-i$  en
$\xi.$ Dire que
$\hbox{\goth{q}}$ est singuli\`ere en $\xi$  signifie que
$f_{0}(\xi)=f_{1}(\xi)=0.$ La conique
$\hbox{\goth{p}}$ associ\'ee a un sens si
$f_{2}(\xi) $ et
$f_{3}(\xi)$ n'ont pas de z\'ero en commun et  a   pour \'equation
$$t^2+2t\phi -\psi=0$$ o\`u la forme lin\'eaire $\phi$ sur $W$ et la
forme quadratique $\psi$ sur $W$ sont d\'efinis par la condition  de
Koszul  dans $H^0({\cal O}_{\Delta}(4))$
$$f_{3}(\xi)\phi + f_{2}(\xi)\psi=f_{4}(\xi) \eqno(7.4)\quad$$

Soit $O$  un point fix\'e dans le plan affine $A.$  Pour tout $g \in
H^0({\cal O}(4))$, on \'ecrit  comme pour $f$
$$g(\xi t +u)= \sum_{i}t^{4-i}g_{i}(\xi)(u) \mod f.$$
L'espace vectoriel tangent \`a $\widetilde{\cal S}_{4}$ au point
$(\hbox{\goth{q}},O)$  est donn\'e par les couples
$(g,\xi)) \in H^0({\cal O}\hbox{\goth{q}}(4)) \times
W$  tels que $$\eqalign{g(O)&=0\cr
g_{1}(O)+ \partial_{O}f_{1}.\xi&=0\cr}$$
Remarquons que  du fait que $f_{1}(z)= \partial_{z}f$,
$\partial_{O}f_{1}$  est aussi le Hessien $ H_{O}f=2 f_{2}(O)$ vu
comme application lin\'eaire $W\hfl{}{} W^*.$
Par d\'erivation de 7.4 au point $(\hbox{\goth{q}},O)$, on
obtient
$$f_{3}(O)\dot{\phi} + f_{2}(O)\dot{\psi}+[g_{3}(O)
+\partial_{O} f_{3}.\xi]\phi+ [g_{2}(O)+ \partial_{O} f_{2}.\xi]\psi=
g_{4}(O)$$
ce qui d\'etermine $\dot\phi $
et $\dot \psi$  en fonction de $g$  et $\xi.$ 
Si on identifie $H^{0}({\cal O}_{\Delta}(i))= S^{i}W^*$  avec  
son image dans $\Hom(W,S^{i-1}W^*) $ on peut \'ecrire
$\partial_{O}f_i.\xi=
{1\over {i!}}\partial_{O}^{i+1}f(\xi)=(i+1) f_{i+1}(O)(\xi).$
La quartique $\hbox{\goth{q}}$  a en $O$ un point double
ordinaire;  d\'esignons par $H_{O}f= 2 f_{2}(O)$  le Hessien de $f$
en
$O,$  vu comme application lin\'eaire inversible $W\hfl{}{} W^*.$  
Les \'equations 
$$\eqalign{
H_{O}f (\xi) &=- g_{1}(O)\cr
f_{3}(O)\dot{\phi} + f_{2}(O)\dot{\psi}&=g_{4}(O)
-[g_{3}(O) +4f_{4}(O).\xi]\phi- [g_{2}(O)+
3 f_{3}(O).\xi]\psi \cr}
$$
permettent de calculer $(\xi,\dot{\phi},\dot{\psi})$ et donc  l'image du
vecteur tangent
$g$  dans $H^0({\cal O}_{\hbox{\goth{p}}}(2)).$

\subsection{Le diviseur ${\cal L}_{2}^*$  de ${\cal S}_{4}^*$}

On s'int\'eresse  au diviseur not\'e ${\cal L}_{2}^*$ de 
l'ouvert ${\cal S}_{4}^*$ des quartiques singuli\`eres telles que la conique
associ\'ee $\hbox{\goth{p}}$ soit singuli\`ere.  Le c\^one tangent $C$ au
point singulier
$O$ d'une quartique
$\hbox{\goth{q}}\in {\cal S}_{4}^*$  coupe la quartique suivant deux
points triples  et l'on a d'apr\`es la formule (7.3) 
$$(C,\hbox{\goth{q})}- (C,\hbox{\goth{q'}})_{O}O = (t+\phi,f_{2}) $$
l'intersection r\'esiduelle est donc  port\'ee par la droite d'\'equation
$t+\phi=0,$  polaire de $O$ par rapport \`a la conique
$\hbox{\goth{p}}$  associ\'ee \`a $O$.  

\begin{prop}
Soit $\hbox{\goth{q}}$  une quartique de ${\cal S}_{4}^*.$  Les
conditions suivantes sont \'equivalentes:

{\rm (i)}  La droite d\'efinie   par l'intersection r\'esiduelle $(C,
\hbox{\goth{q}})- (C,\hbox{\goth{q}})_{0}O$ est tangente \`a la quartique
$\hbox{\goth{q}}.$

{\rm (ii)}  La conique $\hbox{\goth{p}}$ associ\'ee \`a $\hbox{\goth{q}}$ 
est singuli\`ere.

Si ces conditions sont satisfaites, la conique $\hbox{\goth{p}}$ a un
seul point singulier qui est exactement le point de tangence de la
droite port\'ee par l'intersection  r\'esiduelle ci-dessus.

\end{prop}

\begin{proof}
On utilise la formule (7.3). On a vu ci-dessus
$$(f,t+\phi)=(t+\phi,\phi^2+\psi)+ (t+\phi,f_2)$$
L'intersection r\'esiduelle $(t+\phi,\phi^2+\psi)$  est l'intersection
du c\^one quadratique issu de $O$  d\'efini par $\phi^2+\psi \in
H^{0}(\hbox{\goth{m}}^{2}_{O}(2))$  et de la droite d'\'equation
$t+\phi=0,$  qui ne passe pas par $O.$  Cette intersection est un
point double si et seulement si  le discriminant
$\Delta(\phi^2+\psi)$ est nul.  Ce discriminant est aussi celui de
la forme quadratique d\'efinie  sur $V$  par $t^2+2t\phi-\psi$ ce qui
conduit \`a l'\'equivalence de (i) et (ii). Puisque la quartique est
irr\'eductible, le nombre d'intersection $(f,t+\phi)$  est fini, 
et  par suite la forme quadratique
$\phi^2+\psi$  n'est pas nulle. Alors le rang de la conique
$\hbox{\goth{p}}$  est  $2,$  et le noyau de la forme quadratique
correspond au point de tangence de la droite  d'\'equation  $t+\phi=0.$
\end{proof}
\goodbreak
\section{Restriction du morphisme de Barth \`a $D_{2}$ }

Donnons d'abord une description d\'etaill\'ee du diviseur $D_{2}.$
Consid\'erons un syst\`eme coh\'erent $(\Gamma,\Theta)$  repr\'esentant un point
stable de $D_{2}.$ Le
faisceau $\Theta$ sous-jacent n'est jamais stable: il 
contient au moins un sous-faisceau coh\'erent
$\Theta'$  de classe de Grothendieck ${c\over 2}$
et par stabilit\'e, le sous-espace $\Gamma$  ne rencontre pas $H^0(\Theta').$
Il s'\'ecrit  donc  comme extension 
$$0\hfl{}{}(0,\Theta')\hfl{}{} (\Gamma,\Theta)
\hfl{}{}(\Gamma'',\Theta'')\hfl{}{}0\eqno(8.1)\quad$$  o\`u $\Theta'$  et
$\Theta''$  sont des faisceaux purs de classe
${c\over 2}$; ce sont  donc des faisceaux inversibles de degr\'e $2$ sur
des droites $\ell'$ et $\ell''$ respectivement. Le syst\`eme coh\'erent
$(\Gamma'',\Theta'')$  d\'efinit un point de l'espace de modules
$\Syst({c\over 2},2).$ Cet espace de modules est tr\`es facile \`a d\'ecrire:
\begin{lemme} \label{vieux} On a un isomorphisme canonique
$$\syst({c\over 2},2)\simeq
\Hilb^{2}(\proj_2)$$
\end{lemme} 

\begin{proof}
Etant donn\'e un syst\`eme coh\'erent semi-stable
$\Lambda=(\Gamma,\Theta)$ de classe de Grothendieck ${c\over 2}$, 
le faisceau $\Theta$ sous-jacent  est isomorphe \`a un faisceau inversible
${\cal O}_{\ell}(2)$  sur une droite $\ell$ de
$\proj_2$; \`a un tel syst\`eme lin\'eaire $(\Gamma,{\cal O}_{\ell}(2))$ on
associe
 le sous-sch\'ema de $\ell $ des z\'eros
des quadriques singuli\`eres qui appartiennent \`a $\Gamma.$ On obtient ainsi
l'isomorphisme attendu.
\end{proof}

R\'eciproquement, consid\'erons l'ensemble des classes d'isomorphisme
d'extensions (8.1);  quand $\Lambda'=(0,\Theta')$ et
$\Lambda''=(\Gamma'',\Theta'')$ sont fix\'ees, ces extensions sont
class\'ees \`a isomorphisme pr\`es  par l'espace projectif
 $\proj(\Ext^1(\Lambda'',\Lambda'))$  des droites de
$\Ext^1(\Lambda'',\Lambda').$
 De la  suite exacte
$$0\hfl{}{}\Hom(\Theta'',\Theta')\hfl{}{}
\Hom(\Gamma'',H^0(\Theta'))\hfl{}{}
\Ext^1(\Lambda'',\Lambda')\hfl{}{}
\Ext^1(\Theta'',\Theta')\hfl{}{}0 \eqno(8.2)\quad $$ il d\'ecoule  que
l'espace vectoriel $\Ext^1(\Lambda'',\Lambda')$ est de dimension $7.$ 
Consid\'erons maintenant les syst\`emes coh\'erents universels $\Lambda'$ et
$\Lambda''$ param\'etr\'es par $S=\proj_2^* \times
\Hilb^{2}(\proj_2),$  et le fibr\'e vectoriel de rang $7$  sur $S$ d\'efini
par ${\cal V}=\ul{\Ext}_{\pr_{1}}^1(\Lambda'',\Lambda').$
Le fibr\'e en espaces projectif $\proj({\cal V})$ de rang 6 au-dessus de
$ \proj_{2}^*\times\Hilb^2(\proj_{2})$  associ\'e \`a ce  fibr\'e vectoriel
param\`etre alors toutes les extensions (8.1).

\begin{lemme}\label{instables} Soit $I$   le ferm\'e de $ \proj({\cal
V})$  des points
 d\'efinissant  un syst\`eme coh\'erent instable.

 {\rm (i)} L'image de $I$ est contenu dans le ferm\'e des couples
$(\Lambda',\Lambda'')$ tels que le syst\`eme lin\'eaire $\Lambda''=
(\Gamma'',\Theta'')$  ait un point de base.  

{\rm (ii)} Le ferm\'e $I$  est de dimension $4$  et la projection
$I\hfl{}{} \proj_2^*\times\Hilb^2(\proj_2)$ est injective.  

En
particulier une  fibre g\'en\'erale de $ \proj({\cal V})\hfl{}{}
\proj_{2}^*\times\Hilb^2(\proj_{2})$  ne rencontre pas le ferm\'e des
points instables.
 \end{lemme}
\begin{proof} 
Si le syst\`eme coh\'erent $(\Gamma,\Theta)$ est 
instable, il existe un sous-faisceau $F$ de $\Theta$ de
multiplicit\'e 1  tel que
$\Gamma\subset H^0(F)$. Alors $F$  n'est pas contenu
dans 
$\Theta'$ et la projection $F\hfl{}{} \Theta''$ est non nulle; par
suite c'est un plongement. Si c'est un isomorphisme, l'extension
sous-jacente
$$0\hfl{}{} \Theta'\hfl{}{}\Theta\hfl{}{}\Theta''\hfl{}{} 0
\eqno(8.3)\quad$$
serait scindable: le faisceau $F$ serait le graphe d'un morphisme
$a:\Theta''\hfl{}{} \Theta'$ et le sous-espace
$\Gamma$  serait l'image de $\Gamma''$ par  l'application lin\'eaire
$a : H^0(\Theta'')\hfl{}{} H^0(\Theta').$  Mais en vertu de la suite
exacte (8.2), ceci impliquerait que la classe dans
$\Ext^1(\Lambda'',\Lambda') $ de l'extension (8.1)  serait nulle, ce
qui est contraire \`a l'hypoth\`ese.
Il en r\'esulte qu'obligatoirement le morphisme $F\hfl{}{}\Theta''$ 
n'est pas surjectif et a pour image le sous-faisceau $F''$ de
$\Theta''$  des sections qui s'annulent en un point $\alpha \in
\ell''.$  Ceci prouve que le syst\`eme lin\'eaire $\Lambda''$  a un point
de base ce qui d\'emontre  (i).

 Reste \`a
trouver dans ce cas quels sont les points instables dans la fibre
au-dessus de $(\Lambda',\Lambda'').$
Supposons d'abord  $\Theta'\not=\Theta''$. Montrons que dans ce
cas l'extension sous-jacente n'est pas scindable.
Si cette extension  \'etait
scindable,  un tel  scindage fournirait  une  projection
$\Theta\hfl{}{} \Theta'$  de noyau $\Theta''$  et la  restriction 
$F\hfl{}{}
\Theta'$ serait  nulle; dans la somme directe
$H^0(\Theta)=H^0(\Theta')\oplus H^0(\Theta'')$ le sous-espace
$\Gamma$  proviendrait du sous-espace vectoriel $\Gamma''$ de
$\Theta''.$
 Mais alors  ceci signifierait  que la suite exacte de syst\`emes
coh\'erents (8.1)  serait  elle aussi scindable, ce qui est exclus par
hypoth\`ese.   Le support de $\Theta'$  est une droite $\ell'$ d'\'equation
$x=0$. Le faisceau 
$\Theta$  est un faisceau inversible sur la conique $\hbox{\goth{c}}$
support de
$\Theta$ ;   on a
$\Theta|_{\ell'}={\cal O}_{\ell'}(3),$  et la suite exacte 
$0\hfl{}{} \Theta''(-1)\hfl{x}{}
\Theta\hfl{}{}\Theta|_{\ell'}\hfl{}{}0.$ Le morphisme induit 
 $F\hfl{}{} \Theta|_{\ell'}$ est nul.  Donc $\Gamma=\Gamma ''
=H^0(\Theta''(-1))$  et ce sous-espace vectoriel s'identifie au syst\`eme
lin\'eaire des sections de
$\Theta''$  qui sont divisible par $x$  c'est-\`a-dire des sections qui
s'annulent au point singulier  $O$  de  $\hbox{\goth{c}}.$  Dans ce
cas,  dans la fibre, il n'y a donc qu'un point instable: le sous-espace
$\Gamma$  est le noyau du morphisme de restriction
$H^0(\Theta)\hfl{}{}H^0(\Theta|_{\ell'}).$

Supposons maintenant $\Theta'=\Theta'';$  le support de ces faisceaux
est une droite $\ell.$  Montrons qu'alors l'extension (8.2)
est scindable. On sait que $F$ isomorphe
\`a
${\cal O}_{\ell}(1).$  L'extension induite
$$0\hfl{}{} {\cal O}_{\ell}(2)\hfl{}{}\Theta_{1}\hfl{}{} F \hfl{}{}  0$$
est alors scindable; ces extensions sont class\'ees par 
$\Ext^1(F,\Theta')\simeq H^0({\cal O}_{\ell}(2))$, et le morphisme
$\Ext^1(\Theta'',\Theta')\hfl{}{} \Ext^1(F,\Theta') $ induite par
l'inclusion $F\hookrightarrow \Theta''$ s'identifie \`a l'application
lin\'eaire $H^0({\cal O}_{\ell}(1))\hfl{}{} H^0({\cal O}_{\ell}(2))$.
L'inclusion
$F\hookrightarrow \Theta''$ est donn\'ee par
une section non nulle de $H^0({\cal O}_{\ell}(1))$  et l'image de la
classe de l'extension correspondant \`a
$\Theta$  est  nulle. Donc l'extension (8.3)  est elle aussi
scindable. Si on  se donne un tel scindage,  la classe de l'extension
(8.1)  a un repr\'esentant dans l'espace vectoriel 
$\Hom(\Gamma'',H^0(F'))$ qui se calcule en \'ecrivant que  dans le somme
directe $H^0(\Theta)=H^0(\Theta')\oplus H^0(\Theta'')$
le sous-espace $\Gamma$  est le graphe d'une application lin\'eaire
$\alpha:
\Gamma''\hfl{}{} H^0(\Theta')$. 
Le faisceau $F$  est
l'image d'un morphisme de faisceaux $a: F''\hfl{}{}\Theta'$, et  le
sous-espace
$\Gamma$ est l'image de
$\Gamma''=H^0(F'')$  par l'application lin\'eaire 
$ H^0(a): \Gamma''=H^0(F'')\hfl{}{}H^0(\Theta')$; l'application
lin\'eaire $\alpha=H^0(a)$ est donc le repr\'esentant cherch\'e. Puisque
$\Hom(F'',\Theta') \simeq H^0({\cal O}_{\ell}(1))$   l'application
lin\'eaire compos\'ee  $H^0({\cal
O}_{\ell}(1))\hfl{}{}\Hom(\Gamma'',H^0(\Theta'))
\hfl{}{} \Ext^1(\Lambda'',\Lambda')$  est de rang 1
et  son image est une  droite vectorielle: cette classe 
d\'efinit un seul point  $\proj(\Ext^1(\Lambda'',\Lambda'))$; par suite  
    dans la fibre il n'y qu'un seul point
instable.

Il r\'esulte de cette description que se donner un point de $I$  revient \`a
un donner un couple $(\hbox{\goth{c}},O)$  form\'e d'une conique
singuli\`ere et d'un point singulier $O$  de cette conique.
   Dans $\proj({\cal V})$ le
ferm\'e $I$ des syst\`emes coh\'erents instables est donc de dimension 4 et ne
rencontre pas les fibres  au-dessus des points $(\Lambda',\Lambda'')$ 
telles que
$\Lambda''$  soit sans point de base.  Ceci ach\`eve la d\'emonstration.
\end{proof}

L'ouvert $\proj({\cal V})\setminus I$  param\`etre une famille de syst\`eme
coh\'erents semi-stables. Rappelons que les points strictement
semi-stables de $S_{4}$ constituent un ferm\'e de codimension 5. La
propri\'et\'e de module grossier de
$S_{4}$  fournit 
 un morphisme dominant $\proj({\cal V})\setminus I \hfl{}{} D_{2}.$ Il
en r\'esulte en particulier que le diviseur $D_2$ est bien irr\'eductible.
L'\'enonc\'e qui suit entra\^{\i}ne en fait que ce morphisme est birationnel.
  Par 
composition avec le morphisme de Barth on obtient un morphisme
$$\proj({\cal V})\setminus I
\hfl{}{} {\cal C}_{4}$$ not\'e encore $\gamma.$
On consid\`ere l'ouvert  ${\cal S}_{4}^*$ et  dans cet
ouvert le diviseur
${\cal L}_{2}^*$  introduit  dans la section pr\'ec\'edente.
\begin{th} \label{cle}
{\rm (i)}  Si $s\in \proj({\cal V})\setminus I$, la quartique
$\hbox{\goth{q}}_{s}$ associ\'ee \`a $s$ appartient
\`a ${\cal S}_{4}.$ 

{\rm (ii)} L'image de $\gamma: \proj({\cal V})\setminus I\hfl{}{} {\cal
S}_{4}$  rencontre l'ouvert
${\cal S}_{4}^*.$  Au-dessus de cet ouvert, le morphisme $\gamma$  est un
plongement qui a pour image le diviseur  ${\cal L}_{2}^*.$

\end{th}
Il en r\'esulte en particulier que le morphisme $\gamma:D_{2}\hfl{}{}
{\cal C}_{4}$  est g\'en\'eriquement injectif. En plus, cet \'enonc\'e permet
de d\'ecider qu'une quartique singuli\`ere suffisamment g\'en\'erale est ou non
dans l'image de $\gamma.$
\bigskip
\goodbreak
 \begin{proof} Consid\'erons l'ouvert $\proj({\cal V})^*$  des
points de 
$\proj({\cal V})$ repr\'esentant une extension (8.1)  satisfaisant aux
conditions suivantes:

 {\rm (a) } Les supports $\ell' $ de $\Theta'$  et $\ell''$  de
$\Theta''$  sont deux droites distinctes, dont on d\'esigne par $O$ 
l'intersection.

{\rm (b)} L'extension
$0\hfl{}{}\Theta'\hfl{}{}\Theta\hfl{}{}\Theta''\hfl{}{}0$ associ\'ee \`a
(8.1) n'est pas scind\'ee.  En particulier,
$\Theta$ est localement libre sur la conique $\hbox{\goth{c}}$  r\'eunion
de $\ell'$  et $\ell'',$  et on a $\Theta|_{\ell'}={\cal
O}_{\ell'}(3)$  et $\Theta|_{\ell''}={\cal O}_{\ell''}(2). $

{\rm (c)}  Le syst\`eme lin\'eaire $\Lambda= (\Gamma,\Theta)$  est sans
point de base. En particulier, le syst\`eme lin\'eaire $\Lambda''$  n'a pas
de point  de base  et d\'efinit une involution ${\cal I}$ sur la droite
$\ell''$ dont les points fixes $a$  et $b$ sont les points associ\'es dans
$\Hilb^{2}(\proj_2).$

{\rm (d)}  Le point $O$ n'est pas fixe sous l'involution ${\cal I}.$

{\rm (e)} L'unique section $s_{O}$ de $\Gamma$  qui s'annule en
$O$  n'a pas de  z\'ero double  sur la droite $\ell'$ (sauf peut-\^etre $O$ 
 qui alors n'est pas un z\'ero
triple).

{\rm (f)}  Les  sections $s_{a}$ et $s_{b}$ de $\Gamma$ qui
s'annulent aux
 deux points fixes $a$  et $b$ de l'involution ${\cal I}$ sur $\ell''$
n'ont pas de z\'ero double sur la droite $\ell'.$

Les points de $\proj({\cal
V})^{*}$ d\'efinissent  des  syst\`emes coh\'erents semi-stables  d'apr\`es 
le lemme \ref{instables}.  On a alors  un morphisme dominant 
$\proj({\cal V})^{*}\hfl{}{} D_{2}.$ 
L'assertion (i)  r\'esulte du corollaire \ref{singuliers-evidents} et
du fait que si
$s$  est un point de
$\proj({\cal V})^*$  alors la quartique
$\hbox{\goth{q}}_{s}$ associ\'ee est singuli\`ere en $\check{\ell}'$: il
suffit de  constater que $\Theta|_{\ell'}={\cal
O}_{\ell'}(3)$. Pour d\'emontrer (ii), nous devons expliquer comment
trouver l'\'equation de
$\hbox{\goth{q}}_{s}$. Consid\'erons le morphisme  de suites exactes  
associ\'e \`a l'extension (8.1):
$$\diagram{
&&0 &\hfl{}{}& \Gamma \otimes {\cal O}&\hfl{}{}& \Gamma''\otimes
{\cal O} &\hfl{}{}&0\cr
&&\vfl{}{}&&\vfl{}{}&&\vfl{}{}\cr
0&\hfl{}{}& \Theta'  &\hfl{}{}& \Theta &\hfl{}{}& \Theta''&\hfl{}{}&
0\cr }$$
Si $\Lambda''$  est sans point de base, le noyau du morphisme de droite
est un faisceau localement libre $F''$ stable de rang 2 de classes de
Chern
$(-1,2)$  et on a une application lin\'eaire   $$
\Ext^{1}(\Lambda'',\Lambda')\hfl{}{} \Hom(F'',\Theta')$$  induite par
l'homomorphisme de liaison figurant dans le diagramme du serpent associ\'e
\`a ce diagramme \footnote{$(^5)$}{Il est facile de voir que cette
application lin\'eaire est inversible, mais nous n'en aurons pas
besoin.}. Par application du foncteur
$\pr_{1*}\pr_2^*(-)$ on voit d'apr\`es le lemme \ref{droites} que l'on
obtient un diagramme  commutatif de suites exactes
$$\diagram{
&&0 &\hfl{}{}& \Gamma \otimes {\cal O}&\hfl{}{}& \Gamma''\otimes
{\cal O} &\hfl{}{}&0\cr
&&\vfl{}{}&&\vfl{\varphi}{}&&\vfl{\varphi''}{}\cr
0&\hfl{}{}& \hbox{\goth{m}}_{\check{\ell}'}^2(2)  &\hfl{}{}& {\cal G}
&\hfl{}{}&
\hbox{\goth{m}}_{\check{\ell}''}^2(2)&\hfl{}{}& 0\cr }.$$ Le noyau de la
derni\`ere fl\`eche verticale est le faisceau
$\pr_{1*}(\pr_{2}^*(F''))$ isomorphe \`a ${\cal O}(-2).$ Le morphisme de
liaison associ\'e fournit un morphisme ${\cal O}(-2)\hfl{}{}
\hbox{\goth{m}}_{\check{\ell}'}^2(2)$   par suite une section $q_{s}$ de ${\cal
O}(4)$  bien d\'efinie \`a homoth\'etie pr\`es.  

\begin{lemme} \label{singuliere.generale} Si
$\Lambda''=(\Gamma'',\Theta'')$  est sans point de base, la section
$q_{s}$  est non nulle.  La quartique
$\hbox{\goth{q}}_{s}$  d'\'equation $q_{s}=0$  est la quartique de L\"uroth
associ\'e au point de $D_{2}$  d\'efini par $s.$  
\end{lemme}

\begin{proof}
On d\'eroule la suite exacte du serpent associ\'ee au diagramme ci-dessus: on
obtient la suite exacte
$$0\hfl{}{} {\cal O}(-2)\hfl{q_{s}}{}
\hbox{\goth{m}}_{\check{\ell}'}^2(2)\hfl{}{} {\cal Q}\hfl{}{} {\cal
Q}''\hfl{}{} 0$$ o\`u ${\cal Q}=\coker(\varphi) $  et ${\cal Q}''=
\coker(\varphi'').$  Ce dernier faisceau est un faisceau de support le
point 
$\check{\ell}''.$ Ceci entra\^{\i}ne \'evidemment $q_{s}\not=0$; le
faisceau 
${\cal Q}$, dont on sait qu'il est pur de dimension 1, a alors pour 
support sch\'ematique 
  la quartique $\hbox{\goth{q}}_{s}$ d'\'equation
$q_{s}=0.$  D'o\`u le lemme. \end{proof}

On retrouve que la quartique $\hbox{\goth{q}}_{s}$ est  singuli\`ere
en
$\check{\ell}'$.
De plus, le morphisme de Barth $\gamma:
\proj(\Ext^{1}(\Lambda'',\Lambda'))\hfl{}{} {\cal C}_{4}$  est une
homographie (donc injective) sous la seule hypoth\`ese  que le syst\`eme
lin\'eaire
$\Lambda''$  soit sans point de base. V\'erifier que si 
$s \in \proj({\cal V})^*$ la quartique $\hbox{\goth{q}}_{s}$ 
appartient \`a ${\cal S}_{4}^*$  demande une \'etude plus pouss\'ee de
l'\'equation de $\hbox{\goth{q}}_{s}$. 

\begin{lemme} \label{a-f}  Soit $s \in
\proj({\cal V})$  satisfaisant aux conditions {\rm (a)-(f)}. 
Alors la quartique de L\"uroth $\hbox{\goth{q}}_{s}$ appartient \`a
${{\cal S}}_{4}^*.$

\end{lemme}

\begin{proof}
Pour obtenir des renseignements sur  la  quartique
$\hbox{\goth{q}}_{s}$  on utilise le fait que l'application ci-dessus
en restriction \`a une  fibre de
$\proj(\Ext^1(\Lambda'',\Lambda'))$ est  une homographie d\`es que
$\Lambda''$  n'a pas de point de base.  Consid\'erons une extension non triviale $(8.1)$ satisfaisant
aux conditions (a)-(f). Elle d\'efinit un fibr\'e inversible $\Theta$ 
sur la conique
$\hbox{\goth{c}}$.  On d\'esigne par
$x=0$  l'\'equation de
$\ell'$  et
$y=0$  l'\'equation de
$\ell''.$ Le fibr\'e 
 $\Theta$  est de degr\'e 3  sur $\ell'$  et de degr\'e $2$  sur
$\ell''.$

 On peut choisir d'apr\`es la condition (d) une forme lin\'eaire $t$  de
sorte que $(x,y,t)$  soit une base de
$V^*$ et  s'annulant au
point conjugu\'e de $O$  par rapport \`a $a$  et $b$: alors  le syst\`eme
lin\'eaire
$(\Gamma'',\Theta'')$  est d\'efini par les quadriques $xt$  et
$x^2+t^2.$
Consid\'erons une section $\tau$
de
$\Theta(-2)$ non nulle sur
$\ell'';$  on peut encore supposer que $t$ est choisie de sorte  que
$\tau
$ s'annule sur
$\ell'$  au point
$x=t=0.$ On d\'esigne par $[u,v,w]$  les coordonn\'ees homog\`enes dans le
dual de $\proj_2:$  un tel point repr\'esente la droite de $\proj_{2}$
d'\'equation $ux+vy+wt=0.$
\unitlength=13pt
$$\beginpicture(12,9)(0,0)
\put(0,1){{\line(1,0){12}}}
\put(7,0.8){$|$}\put(7,0){$a$}
\put(10,0.8){$|$}\put(10,0){$b$}
  \put(1,0){\line(1,2){4}}
  \put(9,0){\line(-2,3){5}}
  \put(0.3,0.3){0}
  \put(0.5,3){$x=0$}
  \put(2.5,4){$\ell'$}
  \put(4,1.3){$y=0$}\put(4,0.3){$\ell''$}
  \put(6.5,4){$t=0$}
  \put(8,7){$\proj_2$}
\endpicture$$
\bigskip
 Consid\'erons  le  point $s_0 \in \proj(\Ext^1(\Lambda'',\Lambda'))$ 
d\'efini par  le syst\`eme lin\'eaire
$(\Gamma_{0},\Theta)$  engendr\'e par les sections $\sigma_{1}=xt
\tau 
$  et
$\sigma_{2}=(x^2+t^2)\tau.$  Evidemment, ce syst\`eme lin\'eaire a pour
point de base le point
$x=t=0,$ et chacune de ces sections s'annule \`a l'ordre $\geq 3$ en ce
point.  Alors la quartique associ\'ee  $\hbox{\goth{q}}_{0}$   contient
la droite triple d\'efinie par le point [0,1,0], d'\'equation
$v^3=0$ et une autre droite.  La section $\sigma_{1}$  s'annule 
identiquement sur $\ell',$  et elle s'annule aussi au point $[1,0,0].$ 
Donc  la quartique $\hbox{\goth{q}}_{0}$ contient aussi la
droite d'\'equation
$u=0$; elle a  pour \'equation
$v^3 u=0.$  Remarquons que cette quartique est lisse au point
$\check{\ell}'',$  et que la tangente en $\check{\ell}''$  correspond
au second point de $\ell''$ o\`u s'annule la section $\sigma_{1} \in
\Gamma_{0}.$

Consid\'erons maintenant une extension $s_{\infty} \in
\proj(\Ext^1(\Lambda'',\Lambda'))$ quelconque  d\'efinissant un  syst\`eme
coh\'erent $(\Gamma_{\infty},\Theta_{\infty})$  tel que l'extension (8.3)  soit
scind\'ee: autrement dit
$s_{\infty}$ appartient \`a l'hyperplan de l'infini d\'efini par les z\'eros de
la forme lin\'eaire
$$\Ext^1(\Lambda'',\Lambda') \hfl{}{} \Ext^1(\Theta'',\Theta')=k$$
Comme $\Hom(\Theta'',\Theta')=0$
une telle extension est d\'efinie par une  application lin\'eaire non
nulle
$f_{\infty}
\in
\Hom(\Gamma'',H^0(\Theta')),$  unique \`a homoth\'etie pr\`es, qu'on peut
voir aussi comme une application lin\'eaire $f:\Gamma''\hfl{}{}
H^0(\hbox{\goth{m}}^2_{\check{\ell}'}(2)).$  Le morphisme $F''\hfl{}{}
\Theta'$  associ\'e \`a $s_{\infty}$ est alors le compos\'e $F''\hfl{}{}
\Gamma''\otimes {\cal O}\hfl{}{} \Theta'$, o\`u la premi\`ere fl\`eche est l'inclusion
canonique, et la seconde se construit en composant $f_{\infty}$ avec le
morphisme d'\'evaluation. L'\'equation
$q_{\infty}$ de la quartique
$\hbox{\goth{q}}_{\infty}$  associ\'ee s'obtient en appliquant le foncteur
$\pr_{1*}(\pr^*_2(-))$: elle est donn\'ee par  le morphisme compos\'e
$${\cal O}_{\proj_2^*}(-2)\hfl{}{} \Gamma'' \otimes{\cal
O}_{\proj_2^*}\hfl{}{}
\hbox{\goth{m}}_{\check{\ell}'}^2(2).$$  On prend pour base de
$\Gamma''$ les formes quadratiques  $xt$  et $x^2+t^2$ correspondant
aux  sections $\sigma_{1}$  et $\sigma_{2}$ de $\Theta''$  d\'efinies  les
 formules  ci-dessus, en restriction \`a $\ell''$; \`a ces sections
correspondent les formes quadratiques $(-uw, u^2+w^2) $  la quartique
$\hbox{\goth{q}}_{\infty}$ a  donc une  \'equation de la forme
$$(u^2+w^2)A_{1}(v,w)+uwA_{2}(v,w)=0$$
o\`u $A_{1}(v,w)$  et $A_{2}(v,w)$  sont des formes binaires de degr\'e 2.
Cette quartique est singuli\`ere en $\check{\ell}'$  et $\check{\ell}''.$
 Par
lin\'earit\'e, on voit que l'\'equation de la quartique de L\"uroth
$\hbox{\goth{q}}_{s}$ associ\'ee  \`a un point quelconque $s$ de
$\proj(\Ext^1(\Lambda'',\Lambda'))$ est de la forme
$$v^3u+ (u^2+w^2) A_{1}(v,w)+uwA_{2}(v,w)=0$$ o\`u $A_1$ et $A_2$  sont
des formes binaires de degr\'e 2, pourvu que la condition (b)  soit
satisfaite. On voit que le point
$\check{\ell}'=[1,0,0]$  est un point de multiplicit\'e 2  pourvu que la
forme binaire $A_{1}(v,w)$ ne soit pas nulle, et que le point
$\check{\ell}''=[0,1,0]$ est un point lisse: la tangente correspond au
deuxi\`eme z\'ero de la section
$\sigma_{1}$ sur $\ell''.$
\unitlength=13pt
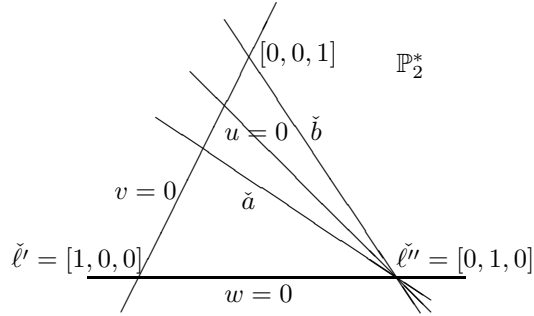
$$\beginpicture(12,10)(0,0)
\put(0,1){{\line(1,0){11}}}
  \put(1,0){\line(1,2){4.5}}
  \put(9,1){\line(-3,2){7}}\put(9,1){\line(3,-2){1}}
  \put(4.5,3){$\check{a}$}
  \put(9,1){\line(-2,3){5}}\put(9,1){\line(2,-3){0.7}}
  \put(9,1){\line(-1,1){6}}\put(9,1){\line(1,-1){1}}
   \put(6.5,5){$\check{b}$}
   \put(0.8,3.3){$v=0$}
  \put(4,0.3){$w=0$}
  \put(4,5){$u=0$}
  \put(9,7){$\proj^*_2$}
  \put(-2.2,1.3){$\check{\ell'}=[1,0,0]$}
  \put(9,1.3){$\check{\ell''}=[0,1,0]$}
  \put(5,7.3){$[0,0,1]$}
\endpicture$$
\bigskip
{\it Le c\^one tangent de $\hbox{\goth{q}}_{s}$ en $\check{\ell}'$.}

Supposons que les conditions (a)-(d) soient satisfaites au
point $s$. Pour comprendre le c\^one tangent  de $\hbox{\goth{q}}_{s}$ au
point singulier $\check{\ell}'$  nous devons \'etudier la forme
quadratique
$A_{1}(v,w).$ Le fibr\'e inversible
$\Theta$ sous-jacent au syst\`eme lin\'eaire
$(\Gamma,\Theta)$ est d\'efini par la m\^eme extension (8.3) que pour le
point $s_{0}$.  D\'esignons par
$\sigma$  l'unique  section  de
$\Theta(-2)$  qui s'annule au point d'intersection $O,$  et li\'ee sur
\`a
$\tau$ par la relation 
${ t \sigma +y \tau =0}.$ 
 Sur $\ell''$  cette section s'annule, et
$\sigma$ peut donc \^etre vue comme une section de $\Theta'(-2).$ 

Puisque $s$ n'a  pas de point de base, $s\not=s_{0}$ et la droite
$s_{0}s$  rencontre l'hyperplan de l'infini en un point
$s_{\infty},$  ce qui d\'efinit  une application  lin\'eaire  $f_{\infty}$  qui est
donn\'ee dans la base de
$\Gamma''$ ci-dessus
  par    deux sections de $\Theta'$ qui s'\'ecrivent sous la
forme
$f_{1}= g_{1}(y,t)\sigma$  et $f_{2}=g_{2}(y,t)\sigma$  o\`u $g_{1}$ et
$g_{2}$ sont des formes binaires,  de sorte que l'espace vectoriel
$\Gamma$ est engendr\'e par les sections 
$xt \tau + g_{1}(y,t)\sigma, (x^2+t^2) \tau + g_{2}(y,t) \sigma.$
Dans l'isomorphisme $H^0(\Theta'') \simeq
H^0(\hbox{\goth{m}}_{\check{\ell}''}^2(2))$,  aux formes binaires $x t$
et
$x^2+t^2$  correspondent les formes binaires $-uw$  et $u^2+w^2$  
respectivement dans
le plan projectif dual. De m\^eme dans l'isomorphisme
$H^{0}(\Theta')
\simeq H^0(\hbox{\goth{m}}_{\check{\ell}'}^2(2))$
 aux  formes
binaires  $g_{i}(y,t)$ correspondent les formes binaires $g_{i}(w,-v)$ 
dans le  plan projectif dual, et l'\'equation de la quartique associ\'ee
est donc de la forme 
$$\lambda v^3 u + (u^2+w^2)g_{1}(w,-v) + uw g_{2}(w,-v) =0
$$
o\`u $\lambda $  est un scalaire non nul que nous allons pr\'eciser.  Pour
ceci, on examine l'intersection de cette quartique  avec la droite de
$\proj^*_2$ d'\'equation $w+u=0:$ cette intersection est donn\'ee par le
sch\'ema des z\'eros de la forme de degr\'e $4$
$$w[-\lambda v^3 + w (2g_{1}(w,-v)- g_{2}(w,-v))].$$  Il s'agit donc
d'interpr\'eter l'intersection du pinceau de  droites de $\proj_2$
passant par le point fixe
$a=[1,0,1]$   avec la	 quartique $\hbox{\goth{q}}_{s}.$ La section de
$\Gamma$ qui s'annule en ce point est $(x-t)^2\tau +
(g_2-2g_{1})(y,t)\sigma$;  elle ne s'annule pas en $O$.
 Par suite,  l'intersection est donn\'ee par le c\^one de
sommet $a$ et d\'efini par le sch\'ema  des z\'eros de la  forme  de degr\'e 3 
sur
$\ell'$  d\'efinie par
$t^3+y(2g_{1}-g_{2})(y,t)$ z\'eros auxquels on doit ajouter la droite
$\ell''.$
  Dans le dual, on obtient ainsi la forme quartique 
$$wv^3-w^2(2g_{1}(w,-v)-g_{2}(w,-v))$$  ce qui prouve que $\lambda=1.$

Remarquons que
les z\'eros de la section qui s'annule en
$0$  sont donn\'es  sur $\ell'$ par l'\'equation $g_{1}(y,t) \sigma=0$ 
donc  d'apr\`es  la condition  (c) la forme binaire $g_{1}$ n'est pas
identiquement nulle: ceci implique  que
$\check{\ell}'$  est un point double sur $\hbox{\goth{q}}_{s}.$

\bigskip\goodbreak
{\it La polaire de $\hbox{\goth{q}}_{s}$ au point singulier
$\check{\ell}'$ , et la conique $\hbox{\goth{p}}$ associ\'ee.}

L'\'equation de la quartique $\hbox{\goth{q}}_{s}$  s'\'ecrit en
choisissant 
$u=0$  comme droite de l'infini \'evitant le point singulier
$\check{\ell}'= [1,0,0]$
$$u^2(g_{1}(w,-v))+ u(v^3+ wg_{2}(w,-v))  +w^2 g_{1}(w,-v) =0$$
Consid\'erons les 2 sections  qui engendrent $\Gamma.$  Sur $\ell'$ elles
s'annulent aux points  d\'efinis par les \'equations 
$$\eqalign{yg_{1}(y,t)&=0\cr
 t^3  - yg_{2}(y,t)& =0 \cr
}$$
D'apr\`es l'hypoth\`ese (c)  le syst\`eme lin\'eaire est sans point de base.
Par cons\'equent, ces deux \'equations n'ont pas de z\'ero commun sur $\ell'.$
Il est \'equivalent de dire que dans le dual les formes binaires 
$wg_{1}(w,-v)$ et $v^3  + wg_{2}(w,-v)$  n'ont pas de z\'ero en commun.
Ainsi, la cubique polaire $\hbox{\goth{q}}'$ de $\ell'$, est
int\`egre.  
La condition (e)  signifie  exactement que la forme quadratique
$g_{1}$  est de discriminant non nul, ce qui implique que $\ell'$
est un point double ordinaire sur la quartique. La condition (f)
implique  que d'apr\`es la proposition \ref{points-singuliers} la
quartique ne contient pas d'autre point singulier que
$\check{\ell}'$ ou
$\check{\ell}''$.  Mais  ces points singuliers sont forc\'ement \`a
l'intersection de
$\hbox{\goth{q}}_{s}$ et de la polaire
$\hbox{\goth{q}}_{s}'$  et on a d\'ej\`a vu que $\check{\ell}''$  n'est
pas un point singulier.  Donc il n'y a pas d'autre point singulier.
Ainsi, la quartique $\hbox{\goth{q}}_{s}$  appartient \`a ${\cal
S}_{4}^*.$

 Avec les notations de la section 7, on voit que $\phi=0, \psi=
w^2$ la conique associ\'ee $\hbox{\goth{p}}$  est alors d\'efinie par
l'\'equation
$u^2-w^2=0$:  c'est donc la  conique singuli\`ere de point singulier le
point $\check{\ell}''$. Cette conique est celle qui est associ\'ee \`a la
forme binaire $x^2-t^2$   sur $\ell'':$  il s'agit bien s\^ur des points
fixes de l'involution sur $\ell''$  associ\'ee au syst\`eme lin\'eaire
$\Lambda''.$
Ceci
ach\`eve la d\'emonstration du lemme \ref{a-f}. \end{proof}

\bigskip

{\it Fin de la d\'emonstration du th\'eor\`eme \ref{cle}}

Consid\'erons  l'image r\'eciproque  $\proj({\cal V})^{**}$ de l'ouvert 
${\cal S}_{4}^*$  par $\gamma: \proj({\cal
V})\setminus I\hfl{}{} {\cal S}_{4}.$   Puisque $\proj({\cal V})^*$ est
dense dans
$\proj({\cal V})^{**}$  cette image r\'eciproque est aussi l'image r\'eciproque
de
${\cal L}_{2}^*.$
D\'esignons par  $\Hilb^2_*(\proj_{2})$  le compl\'ementaire de la diagonale
dans $\Hilb^2(\proj_{2}).$ On a    une application r\'eguli\`ere
${\cal L}_{2}^*\hfl{}{}\proj^*_2
\times
\Hilb^{2}_*(\proj_2)$  qui associe \`a $\hbox{\goth{q}}$  le couple form\'e
du point singulier $\omega$  de $\hbox{\goth{q}}$, et  du point de
$\Hilb^2_{*}(\proj_{2})$  d\'efini par la conique singuli\`ere $\hbox{\goth{p}}$
de rang 2 associ\'ee \`a
$\hbox{\goth{q}}.$  Cette paire d\'efinit un couple  de syst\`emes coh\'erents
$(\Lambda',
\Lambda'')$  de support distincts et tels que  $\Lambda''$ soit  sans point
de base.  On a donc un diagramme commutatif
$$\diagram{
\proj({\cal V})^{**}&&\hfl{\gamma}{}&& {\cal L}_{2}^*\cr
&\sefl{}{}&& \swfl{}{}&\cr
&&\proj^*_2 \times
\Hilb^{2}_*(\proj_2)&&\cr
}$$
Au-dessus  de l'ouvert 
$\Hilb^{2}_*(\proj_2)$ le morphisme $\gamma$  est un plongement propre
puisque cela signifie exactement  que $\Lambda''$ n'a pas de point de base. 
  Le diviseur  ${\cal L}_{2}^*$  est irr\'eductible de
 dimension 12  et  $\proj({\cal V})^{**}$    est aussi de dimension 12: 
c'est donc un isomorphisme. \cqfd

\begin{prop}
L'ouvert $\proj({\cal V})^{**}$  co\"{\i}ncide avec l'ouvert $\proj({\cal
V})^*$.
\end{prop}
\begin{proof}
 En effet, soit $s \in
\proj(\Ext^1(\Lambda'',\Lambda'))$  un point d\'efinissant un syst\`eme
coh\'erent semi-stable $(\Gamma,\Theta)$ tel  que  la quartique associ\'ee
$\hbox{\goth{q}}_{s}$  appartienne \`a ${\cal S}_{4}^*.$ Cette quartique
\'etant int\`egre, ce syst\`eme coh\'erent est en fait stable, et  n'a pas de
point de base. Par continuit\'e,
le support $\hbox{\goth{c}}$ de  
$\Theta$  est une conique singuli\`ere.
 Supposons que cette conique soit une droite  double $\ell^2.$  Alors
$\Theta$ est extension 
$$0\hfl{}{}{\cal O}_{\ell}(2)\hfl{}{}\Theta\hfl{}{}{\cal
O}_{\ell}(2)\hfl{}{}0 \eqno (8.4) \quad$$ Si cette extension n'est pas
scindable
$\Theta$ est singulier en un seul point $O,$  et dans le pinceau
de droites passant par $O,$  seule la droite $\ell$  serait de saut 
d'apr\`es la proposition \ref{points-singuliers}.  Ceci signifie  que la
droite de $\proj_2^*$  d\'efinie par $O$  rencontre $\hbox{\goth{q}}_{s}$
au seul point $\check{\ell}.$  Ceci contredit le fait que
$\hbox{\goth{q}}_{s}$  n'a qu'un point double ordinaire.  Si
l'extension $\Theta$  est scind\'ee, $\Theta$  est un ${\cal
O}_{\ell}-$module localement libre de rang 2  et de degr\'e $4$ et le
syst\`eme coh\'erent $(\Gamma,\Theta)$ aurait des points de base. Donc  la
conique
$\hbox{\goth{c}}$ est r\'eduite et la condition (a)  est satisfaite. La 
 suite exacte (8.3)  ne peut \^etre  scindable car alors   chacune des
deux droites contenues dans
$\hbox{\goth{c}}$ d\'efinirait un point singulier de
$\hbox{\goth{q}}_{s}$.  Donc la condition (b) est satisfaite.  La
condition (c)  est vraie parce $\hbox{\goth{q}}_{s}$  est int\`egre.
Par continuit\'e, la conique $\hbox{\goth{p}}$ associ\'ee \`a
$\hbox{\goth{q}}_{s}$  est une conique singuli\`ere dont le point
singulier est le $\check{\ell''}.$  Puisque cette conique ne rencontre
pas  le point $\check{\ell'}$, le point $O$  n'est pas un des points
fixes    de l'involution ${\cal I}$  sur $\ell'',$  ce qui donne la
condition (d). 
Comme nous l'avons vu dans la preuve du lemme \ref{singuliere.generale}
la condition (e)  est vraie parce cela traduit exactement  que le point
$\check{\ell}'$  est un point double ordinaire de $\hbox{\goth{q}}_{s}$
et la condition (f)  r\'esulte du fait que c'est le seul point singulier
de
$\hbox{\goth{q}}_{s}.$ \end{proof}

La vari\'et\'e des quartiques de L\"uroth  singuli\`eres 
de type II  est, par d\'efinition,  l'image ${\cal
L}_{2}=\gamma(D_{2})$  du diviseur
$D_{2}.$ Puisque $D_{2}$ est irr\'eductible,  il en est de m\^eme de ${\cal
L}_{2}.$ 
\begin{cor}\label{cor87}
La vari\'et\'e ${\cal L}_{2}$ des quartiques de L\"uroth singuli\`eres de type
$II$ est l'adh\'erence dans ${\cal C}_{4}$ de l'hypersurface de
${\cal S}_{4}^*$ d\'efinie par l'\'equation
$$\Delta(\phi^2+\psi)=0$$
o\`u $\Delta$  est le discriminant.  En particulier, le morphisme ${\cal
L}_{2}\hfl{}{} {\cal C}_{2}$  d\'efini par $\hbox{\goth{q}}\mapsto
\hbox{\goth{p}}$  n'est pas dominant.
\end{cor}

{\it D\'emonstration du lemme \ref{diviseur2}}

  Consid\'erons un point de $\proj({\cal V})^*$  d\'efinissant une extension
(8.1). 
  Le syst\`eme coh\'erent
$(\Gamma'',\Theta'')$  d\'efinit  d'apr\`es [16] un faisceau stable $F''$ de
rang 2  et classes de Chern $(1,2)$: c'est  le dual du noyau du
morphisme  d'\'evaluation
$\Gamma''\otimes {\cal O}\hfl{}{} \Theta''.$ 
Le faisceau semi-stable $F$  de
$M_{4}$ associ\'e \`a 
$(\Gamma,\Theta)$  dans la correspondance de la section 4 s'ins\`ere dans une
suite exacte $$0\hfl{}{} F''\hfl{}{} F(1) \hfl{}{} {\cal
O}_{\ell'}(-1)\hfl{}{} 0$$ ce qui implique que $\ell'$  est une droite
de saut d'ordre 2 pour
$F$, \ie $\dim H^{1}(F(-1)|_{\ell'})=2.$ Il en r\'esulte que l'image 
 par la contraction
$\pi:S_{4}\hfl{}{}M_{4}$ du diviseur $D_{2}$ est  contenu dans   le
diviseur des classes de faisceaux semi-stables
$F$  qui ont au moins une droite de saut d'ordre $\geq 2.$
R\'eciproquement, si $F$  poss\`ede une droite de saut  $\ell'$ d'ordre $\geq 2$  
et n'est pas sp\'ecial,  le conoyau  du morphisme d'\'evaluation
$\ev:H^0(F(1))\otimes {\cal O}\hfl{}{} F(1)$ poss\`ede un quotient isomorphe \`a
${\cal O}_{\ell'} (-1)$  et par suite son  support sch\'ematique est une
conique d\'eg\'en\'er\'ee.  Donc le syst\`eme coh\'erent $(\Gamma,\Theta)$  
associ\'e \`a $F$  appartient au diviseur $D_{2}.$
\cqfd

\section{Fin de la d\'emonstration, dans le cas $n=4.$}

Consid\'erons le ferm\'e $\gamma(\partial S_{4} \cup \hbox{\goth{e}})$; ce
ferm\'e est constitu\'e de quartiques r\'eductibles ou de quartiques ayant un
point triple, et  est donc de dimension $\leq 11$; soit ${\cal C}_{4}^{*}$
le compl\'ementaire de ce ferm\'e dans ${\cal C}_{4}$. On sait  que
$\gamma$ est d'apr\`es [13] un morphisme fini   au-dessus de l'ouvert
${\cal C}_{4}^{*}$.  Les diviseurs
$D_{1}$  et
$D_{2}$  rencontrent l'image r\'eciproque de cet ouvert ${\cal C}_{4}^*$. 
Donc les images
${\cal L}_{1}$ et ${\cal L}_{2}$  de
$D_{1}$  et
$D_{2}$  sont deux diviseurs de  l'hypersurface ${\cal L}=\beta(M_{4})$ des
quartiques de L\"uroth appel\'es respectivement diviseurs des quartiques de
L\"uroth singuli\`eres de type I et type II; on va voir que ces deux
diviseurs sont distincts; ce sont les seules composantes irr\'eductibles de
${\cal L}\cap {\cal S}_{4}.$  

On doit v\'erifier 
que le diviseur 
${\cal L}_{1}$  n'est pas contenu dans ${\cal L}_{2}.$  Pour ceci on
consid\`ere la quartique associ\'ee \`a un point de $D_{1}.$  En g\'en\'eral, cette
quartique appartient \`a ${\cal S}^*_{4}$, ce qui d\'efinit  d'apr\`es la section
8  une conique $\hbox{\goth{p}}.$ Il suffit d'apr\`es le corollaire
\ref{cor87} de prouver le r\'esultat suivant:
\begin{prop} \label{D1}
L'application rationnelle  $ D_{1}\rightsquigarrow  {\cal C}_{2}$ 
qui associe \`a un syst\`eme  lin\'eaire la conique $\hbox{\goth{p}}$  
est dominante.  

\end{prop}

On doit d'autre part v\'erifier l'\'enonc\'e suivant: 

\begin{prop} \label{rang}
 Le rang de $d\gamma$  en un point 
g\'en\'eral de $D_{2}$  est $13.$  

\end{prop}

D\`es lors, la fibre g\'en\'erale du morphisme $\gamma: S_{4}\hfl{}{} {\cal
C}_{4}$  est r\'eduite \`a un point au-dessus d'un point g\'en\'eral de ${\cal
L}_{2}$, et $\gamma$ non  ramifi\'e en ce point. Il en r\'esulte (\cf
Hartshorne, Chapitre II, lemme 7.4) que le morphisme 
$\gamma$  et  le  morphisme de Barth $\beta$ sont g\'en\'eriquement 
injectifs. Ceci prouve le th\'eor\`eme \ref{main}  dans le cas $n=4.$
\bigskip\goodbreak
\subsection{D\'emonstration de la proposition 9.2}

Il suffit de trouver un point $s\in D_{2}$ et un vecteur tangent
$\tau
\in T_{s}S_{4}$  tel que 

--- (1) le morphisme $\gamma: D_{2}\hfl{}{} {\cal C}_{4}$  soit un
plongement au voisinage de $s;$ 

--- (2) le vecteur $d_{s}\gamma
(\tau)$   n'appartienne pas \`a l'espace tangent \`a $T_{\beta(s)}{\cal
L}_{2}.$ 

Consid\'erons la famille param\'etr\'ee par $\proj_1$  d\'efinie par le
conoyau du morphisme de faisceaux localement libres sur
$\proj_1\times
\proj_2$
$$\psi:{\cal O}(-1,1)\oplus {\cal O}(0,1)\hfl{}{} {\cal O}^2(0,2)$$
dont la matrice est donn\'ee par
$$\pmatrix{\alpha x&t\cr
\beta t&y}$$
Le conoyau de $\psi$ d\'efinit une famille de faisceaux semi-stables de
classe de Grothendieck $c=2h+4h^2,$  dont le support sch\'ematique est le
sous-sch\'ema
$C$ de
$\proj_1\times \proj_2$  d\'efini par l'\'equation \`a valeurs dans ${\cal
O}(1,2)$
$$\alpha xy-\beta t^2 =0.$$  Le morphisme $\proj_1\hfl{}{} {\cal C}_2$
ainsi obtenu est transverse  au point
$\infty=[1,0]$ au diviseur des coniques singuli\`eres. D\'esignons par $\sigma$ 
et
$\tau$  les deux sections canoniques de $\Theta(-2)$  obtenues; au-dessus
du point $\infty$ elles s'annulent  respectivement  sur la droite $\ell''$
d'\'equation $y=0$  pour $\sigma,$  et au point  $x=t=0$  pour $\tau.$
 On se place dans ce qui suit au voisinage de $\infty,$  et on pose
$\epsilon= \beta/\alpha.$ La famille param\'etr\'ee de faisceaux ${\cal G}$
obtenue sur $\proj^*_2$  est facile \`a calculer: en effet, le module
gradu\'e $\oplus_{i} H^{0}(\Theta(i))$ sur l'anneau des polyn\^omes
  est engendr\'e par $H^0(\Theta(-2))$, et 
les sections $\sigma$  et $\tau $  satisfont 
aux relations  
$x\sigma +\epsilon t \tau=0; t\sigma + y \tau=0.$ 
Alors on peut prendre pour rep\`ere pour le fibr\'e
$\pr_{2}^*(\Theta(-1))$  les sections  $y\sigma,t\sigma,x\tau, t\tau.$ 
De m\^eme, pour le fibr\'e $\pr_{2*}(\Theta)$  on peut prendre pour rep\`ere
$(y^2\sigma,yt\sigma,t^2\sigma,x^2\tau,xt\tau,t^2\tau)$. Une droite  de
$\proj_{2}$  a pour \'equation $ux+vy+wt=0;$  le point associ\'e dans le plan
projectif dual a pour coordonn\'ees homog\`enes $[u,v,w].$
Le morphisme canonique
$$\pr_{1*}(\Theta(-1))\boxtimes {\cal O}_{\proj_{2}^*}(-1)\hfl{}{}
\pr_{1*}(\Theta)\boxtimes {\cal O}_{\proj_2^*}$$  a pour matrice
$$\pmatrix{
v& 0&0&0\cr
w&v&0&0\cr
\epsilon u&w&0&-v\cr
0&0&u&0\cr
0&0&w&u\cr
0&-\epsilon u&\epsilon v&w\cr
}$$
La famille de faisceaux ${\cal G}_{\epsilon}$ associ\'ee est le conoyau de ce
morphisme. Choisissons maintenant un  sous-fibr\'e de rang 2
$\Gamma\subset\pr_{1*}({\cal \Theta})$  tel que  au-dessus
de
$\epsilon=0$,
$(\Gamma,\Theta)_{0}$ n'appartienne pas \`a
$D_{1}.$  Pour ceci, on  consid\`ere le sous-espace vectoriel
$\Gamma_{\epsilon} \subset H^0(\Theta_{\epsilon})$ engendr\'e par les
sections
$s_+= (y^2+t^2)\sigma + (x+t)^2\tau$  et $s_{-}=-(y^2+t^2)\sigma +(x-t)
^2\tau.$  Le sch\'ema des z\'eros de $s_{+}$ est obtenu en \'ecrivant que la
matrice
$$\pmatrix{
y^2+t^2 & x &t\cr
(x+t)^2&0&y\cr
}$$
est de rang $\leq $1, ce qui donne sur la droite $\ell''$, $(x+t)^2=0$; de m\^eme
sur la droite $\ell'$  d'\'equation $x=0,$ 
$t^3-t^2y-y^3=0$.  De m\^eme, la section $s_{-}$ s'annule sur la droite
$\ell''$, au point double d\'efini par $(x-t)^2=0$  et sur $\ell'$  en la
cubique
$t^3+t^2y+y^3=0$.
 Il en r\'esulte que le
syst\`eme lin\'eaire est sans point de base. On v\'erifie que le syst\`eme
lin\'eaire obtenu pour $\epsilon=0$  satisfait  en fait  \`a toutes les
 conditions (a)-(f)  de la section 9, et donc appartient \`a l'ouvert
$\proj({\cal V})^{*}$; par suite, d'apr\`es
 le th\'eor\`eme
\ref{cle}, la condition (1) est satisfaite au point de
$D_{2}$ ainsi d\'efini.  Calculons maintenant la famille de quartiques  de
L\"uroth
$\hbox{\goth{q}}_{\epsilon}$ associ\'ees: on prend pour base dans
$\Gamma$  les sections
${1\over 2}(s_{+}\pm s_{-}).$ La famille 
$\hbox{\goth{q}}_{\epsilon}$  de quartiques de saut associ\'ees  est
obtenue en calculant le d\'eterminant du morphisme  canonique sur $k
\times \proj_2^*$
$$\Gamma \boxtimes {\cal O}_{\proj_2^*}\oplus
\pr_{1*}(\Theta(-1))\boxtimes {\cal O}_{\proj_2^*}(-1)
\hfl{}{} \pr_{1*}(\Theta)$$  ce qui donne 
$$\det\pmatrix{
1&0&v& 0&0&0\cr
0&0&w&v&0&0 \cr
1&0&\epsilon u&w&0&-v\cr
0&1&0&0&u&0\cr
2&0&0&0&w&u\cr
0&1&0&-\epsilon u&\epsilon v&w\cr
}=0.$$
On obtient
$$[(u^2+w^2)(v^2+w^2)+2uv^3]-\epsilon (vu^3+3uvw^2+uv^3+2v^4)
+\epsilon^2u^2v^2=0.$$
La quartique de L\"uroth $\hbox{\goth{q}}_{0}$ d'\'equation
$(u^2+w^2)(v^2+w^2)+2uv^3=0$ est singuli\`ere au seul point
$\ell'=[1,0,0]$ et appartient \`a ${\cal L}_{2}^*.$ On prend comme droite
de l'infini la droite
$\Delta$ d\'efinie par l'\'equation
$u=0.$ La question qui nous pr\'eoccupe est de montrer que     le
vecteur tangent  d\'efini par 
$g= vu^3+3uvw^2+uv^3+2v^4$ n'appartient pas \`a
$T_{\hbox{\goth{q}}_{0}}{\cal L}^*_{2}$.  Il suffit de se convaincre que la
conique  d'\'equation $2u\dot \phi - \dot\psi =0$  ne passe pas par le
point singulier de la conique  $\hbox{\goth{p}}_{0}$ associ\'ee \`a
$\hbox{\goth{q}}_{0}.$ L'\'equation de
$\hbox{\goth{q}}_{0}$ s'\'ecrit
$$f= u^2f_2+ uf_3+f_4=0$$
o\`u  les formes binaires $f_{i}$  de degr\'e $i$  sont d\'efinies par
$$\eqalign{
f_2&= v^2+w^2\cr
f_3&=2v^3\cr
f_4&
=w^2(v^2+w^2)\cr
}$$
L'application lin\'eaire
$$H^0({\cal O}_{\Delta}(2)) \oplus H^0({\cal
O}_{\Delta}(1))\hfl{(f_2 \   f_3)}{} H^{0}({\cal O}_{\Delta}(4))$$
est inversible, puisque $f_{2}$  et $f_{3}$ n'ont pas de z\'ero en
commun. La conique $\hbox{\goth{p}}_{0}$  associ\'ee s'obtient en
r\'esolvant l'\'equation
$$2v^3 \phi+(v^2+w^2)\psi= w^2(v^2+w^2)$$ 
qui a pour seule solution $\phi=0, \psi= w^2.$  Donc la conique
$\hbox{\goth{p}}_{0}$ a pour \'equation $u^2-w^2=0$: c'est bien  comme
nous l'avons vu dans la section 8 une conique singuli\`ere r\'eunion de 
deux droites  passant par
$\ell''=[0,1,0]$  correspondant aux deux points doubles du syst\`eme lin\'eaire
sur la droite $\ell''.$ On a ici
$g=u^3v+u(3vw^2 +v^3)+ 2v^4.$ Ecrivons $\xi= (\dot v,\dot w);$
l'\'equation \`a r\'esoudre est 
$$\eqalign{2v\dot v&= -v\cr
2\dot w&= 0\cr
2v^3 \dot \phi + (v^2+w^2) \dot \psi&= 2v^4- 6 w^2v^2\dot v 
\cr }$$
ce qui donne comme solution $\dot v= -\ds{1\over 2}, \dot w=0,\dot
\phi=-\ds{v\over 2}$  et
$\dot
\psi= 3v^2.$ La forme quadratique obtenue est donc $-uv-3v^2$  et ne
passe pas par le point $\ell''.$  Donc ceci d\'efinit un vecteur  de
$H^0({\cal O}_{\hbox{\goth{p}}_{0}}(2))$ qui n'est pas tangent \`a 
l'hypersurface  des coniques singuli\`eres.   \cqfd
\bigskip
\subsection{D\'emonstration de la proposition 9.1}

Consid\'erons la quartique de Luroth $\hbox{\goth{q}}_{0}$ singuli\`ere de
type I
  associ\'ee au syst\`eme lin\'eaire  
sur la conique $\hbox{\goth{c}}$  d'\'equation
$xy-t^2=0$  (correspondant \`a $\epsilon=1$ dans la famille ci-dessus)
d\'efinie par  les sections 
$t^2(\sigma +
\tau)
$ et
$-y^2 \sigma + x^2\tau$. Ce syst\`eme lin\'eaire appartient \'evidemment au
diviseur $D_{1}$  et est sans point de base: la droite $t=0$
correspond, dans le dual, au seul point singulier  de la quartique de
L\"uroth. Cette quartique de L\"uroth singuli\`ere de type I  a pour \'equation
$$\det\pmatrix{
0&-1&v& 0&0&0\cr
0&0&w&v&0&0 \cr
1&0& u&w&0&-v\cr
0&1&0&0&u&0\cr
0&0&0&0&w&u\cr
1&0&0&- u& v&w\cr
}=0.$$
c'est-\`a-dire 
$w^2f_2+wf_3+f_4=0$
o\`u $f_{2},f_{3},f_{4}$  sont des formes binaires  de
degr\'e 2, 3, 4 d\'efinies par 
$$\eqalign{f_2&=u^2+v^2\cr
f_3&=u^3+v^3\cr
f_{4}&= -uv[u^2+v^2]\cr
}$$ 
Soit
$\Delta$  la droite de l'infini, d\'efinie par l'\'equation $w=0.$ Les formes
binaires 
$f_{2}$ et
$f_{3}$  n'ont pas de racine commune.
Ainsi, cette quartique appartient \`a ${\cal S}_{4}^*.$
  La conique $\hbox{\goth{p}}$ associ\'ee  a pour 
\'equation 
$w^2+uv=0.$ C'est donc une conique non singuli\`ere. On passe  facilement 
de la conique duale
$\check{\hbox{\goth{c}}}$ de la conique $\hbox{\goth{c}}$ \`a la conique
$\hbox{\goth{p}}.$  En effet, la conique duale 
$\check{\hbox{\goth{c}}}$ a pour \'equation $w^2-4uv=0$, et 
$\hbox{\goth{p}}$ est l'image de 
$\check{\hbox{\goth{c}}}$ par l'homographie $h:\proj^*_2\hfl{}{}\proj_2^*$ 
caract\'eris\'ee par les propri\'et\'es suivantes:
$h$ laisse fixe le point
$O=[0,0,1]$  et les points de la droite d'\'equation $w=0,$ et le birapport 
 $[\infty,0,m,h(m)]$  est ${i\over 2}.$  
Pour obtenir l'\'enonc\'e, il suffit donc de faire varier la conique
$\hbox{\goth{c}}$.
\end{proof}

 M\^eme si la conique  associ\'ee $\hbox{\goth{p}}$  est une conique lisse,
 il n'est pas toujours aussi simple de construire une
homographie qui  transforme la conique
$\check{\hbox{\goth{c}}}$ en la conique  $\hbox{\goth{p}}.$ 
Quand la
conique
$\hbox{\goth{c}}$  est fix\'ee, la conique $\hbox{\goth{p}}$ d\'epend du
choix du syst\`eme lin\'eaire  sur $\hbox{\goth{c}}.$  Ceci rend
difficile la reconstruction directe de la conique $\hbox{\goth{c}}$  \`a
partir de la quartique $\hbox{\goth{q}}.$  Il peut aussi arriver  que la
conique
$\hbox{\goth{p}}$  soit singuli\`ere: ceci  fournit
alors des exemples explicites (\cf \S 9.3) de fibr\'es stables non isomorphes
qui ont m\^eme quartique de L\"uroth.

\begin{cor} L'hypersurface ${\cal L}$  de ${\cal C}_{4}$  des quartiques
de L\"uroth est de degr\'e 54.

\end{cor}

\begin{proof} Cet corollaire est un cas particulier du corollaire 1.2 
qui se d\'emontre exactement par la m\^eme m\'ethode.
Soit ${\cal D}=\gamma^*({\cal O}(1)) $ est le fibr\'e
d\'eterminant de  Donaldson,  et le nombre de Donaldson 
$q_{13}=\int_{M_{4}}c_{1}({\cal D})^{13}$  est 54.  On a alors 
$\beta_{*}(1)= 54 h$, et par suite, puisque $\beta$  est g\'en\'eriquement
injectif, la classe fondamentale de l'image ${\cal L}$ est $54 h.$
\end{proof}

\begin{cor} \label{degreII} La sous-vari\'et\'e  ${\cal L}_{2}$ des
quartiques de L\"uroth de type II est de degr\'e
$15.27= 5.3^{4}.$
\end{cor}
\begin{proof}
Soit $h$  la classe d'une section hyperplane dans ${\cal C}_{4}.$
 On a $D_{1}= 12 \hbox{\goth{d}} \mod \hbox{\goth{e}}$ d'apr\`es la formule
\ref{calcul}. Il r\'esulte de la formule de projection et du fait que
$\beta_{*}(1)=54 h$  que 
$\beta_{*}([D_{1}])= 12. 54 h^2.$  D'autre part, ${\cal
L}\cap {\cal S}$  est de degr\'e $27.54.$ D'autre part,  dans $A^{2}({\cal
C}_{4})$  on a d'apr\`es la formule (6.2)
$$[{\cal L}\cap {\cal S}]= \beta_{*}([D_{1}])+2 \beta_{*}([D_{2}])$$
Il en r\'esulte que la vari\'et\'e
${\cal L}_{2}$  est de  degr\'e $15.27=5.3^{4}.$ 
\end{proof}
 Le degr\'e de ${\cal L}_{1}$  divise \'evidemment $12.54=2^3.3^{4}.$  Faute de
pouvoir pr\'eciser   le degr\'e de $\beta:D_{1}\hfl{}{} {\cal L}_{1},$  nous
ignorons quel est ce degr\'e; autrement dit nous ignorons quelle est la
multiplicit\'e de
${\cal L}\cap {\cal S}$ le long de ${\cal L}_{1}.$

\subsection{Deux fibr\'es stables de m\^eme quartique de L\"uroth}
De tels exemples sont d\'ej\`a mentionn\'es par Barth [3]: les quartiques
dites desmiques de Humbert sont les courbes de saut de 6 fibr\'es stables
distincts.  L'exemple  pr\'esent\'e ici est explicite. Il s'agit de
constater qu'il existe des points de $D_{1}\setminus D_{2}$ dont l'image
par $\gamma$ appartient \`a ${\cal L}_{1}\cap {\cal L}_{2}.$

On consid\`ere comme au paragraphe \S 9.2 le syst\`eme lin\'eaire sur la
conique d'\'equation $xy-t^2=0$ d\'efini par les sections $t^2(\sigma+\tau)$ 
et
$(-y^2+c t^2)\sigma + (x^2-ct^2)\tau$ o\`u
$c$  est un param\`etre.  L'\'equation de la quartique $\hbox{\goth{q}}_{c}$
associ\'ee est donn\'ee par
$$\det\pmatrix{
0&-1&v& 0&0&0\cr
0&0&w&v&0&0 \cr
1&c& u&w&0&-v\cr
0&1&0&0&u&0\cr
0&0&0&0&w&u\cr
1&-c&0&- u& v&w\cr
}=0$$
L'\'equation de la quartique de L\"uroth $\hbox{\goth{q}}_{c}$ associ\'ee est 
$w^2f_{2}+wf_{3}+f_{4}=0,$ o\`u on a pos\'e
$$\eqalign{f_{2}&=u^2+v^2 \cr
f_3&=u^3+v^3\cr
f_{4}&=-u^3v-2cu^2v^2-uv^3\cr}$$
ce qui donne pour  l'\'equation $f_{3} \phi+f_{2}\psi=f_{4}$  la solution 
$\phi=c(u+v)$ et
$\psi=-[ c u^2+(1+c)uv+cv^2].$  Le discriminant de la forme quadratique
$w^2+2t\phi -\psi$  est alors $\Delta=-{1\over 4}(4c+1)(c-1)^2.$  Pour
$c\not=1,$  le syst\`eme lin\'eaire est sans point de base, et pour $c=-{1\over
4}$  on obtient une quartique de L\"uroth singuli\`ere \'evidemment  de type
I; montrons qu'elle  est de type II. On peut v\'erifier que pour
$c=-{1\over 4}$ le syst\`eme lin\'eaire ne contient qu'une section ayant deux
z\'eros doubles. La quartique
$\hbox{\goth{q}}_{-{1\over 4}}$  est donc lisse
  en dehors du point
$O=[0,0,1]$ et ce point singulier est \'evidemment un point double
ordinaire. La cubique polaire du point singulier $O
$ \'etant \'evidemment
int\`egre, la quartique $\hbox{\goth{q}}_{-{1\over
4}}$ appartient au diviseur ${\cal L}_{2}^*:$ d'apr\`es la section \S 8,
elle provient de $D_{2}$.  
Plus pr\'ecis\'ement, la solution est  alors dans ce cas
$\phi=-{1\over 4}(u+v),
\psi={1\over 4}(u^2-3uv+v^2),$ et $\phi^2+\psi= -{5\over {16}}(u-v)^2$. La
conique   $\hbox{\goth{p}}$ associ\'ee est alors la r\'eunion de deux
droites  qui se coupent au point d\'efini par les \'equations
$4w=u+v,u=v$  c'est-\`a-dire au point
$[2,2,1]
$ du plan projectif dual. Ce point repr\'esente la droite  de $\proj_2$
d'\'equation
$2x+2y+t=0.$
On peut alors reconstruire d'apr\`es la section \S 8 un syst\`eme lin\'eaire 
sans point de base,  appartenant \`a
$D_2$ et dont le support est la conique singuli\`ere d\'efinie par l'\'equation
$t(t+2x+2y)=0$ et qui a pour quartique de L\"uroth
$\hbox{\goth{q}}_{-{1\over 4}}$.

  Par suite, on obtient deux fibr\'es
vectoriels stables et non isomorphes qui ont m\^eme image dans ${\cal
C}_{4}.$
\bigskip

\bigskip\goodbreak
\centerline{\capit La r\'ecurrence}\bigskip
Nous supposons d\'esormais  $n\geq 5$, et que le th\'eor\`eme \ref{main}  est d\'emontr\'e
en degr\'e $c_2=n-1.$
\section{La restriction du morphisme de Barth au bord $\partial M_{n}$ }

Rappelons les notations: on d\'esigne par
$U_{n}$ l'ouvert de  l'espace de modules $M_{n}$  des classes de
fibr\'es vectoriels stables  de rang 2, de classes de Chern
$(0,n)$ sur le plan projectif, 
par ${\cal C}_{n}$  le syst\`eme lin\'eaire
complet des courbes de degr\'e $n$  dans le plan projectif dual, et par 
$$\beta:M_{n}\hfl{}{} {\cal C}_{n}$$  le
morphisme de Barth. La fronti\`ere $\partial M_{n}$  est un diviseur
irr\'eductible de
$M_{n}$;  la restriction  de $\beta$ au bord  $\partial M_{n}$
a ses fibres 
de  dimension $\geq 1$, et g\'en\'eriquement de dimension 1 [23]. L'image 
$B:=\beta(\partial M_{n})$   est  donc un
ferm\'e irr\'eductible  de dimension $4n-5$  de ${\cal C}_{n}.$  On se
propose de d\'emontrer que l'hypoth\`ese de r\'ecurrence entra\^{\i}ne qu'au-dessus
d'un   ouvert non vide de
$B,$  les fibres de $\proj_1$ sont r\'eduites \`a
$\proj_1.$

\begin{lemme} Un faisceau $F$ 
 semi-stable de rang 2 et classes de Chern $(0,n)$ singulier en au plus un
point  est
$\mu-$stable. 
\end{lemme}
\begin{proof}
Si $\lg(\underline{\Ext}^{1}(F,{\cal O}))\leq 1$, le faisceau $F$ n'a pas
de sous-faisceau $F'$ de rang 1, de classe de Chern $c_{1}=0,$ 
 singulier en plus  d'un point. 
  Un tel faisceau est donc obligatoirement $\mu-$stable.\end{proof}

Soit
  l'ensemble $D$ des  points de $M_{n}$ repr\'esentant des faisceaux
semi-stables
$F$ tels que
$
\lg(\ul{\Ext}^{1}(F,{\cal O}))=1,$ c'est-\`a-dire des faisceaux singuliers en
 exactement un point.  Cet ensemble  $D$  est l'intersection du bord
$\partial M_{n}$  avec un ouvert contenu dans l'ouvert des classes
de faisceaux
$\mu-$stables. Si le faisceau $F$  est singulier en exactement un point
 le bidual $E= F^{**}$ est encore $\mu-$stable et  localement
libre; 
 le support du faisceau
$\ul{\Ext}^1(F,{\cal O})$ d\'efinit un point $x$. On obtient un
morphisme 
$$\pi:D\hfl{}{} U_{n-1} \times \proj_2$$
d\'efini par $F\mapsto (F^{**}, x).$
C'est un morphisme surjectif et lisse dont les fibres sont des  droites
projectives.  La courbe $\beta_{F}$ des droites de saut de $F$  est alors
la r\'eunion (que nous \'ecrirons somme, dans la terminologie usuelle des
diviseurs) 
 de
la courbe des droites de saut de $E=F^{**}$  et de la droite $\check{x}$
de $\proj_2^*$  d\'efinie par le point $x.$  On a donc un diagramme commutatif
$$\eqalign{\diagram{D&\hookrightarrow& M_{n}\cr
\vfl{\pi}{}&&\vfl{}{\beta} \cr
 U_{n-1}\times \proj_2&\hfl{\psi}{}& {\cal C}_{n}\cr 
}}\eqno(10.1)\quad$$
o\`u le morphisme $\psi$  associe au couple $(E,x)$ la courbe $\beta_{E}
+\check{x}.$

\begin{prop} \label{P1}
Il existe un ouvert non vide ${\cal C}_{n}^*$ de
${\cal C}_{n}$ satisfaisant aux conditions suivantes

{\rm (a)} l'ouvert $\beta^{-1}({\cal C}_{n}^*)$  rencontre $D$  et ne
contient pas de point  singulier;

{\rm (b)}  tout point de $M_{n}$  repr\'esentant un faisceau $F$ dont la
courbe des droites de saut $\beta_{F}$  est une courbe  r\'eductible de ${\cal
C}_{n}^*$ appartient \`a $D.$

\end{prop}

\begin{cor}\label{cartesien}
Au-dessus d'un ouvert convenable  de
${\cal C}_{n}$ rencontrant l'image de $\psi$,   le diagramme (10.1)
ci-dessus est cart\'esien (ensemblistement).
\end{cor}
\begin{proof}  
On peut trouver  un ouvert ${\cal
C}_{n}^*\subset {\cal C}_{n}$  satisfaisant aux conditions (a)
et (b) ci-dessus et satisfaisant en outre aux conditions suivantes:

{\rm (c)}  les points de ${\cal C}^*_{n}$ repr\'esentent des
courbes $\hbox{\goth{q}}$ qui
 n'ont pas de composante 
de degr\'e $i$  avec
$1<i<n-1;$

{\rm (d)} le morphisme $\psi:U_{n-1}\times
\proj_2\hfl{}{}{\cal C}_{n}$  d\'efini par  $(E,x)\mapsto
\beta_{E}+\check{x}$ est un plongement ferm\'e au-dessus de l'ouvert 
${\cal C}_{n}^{*}.$

Trouver un  ouvert satisfaisant aux  conditions (a)-(c) est facile 
parce que l'ouvert d\'ej\`a construit dans  la proposition \ref{P1}
rencontre l'ouvert des courbes satisfaisant \`a la condition (c): il suffit
en effet de remarquer qu'il existe un point de
$U_{n-1}$  repr\'esentant un fibr\'e $E$ dont la courbe des droites de saut
est irr\'eductible. Ceci r\'esulte d\'ej\`a du r\'esultat de Barth qui
dit que l'image de $\beta$  rencontre l'ouvert des courbes lisses,
r\'esultat qu'on retrouve facilement, comme on l'a d\'ej\`a vu,  \`a partir de la
proposition
\ref{points-singuliers}.  Pour r\'ealiser la condition (c), il suffit
alors d'enlever \`a cet ouvert le ferm\'e des points repr\'esentant des
courbes dont une des composantes est de degr\'e $i$, o\`u $1<i<n-1.$ Pour
r\'ealiser la condition (d)  il suffit d'appliquer l'hypoth\`ese de
r\'ecurrence, compte-tenu du fait que l'application ${\cal C}_{n-1}\times
\proj_2\hfl{}{} {\cal C}_{n}$  est un plongement propre au-dessus de 
l'ouvert de ${\cal C}_{n}$ des courbes  satisfaisant \`a la   
condition (c).  
Le corollaire r\'esulte alors trivialement de la proposition
\ref{P1}.
\end{proof}

On verra m\^eme qu'au dessus d'un ouvert convenable de ${\cal C}_{n}^*$ ce
diagramme est cart\'esien, au sens des sch\'emas.  Cet \'enonc\'e entra\^{\i}ne
compte-tenu de l'hypoth\`ese de r\'ecurrence, que les fibres au-dessus d'un
ouvert convenable de $B$ sont r\'eduites \`a des droites projectives. La
d\'emonstration de la	 proposition
\ref{P1} utilise le r\'esultat suivant d\^u \`a Str\o mme [23].
\begin{lemme}  \label{Stromme} L'ensemble des points de $U_{n}$ 
correspondant  aux fibr\'es dont la courbe des droites de saut contient
une droite est un ferm\'e de $U_{n}$ de codimension $\geq n-1.$
\end{lemme}

{\it D\'emonstration de la proposition  \ref{P1}}

Soit  $R$  le ferm\'e  de ${\cal C}_{n}$  des points repr\'esentant des courbes 
dont une des composantes est une droite. Le lemme de Str\o mme montre que 
l'image r\'eciproque
$\beta^{-1}(R)$ rencontre $U_{n}$ en codimension
$\geq n-1.$  Par suite  l'image par $\beta$ de $A'=\beta^{-1}(R)\cap U_{n}$ 
est contenue dans un ferm\'e de dimension $\leq 3n-2.$ De plus, l'image du ferm\'e
$A''$ de $M_{n}$ correspondant aux faisceaux
$F$ tels que
$\lg(\ul{\Ext}^1(F,{\cal O}))\geq 2$  est contenue dans le ferm\'e des courbes
qui s'\'ecrivant $\beta_{G} + C$ o\`u $G$  est un
faisceau  de rang 2, $\mu-$semi-stable, de classe de
Chern
$(0,n-2)$  et
$C$  une conique  singuli\`ere.
  L'ensemble
de ces courbes  est contenu dans un ferm\'e de dimension
$\leq 4n-7.$  L'ouvert  $${\cal C}_{n}^{*}= {\cal
C}_{n}-\beta(\overline{A'}\cup A'')$$ 
satisfait aux conditions demand\'ees. \cqfd

\begin{lemme} \label{bord}
Le long des fibres de $\pi$  le faisceau inversible ${\cal O}(\partial M_{n})$ 
est isomorphe \`a ${\cal O}_{\proj_1}(-2).$

\end{lemme}

\begin{proof}
On sait que le groupe de Picard de $M_{n}$  est $\entr^2$; on le d\'ecrit en
consid\'erant l'orthogonal $c^{\bot}$ de $c=2-nh^2$ dans le groupe de
Grothendieck $K(\proj_2)$ pour la forme quadratique
$u\mapsto \chi(u^2)$ et le morphisme   naturel 
$$\lambda_{M_{n}}: c^{\bot}\hfl{}{} \Pic(M_{n})$$
d\'efini par le d\'eterminant de la cohomologie [16].
On sait que cet homomorphisme  est un isomorphisme.
Consid\'erons le
sous-groupe de  base ${\cal D}$  et ${\cal B}$ 
o\`u ${\cal D}$  est le fibr\'e d\'eterminant de Donaldson, et ${\cal B}$ le
fibr\'e d\'eterminant, d\'efini par
${\cal B}=-\lambda_{M_{n}}(u)$
o\`u $u\in K(\proj_2)$  est la classe de Grothendieck  d\'efinie par
$u=2+2h+(n-4)h^2.$ Ce sous-groupe est d'indice 2 ou 1  suivant  que
$n$ est pair ou impair. Dans
$\Pic(M_{n})$, on a, en notations additives
$\partial M_{n}=5{\cal D}-{\cal B}$ (\cf [17]). Fixons un  point de
$U_{n-1}$  correspondant \`a un   faisceau localement libre stable
$E$  et un point $x
\in
\proj_{2}.$ Pour calculer ${\cal B}$ le long de la droite
$\proj_{1}=\proj(E_{x})$, on introduit la famille 
$F$ param\'etr\'ee par
$\proj(E_{x})$  d\'efinie par  la  suite exacte  canonique  sur $\proj(E_{x})
\times \proj_2$
$$0\hfl{}{} F\hfl{}{}{\cal O}\boxtimes E \hfl{}{} {\cal O}(1)
\boxtimes k(x)\hfl{}{} 0 \eqno(10.2)\quad$$
o\`u $k(x)$  est le faisceau structural du point $x.$
La famille $F$  ainsi obtenue fournit par la propri\'et\'e de
module grossier un morphisme   $f_{F}:\proj(E_{x})\hfl{}{} M_{n}$  qui
identifie la fibre \`a la droite projective $\proj(E_{x}).$ 
 De la suite exacte ci-dessus, et de la
caract\'erisation de ${\cal B}$,  il en d\'ecoule, en introduisant  les
projections  naturelles $\pr_{1}$ et $\pr_{2}$ sur $\proj(E_{x})$  et
$\proj_2$  respectivement
$$f_{F}^*({\cal B})=-\lambda_{M_{n}}(u)=
\det
\pr_{1!}({\cal O}(1) \boxtimes k(x).\pr_{2}^*(u))$$
et par suite  ${\cal B}= {\cal O}(2)$. Puisque
${\cal D}$ est trivial
 sur  les fibres de $\beta,$
  ceci conduit au 
r\'esultat.\end{proof}
\goodbreak
\section{Construction d'une surface $\Sigma$ rencontrant $\partial
M_{n}$  tranversalement le long d'une fibre de $\pi.$}
\subsection{Faisceaux de Hulsbergen}

\begin{defi} Un faisceau semi-stable $F$  de rang 2 et classes de Chern
$(0,n)$  est appel\'e faisceau de Hulsbergen si $\dim H^{0}(F(1))\geq 1.$
\end{defi}

On consid\`ere  le sch\'ema universel $\Xi \subset \Hilb^{n+1}(\proj_2)
\times \proj_2$ et le diagramme
$$\diagram{\Xi &\hfl{\pr_{2}}{}& \proj_2\cr
\vfl{\pr_{1}}{}&&\cr
\Hilb^{n+1}(\proj_2)&&\cr}$$
L'image directe ${\cal R}= R^{1}\pr_{1^*}({\cal I}_{\Xi}(-1))=
\pr_{1^*}({\cal O}_{\Xi}(-1))$ est un faisceau localement libre de
rang $n+1$ sur le sch\'ema de Hilbert. L'espace total du fibr\'e en espaces 
projectifs (au  sens de Grothendieck)
$S=\proj({\cal R)}$  param\`etre les syst\`emes coh\'erents
semi-stables $(\Gamma,F(1)),$  o\`u
$F$  est un faisceau  de rang 2, de classes de Chern $(0,n)$  et
$\Gamma\subset H^0(F)$ un sous-espace vectoriel de dimension 1, pour une
notion  ad\'equate de semi-stabilit\'e [17]. Consid\'erons le fibr\'e quotient
canonique ${\cal O}_{\proj({\cal R})}(1)$ sur $\proj({\cal R}).$ Il
existe sur
$S \times \proj_2$ une extension universelle 
$$0\hfl{}{} {\cal O}(1,-1)\hfl{}{} F \hfl{}{} {\cal
I}_{\Xi_{S}}(0,1)\hfl{}{}0. \eqno(11.1)\quad$$ o\`u $\Xi_{S}$  est le
sous-sch\'ema de codimension 2 obtenu \`a partir de $\Xi$  par changement de
base. Cette famille 
$F$ d\'efinit une application rationnelle $S \rightsquigarrow M_{n}.$ 
L'adh\'erence de l'image est constitu\'ee des points de $M_{n}$ qui peuvent 
se repr\'esenter par un faisceau de Hulsbergen; ce ferm\'e est
irr\'eductible de dimension $3n+2$ (\cf [17]).
\subsection{Famille \`a un param\`etre de faisceaux de Hulsbergen de $D$
associ\'ee \`a un couple $(E,x_0) \in P_{n-1} \times \proj_2$ g\'en\'eral.}

Consid\'erons un fibr\'e de Poncelet $E$ fix\'e, de classes de Chern
$(0,n-1)$ associ\'e \`a un syst\`eme lin\'eaire  
$(\Gamma,\Theta)$ sans point de base   sur une conique lisse
$C$  de $\proj_2$.  Le faisceau $\Theta	$  est  un
faisceau inversible de degr\'e $n$  sur la conique $C.$ On a
donc  une suite exacte 
$$0\hfl{}{}\Gamma^* \otimes {\cal O}\hfl{}{} E(1)\hfl{}{}
\check{\Theta}\hfl{}{}0 \eqno (11.2)\quad$$ 
de sorte que si $x_{0}$  est un point choisi en dehors de
$C$ la fibre $E_{x_{0}}$ de $E$  s'identifie \`a
$\Gamma^*.$  Consid\'erons la droite projective
$\proj=\proj(\Gamma^*)=\proj(E_{x_0}).$ Sur $\proj \times \proj_2,$ 
la suite exacte (10.2)  d\'efinit une famille plate de faisceaux
 param\'etrant la fibre de $\pi:D\hfl{}{} U_{n-1}\times
\proj_2$  au-dessus du point
$(E,x_{0}).$  Ces faisceaux sont \'evidemment des faisceaux  de
Hulsbergen.

\begin{lemme}  Le faisceau  $F$  d\'efinit une famille de faisceaux de
Hulsbergen de $M_{n}$ qui provient d'un morphisme $\varphi:\proj\hfl{}{}
\proj({\cal R}).$ 
\end{lemme}

\begin{proof}
Consid\'erons le 
sous-sch\'ema
$\Delta$ de
$\proj_1
\times C$, plat et fini au-dessus de $\proj_{1},$  d\'efini par le
sch\'ema des z\'eros de la section canonique de  ${\cal O}(1)\boxtimes
\Theta$ associ\'ee au syst\`eme lin\'eaire. Le choix  d'un g\'en\'erateur de
$\wedge^2
\Gamma^*$ fournit un isomorphisme $\Gamma\simeq \Gamma^*$;  on dispose
alors  de la   suite exacte d'Euler sur
$\proj$ 
$$0\hfl{}{}{\cal
O}_{\proj}(-1)\hfl{}{}
{\cal O}_{\proj}\otimes\Gamma^
*\hfl{}{} {\cal O}_{\proj}(1)\hfl{}{}0$$ 
et par suite d'un morphisme injectif
${\cal O}_{\proj}(-1)\boxtimes {\cal O}\hookrightarrow {\cal O}\boxtimes
E(1)
$ dont le conoyau est isomorphe \`a ${\cal I}_{\Delta}(1,2).$ Consid\'erons
l'id\'eal
${\cal I}_{Z} $ du  sous-sch\'ema  de 
$Z= \Delta\cup
 (\proj\times\{x_{0}\})$; ce sous-sch\'ema de dimension 
relative 
$0$  et de degr\'e $n+1$ au-dessus de la droite projective $\proj$. On en
d\'eduit le diagramme commutatif de suites exactes au-dessus de $\proj\times
\proj_2:$
$$\diagram{
&&&&0&&0&&\cr
&&&&\vfl{}{}&&\vfl{}{}&&\cr
0&\hfl{}{}&{\cal O}(-1,-1)&\hfl{}{}& F&\hfl{}{}& {\cal
I}_{Z}(1,1)&\hfl{}{}&0
\cr 
&&||&&\vfl{}{}&&\vfl{}{}&&\cr
0&\hfl{}{}&{\cal O}(-1,-1)&\hfl{}{}&{\cal O}\boxtimes E&\hfl{}{}& {\cal
I}_{\Delta}(1,1)&\hfl{}{}&0\cr
&&&&\vfl{}{}&&\vfl{}{}&&\cr
&&&&{\cal O}_{\proj}(1)\boxtimes k(x_0)&\simeq&{\cal
O}_{\proj}(1)\boxtimes k(x_0)&&\cr &&&&\vfl{}{}&&\vfl{}{}&	&\cr
&&&&0&&0&&\cr
}\eqno(11.3)\quad$$

 La
propri\'et\'e universelle du sch\'ema de Hilbert  fournit un plongement 
$\psi:\proj\hfl{}{} \Hilb^{n+1}(\proj_2)$  
tel que $(\psi\times \id)^{-1}(\Xi)= Z. $  Consid\'erons la suite exacte 
$$0\hfl{}{}{\cal O}(-1,-1)\hfl{}{} F\hfl{}{} {\cal
I}_{Z}(1,1)\hfl{}{}0
 \eqno \quad (11.4)\quad $$ figurant dans le diagramme ci-dessus. Si on
restreint cette suite exacte \`a la fibre $\{\eta\} \times \proj_2$
au-dessus du point
$\eta \in
\proj$ on obtient encore par platitude une suite exacte, et  cette suite
exacte n'est pas  scindable puisque le faisceau associ\'e  $F_{\eta}$ n'a
qu'un point singulier. Ainsi, l'extension (11.4) d\'efinit  un
morphisme
$\varphi:\proj
\hfl{}{}
\proj({\cal R})$ au-dessus de $\psi$ caract\'eris\'e par le fait que
l'extension  obtenue \`a partir de  (11.1)  par le changement de base
$\varphi$  est l'extension (11.4).
\end{proof}

 On se propose d'\'etendre le
morphisme
$\varphi$  en un plongement $\phi: \Sigma\hfl{}{} \proj({\cal R})$ d'une
surface d'Hirzebruch  au-dessus de $\proj.$  Commen\c{c}ons par l'\'etude de
$\psi^*({\cal R}):$

\begin{lemme}  Soit $\Delta$  le sch\'ema des z\'eros de la section canonique
de ${\cal O}_{\proj}(1)\boxtimes \Theta,$ et $t\in H^{0}({\cal
O}_{\proj_2}(1))$ une forme lin\'eaire ne s'annulant pas en $x_{0}.$

{\rm (i)} Le choix de $t$ induit un isomorphisme $$\psi^*({\cal R}) = 
\pr_{1*}({\cal O}_{\Delta}(-1))
\oplus {\cal O}_{\proj_{1}}.$$

{\rm (ii)} Soit $\Gamma^{\bot}\subset 
H^{0}(\Theta)^*$ l'orthogonal de $\Gamma.$  Le
faisceau
$\pr_{1*}({\cal O}_{\Delta}(-1))$  poss\`ede un faisceau inversible quotient $L$ isomorphe \`a ${\cal O}_{\proj}(-2)$ et un
seul;  le noyau de la surjection 
$$\pr_{1*}({\cal O}_{\Delta}(-1))\hfl{j}{} L$$
est canoniquement isomorphe \`a $\Gamma^{\bot}\otimes {\cal
O}_{\proj}(-1)\simeq {\cal O}_{\proj}(-1)^{n-1}.$

{\rm (iii)}  Dans l'isomorphisme $\Ext^1({\cal I}_{Z}(0,2),L)\simeq
\Hom(\psi^*({\cal R}),L)$
l'\'el\'ement de $\Ext^1({\cal I}_{Z}(0,2),L)$ sur
$\proj\times
\proj_2$ correspondant \`a la projection
$\psi^*({\cal R})\hfl{}{}L$ est \`a homoth\'etie pr\`es 
l'extension associ\'ee \`a la suite exacte 
$$0\hfl{}{}{\cal O}(-1,-1)\hfl{}{} F\hfl{}{} {\cal
I}_{Z}(1,1)\hfl{}{}0
 $$
figurant dans le diagramme {\rm (11.3)}. En particulier,
$\varphi^*({\cal O}_{\proj({\cal R})}(1))\simeq{\cal O}_{\proj}(-2).$
\end{lemme}

\begin{proof}
On a \'evidemment ${\cal R}\simeq \pr_{1*}({\cal O}_{\Xi}(-1))$  et par 
changement de base $f^*({\cal R})\simeq \pr_{1*}({\cal O}_{Z}(-1)).$  Il en
d\'ecoule (i). De la  r\'esolution 
$$0\hfl{}{}{\cal O}(-1)\boxtimes \Theta^*\hfl{}{} {\cal O}_{\proj\times C}
\hfl{}{} {\cal O}_{\Delta}\hfl{}{} 0 $$
on tire  la suite exacte  $$0\hfl{}{}\pr_{1*}({\cal
O}_{\Delta}(-1))\hfl{}{} R^1\pr_{1^*}({\cal O}_{\proj}(-1)\boxtimes
\Theta^*(-1))\hfl{}{} R^1\pr_{1*}({\cal O}_{\proj \times C}(0,-1))\hfl{}{}
0$$ ce qui montre apr\`es dualit\'e de Serre que
${\cal N}=\pr_{1*}({\cal O}_{\Delta}(-1))$ est le noyau du morphisme
naturel
${\cal O}_{\proj}(-1)
\otimes H^{0}(\Theta)^*\hfl{}{} {\cal O}_{\proj}$  induit par la
projection canonique
$H^0(\Theta)^*\otimes {\cal O}_{\proj}\hfl{}{} {\cal
O}_{\proj}(1).$ La filtration de Harder-Narasimhan de ${\cal N}$ fournit 
la suite exacte attendue: le quotient est le noyau du morphisme
$\Gamma^*\otimes {\cal O}_{\proj} (-1)\hfl{}{}{\cal O}$, lequel
s'identifie, apr\`es avoir choisi un g\'en\'erateur de $\wedge^2\Gamma^*$, \`a
${\cal O}_{\proj}(-2)$ ce qui d\'emontre (ii).   L'isomorphisme $\Ext^1({\cal
I}_{Z}(0,2),L)\simeq
\Hom(f^*({\cal R}),L)$  r\'esulte du th\'eor\`eme de dualit\'e de Serre relatif
pour la projection $\pr_1:\proj\times \proj_2\hfl{}{}\proj.$ Cet
isomorphisme prouve qu'il n'y a qu'une extension non scind\'ee \`a isomorphisme
pr\`es. Il en d\'ecoule l'assertion (iii). \end{proof}

\goodbreak
\subsection{La surface de Hirzebruch $\Sigma$}

Consid\'erons le fibr\'e vectoriel de rang 2 sur $\proj$  d\'efini par 
$W= {\cal O}_{\proj}(-2) \oplus {\cal O}_{\proj}$. C'est un  fibr\'e vectoriel
quotient de rang 2  de
$\psi^*({\cal R})$  et par cons\'equent  la surface de Hirzebruch
$\Sigma=\proj(W)$ 
 se plonge dans $\proj({\cal R}),$ de sorte que l'on a un diagramme
commutatif 
$$\diagram{
\Sigma&\hfl{\phi}{}&\proj({\cal R})\cr
\vfl{\pi}{}&\nefl{\varphi}{}&\vfl{}{\pi}\cr
\proj&\hfl{\psi}{}&\Hilb^{n+1}(\proj_2)\cr
}$$
o\`u $\phi$  est le plongement ci-dessus.
La surface $\Sigma$ param\`etre une extension  obtenue \`a partir de
l'extension (11.1)  par changement de base, et donc une famille
${\cal F}$ de faisceaux $\mu-$semi-stables de rang 2 et classe de Chern
$(0,n).$  Plus pr\'ecis\'ement:

\begin{prop} Le faisceau ${\cal F}$   d\'efinit  une famille plate de
faisceaux stables param\'etr\'ee par la surface de Hirzebruch $\Sigma$.
\end{prop}

\begin{proof}
La surface $\Sigma$ a deux sections canoniques $\Sigma_{0}$  et
$\Sigma_{\infty}$  correspondant aux quotients ${\cal O}(-2)$  et ${\cal O}.$ 
 Le long de la section
$\Sigma_{0}$, le morphisme compos\'e  $\proj \simeq
\Sigma_{0}\hfl{}{}\proj({\cal R})$ co\"{\i}ncide avec le morphisme $\varphi$ 
consid\'er\'e ci-dessus.  
 D\'esignons  par $\Sigma^*$  l'ouvert compl\'ementaire de ces deux
sections, et  par $Z_{\Sigma}$  le sous-sch\'ema de $\Sigma\times \proj_2$ 
obtenu \`a partir de $Z\subset \proj\times
\proj_2$
  par le changement de base $\Sigma\hfl{}{}\proj.$
 La projection $\phi^*\pi^*({\cal R})\hfl{}{}{\cal O}_{\proj}(1) $ 
induit une section $\sigma$ du faisceau inversible le
long de $Z_{\Sigma}$  d\'efini par $\ul{\Ext}^1({\cal
I}_{Z_{\Sigma}}(0,2),\pr_{1}^*({\cal O}_{\Sigma}(1))).$ Ce faisceau
inversible est isomorphe \`a
$\ul{\Ext}^2({\cal O}_{Z_{\Sigma}}(2),\pr_{1}^*({\cal O}_{\Sigma}(1)))$  et 
 s'obtient donc, apr\`es tensorisation par ${\cal O}_{\Sigma}(1)$, par le
changement de base
$\Sigma\hfl{}{}
\proj$ \`a partir du faisceau inversible sur $Z$
$$\eqalign{\ul{\Ext}^{2}({\cal O}_{\Delta},{\cal O}(0,-2))\oplus
\ul{\Ext}^{2}({\cal O}_{\proj
\times
\{x_{0}\}},{\cal O}(0,-2))
&= ({\cal O}
(1)\boxtimes \Theta)|_{\Delta} \oplus {\cal O}_{\proj \times
\{x_{0}\}}
\cr}$$
Compte-tenu du fait que le faisceau ${\cal O}(-1)\boxtimes \Theta$  est
trivial sur
$\Delta$, ce faisceau s'identifie aussi \`a ${\cal O}(2,0)|_{\Delta}
\oplus {\cal O}_{\proj
\times
\{x_{0}\}}
.$ Il en r\'esulte que  cette section
ne s'annule pas sur  l'ouvert
$\Sigma^*\times
\proj_2,$  
et donc le faisceau
${\cal F}|_{\Sigma^{*}\times\proj_2^*}$  est  localement libre. Pour
$\eta\in \Sigma^*$, le faisceau ${\cal
F}_{\eta}$ est donc un faisceau  localement libre  de rang 2, de classes
de Chern $(0,n)$  qui  n'a pas de sections;
il en r\'esulte que c'est   un fibr\'e stable.   Reste l'\'etude de ${\cal F}$
le long de la section
$\Sigma_{\infty}.
$
\begin{lemme}  Pour $\eta\in \Sigma_{\infty},$  le faisceau ${\cal F}_{\eta}$ 
est un faisceau stable singulier aux points du diviseur 
$\Delta_{\eta}$ de $ C.$
\end{lemme}
\begin{proof}  De la suite exacte  $0\hfl{}{}{\cal I}_{Z_{\eta}}
\hfl{}{}  {\cal I}_{\Delta_{\eta}}\oplus {\cal
I}_{{\{x}_0\}}\hfl{}{}{\cal O}\hfl{}{}0$  sur $\proj_2$ d\'ecoule
l'isomorphisme 
$$\Ext^1({\cal I}_{Z_{\eta}}(1),{\cal O}(-1))\simeq \Ext^1({\cal
I}_{\Delta_{\eta}}(1),{\cal O}(-1))\oplus \Ext^1({\cal
I}_{\{x_{0}\}(1)},{\cal O}(-1)).$$  Dire  que $\eta$  appartient \`a la
section $\Sigma_{\infty}$  signifie  que l'extension provient de 
$\Ext^1({\cal I}_{\{x_{0}\}}(1),{\cal
O}(-1))$  dans cet isomorphisme.  Le faisceau $G$ associ\'e \`a une telle 
extension est alors le faisceau trivial  de rang 2. Une section $t$  de
${\cal O}(1)$ ne s'annulant pas sur $\Delta_{\eta}$ fournit un
morphisme surjectif
$G\hfl{}{}{\cal O}_{\Delta_{\eta}}$ dont ${\cal F}_{\eta}$  est le noyau, et on obtient le
diagramme commutatif de suites exactes: 
$$\diagram{
&&&&0&&0&&\cr
 &&&&\vfl{}{}&&\vfl{}{}&&\cr
0&\hfl{}{}& {\cal O}(-1)&\hfl{}{}& {\cal
F}_{\eta}&\hfl{}{} &{\cal I}_{Z_{\eta}}(1)&\hfl{}{}& 0\cr
&&\vfl{}{\wr}&&\vfl{}{}&&\vfl{}{}&&\cr
0&\hfl{}{}& {\cal O}(-1) &\hfl{}{} &G&\hfl{}{} &{\cal
I}_{\{x_{0}\}}(1)&\hfl{}{}& 0 \cr
&&&&\vfl{}{}&&\vfl{}{}&&\cr
&&&&{\cal O}_{\Delta_{\eta}}&\hfl{\sim}{t} &{\cal
O}_{\Delta_{\eta}}(1)&\cr &&&&\vfl{}{}&&\vfl{}{}\cr
&&&&0&&0\cr
}$$
Le choix d'un isomorphisme $G \simeq {\cal O}^2$  d\'etermine des formes 
lin\'eaires $(x,y)$ ind\'ependantes s'annulant au point $x_{0}.$ 

Consid\'erons un sous-faisceau $F'$  de ${\cal F}_{\eta}$  tel que le quotient
$F''={\cal F}_{\eta}/F'$  soit sans torsion.  Alors $c_{1}(F')\leq 1$  et
en cas d'\'egalit\'e, la relation $c_{2}(F')+c_{2}(F'')=n+1$ montre que 
$c_{2}(F')\leq n+1.$  Si $c_{1}(F')=1$, le morphisme induit
$F'\hfl{}{}{\cal I}_{Z_{\eta}}(1)$  serait  non nul, et donc $c_{2}\geq n+1.$ 
Alors ce morphisme serait un isomorphisme, et la premi\`ere ligne serait
scind\'ee, ce qui est absurde parce que ${\cal F}_{\eta}$  est
\'evidemment  localement libre  au voisinage de $\{x_{0}\}.$ 
Donc $c_{1}(F')\leq 0.$  

Supposons  que $c_{1}(F')=0.$ Le faisceau $F'$  est
l'intersection de ${\cal F}_{\eta}$ avec un sous-faisceau ${\cal O
}\hookrightarrow{\cal O}^2,$ donn\'e par un vecteur non nul $(\lambda,\mu)\in
k^2.$  Ainsi,  le faisceau
$F'$ est le noyau du morphisme ${\cal O}\hfl{\lambda x+\mu y}{} {\cal
O}_{\Delta_{\eta}}(1)$. Par suite,
$$\eqalign{c_{2}(F')&= \lg (\ker {\cal O}_{\Delta_{\eta}}\hfl{\lambda x+\mu
y}{}{\cal O}_{\Delta_{\eta}}(1)) \cr 
&\geq n-2}$$
Ainsi, l'hypoth\`ese $n\geq 5$  entra\^{\i}ne $c_{2}(F')>{n\over 2}.$   Donc le
faisceau
${\cal F}_{\eta}$  est stable. \end{proof}

La propri\'et\'e de module grossier de $M_{n}$  permet d'associer \`a la famille
${\cal F}$  un morphisme $f:\Sigma\hfl{}{} M_{n}$  qui ne rencontre la
fronti\`ere $\partial M_{n}$ que le long des sections $\Sigma_{0}$  et
$\Sigma_{\infty}.$  Plus pr\'ecis\'ement, l'image r\'eciproque de $D$  est, au moins
ensemblistement, la section $\Sigma_{0}.$  L'\'enonc\'e suivant montre que
c'est l'image r\'eciproque au sens des sch\'emas:

\begin{prop}\label{transverse}  Le morphisme $f:\Sigma\hfl{}{} M_{n}$ 
est transverse \`a $D$  le long de $\Sigma_{0}.$
\end{prop}

\begin{proof} Pour toute famille plate 
$({\cal E}_{s})_{s \in S}$  de faisceaux stables  de
$M_{n}$, param\'etr\'ee par une vari\'et\'e alg\'ebrique  lisse $S$,   l'image
r\'eciproque du diviseur 
$\partial M_{n}$  est  le sous-sch\'ema  d\'efini par l'id\'eal de Fitting du
faisceau
$\pr_{1^*}(\ul{
\Ext}^1({\cal E},{\cal O}(0,-1))$ (Str\o mme [23], \S 4). Si on
applique ceci \`a notre faisceau ${\cal F}$  on obtient sur
$\Sigma\times \proj_2$ la pr\'esentation
$${\cal O}(-1,0)\hfl{\omega}{}\ul{\Ext}^2({\cal
O}_{Z_{\Sigma}},{\cal O}(0,-2))\hfl{}{}\ul{\Ext}^1({\cal F},{\cal
O}(0,-1))\hfl{}{} 0
$$
Le faisceau $\ul{\Ext}^2({\cal O}_{Z_{\Sigma}},{\cal O}(0,-2)))$  est
localement libre de rang 1  sur
$Z,$ et le morphisme
$\omega$ s'identifie \`a la section canonique de $\ul{\Ext}^2({\cal
O}_{Z_{\Sigma}},{\cal O}(0,-2))\otimes {\cal O}_{\Sigma}(1,0)$ sur
$\Sigma$ d\'ej\`a utilis\'ee; si $\Delta_{\Sigma}$  d\'esigne  l'image
r\'eciproque de
$\Delta$  dans $\Sigma\times \proj_2,$  on sait que
$$\ul{\Ext}^2({\cal O}_{Z_{\Sigma}},{\cal O}(0,-2))=
{\cal O}_{\Delta_{\Sigma}}(2,0)\oplus {\cal O}_{\Sigma\times \{x_{0}\}}$$
et  sur le compl\'ementaire de $\Sigma_{\infty}$ la premi\`ere   
composante $\omega_{1}$ de
$\omega$ ne s'annule pas.  La deuxi\`eme composante $\omega_{2}$ de $\omega$
est obtenue en restreignant l'image r\'eciproque du morphisme  compos\'e 
${\cal O}_{\Sigma}\hookrightarrow \pi^*(W)\hfl{}{} {\cal
O}_{\Sigma}(1)$ \`a $\Sigma_{0} \times \{x_{0}\}.$ Le compl\'ementaire de
$\Sigma_{\infty}$ est isomorphe
 \`a l'espace total du fibr\'e ${\cal
O}_{\proj}(-2),$ la section $\Sigma_{0}$  s'identifiant \`a la section
nulle, et le fibr\'e
${\cal O}_{\Sigma}(1)$ est isomorphe \`a l'image r\'eciproque de ${\cal
O}_{\proj}(-2)$ sur cet ouvert;  le morphisme $\omega_{2}$ ci-dessus est
donn\'e par l'identification
${\cal O}_{\proj}(-2)\hfl{}{} {\cal O}_{\proj\times
\{x_{0}\}}(-2,0)$ des espaces totaux de ces fibr\'es inversibles.
En dehors de $\Sigma_{\infty}$, la section $\omega_{2}$ ne s'annule 
que sur $\Sigma_{0}$ et ceci transversalement \`a la section nulle.  Par
suite,  le faisceau 
$\pr_{1*}(\ul{\Ext}^1({\cal F},{\cal
O}_{\Sigma \times \proj_2}(0,-1)))$  est  de longueur relative 1  le long de
la section
$\Sigma_{0}.$ Il en d\'ecoule qu'au sens des sch\'emas,
  on a $f^{-1}(\partial M_{n})= \Sigma_{0}.$ On retrouve le fait que
$\partial M_{n}$  est lisse en tout point de l'image de $f;$ de plus, 
puisque
$\Sigma_{0}$ est lisse, le morphisme $f$  est transverse \`a $\partial M_{n}.$
\end{proof}

\goodbreak
\section {Image de $\Sigma$  par le morphisme de Barth}

La famille de faisceaux ${\cal F}$ param\'etr\'ee par la surface  $\Sigma$ est
une famille de faisceaux stables de classes de Chern $(0,n)$. On peut pour
cette famille consid\'erer la courbe des droites de saut. On note $\gamma:
\Sigma \hfl{}{} {\cal C}_{n}$  le morphisme compos\'e $\beta\circ f$ ainsi
obtenu.   
\begin{lemme}
Le long des fibres de
$\Sigma \hfl{}{}\proj$  le morphisme $\gamma$  est
une homographie et donc est injectif.
\end{lemme}
\begin{proof}
Soit $\eta \in \proj,$  et $Z_{\eta}$  la fibre de $Z$  au-dessus de $\eta.$
Consid\'erons, pour $t\in \Sigma_{\eta},$ la suite exacte 
$$0\hfl{}{} {\cal O}(-1)\hfl{}{} {\cal F}_{t}\hfl{}{} {\cal
I}_{Z_{\eta}}(1)\hfl{}{}0$$
 Consid\'erons le diagramme standard
$$\diagram{
D&\hfl{\pr_2}{}& \proj_2\cr
\vfl{\pr_1}{}&&\cr
\proj_2^*\cr&&
}$$ La courbe des droites de saut  $\gamma(t)=\beta_{F}$  du
faisceau
$F={\cal F}_{t}$ est  d\'efinie par 
 l'id\'eal de Fitting du  faisceau
$R^{1}{\pr_{1*}}(\pr_{2}^*(F(-1)).$  On a la suite
exacte
$$0\hfl{}{}\pr_{1*}(\pr_2^*({\cal I}_{Z_{\eta}}))\hfl{}{}{\cal
O}(-1)\hfl{}{} R^{1}{\pr_{1*}}(\pr_{2}^*(F(-1)))\hfl{}{}
R^1\pr_{1*}(\pr_2^*({\cal I}_{Z_{\eta}}))\hfl{}{} 0
$$
Le support de $R^1\pr_{1*}(\pr_2^*({\cal I}_{Z}))$  est form\'e des  points
correspondant aux droites de $\proj_2$ qui rencontrent $Z$ suivant un
sous-sch\'ema de longueur
$\geq 2;$   ce sous-sch\'ema est fini. Il en
d\'ecoule que
$\pr_{1*}(\pr_2^*({\cal I}_{Z}))$ est un faisceau sans torsion de rang 1  et
de classe de Chern $c_{1}= -n-1. $ L'id\'eal de Fitting de
$R^{1}{\pr_{1*}}(\pr_{2}^*(F(-1)))$  co\"{\i}ncide avec $\pr_{1*}(\pr_2^*({\cal
I}_{Z_{\eta}}))(1)$  en dehors d'un ferm\'e de codimension $2$.  Par suite,
l'\'equation de $\beta_{F}$  d\'epend lin\'eairement de la classe  $\omega \in
W_{\eta}^*
\subset \Ext^1({\cal I}_{Z_{\eta}},{\cal O}(-2))$  qui d\'efinit l'extension. 
Cette \'equation est non nulle si $\omega$  est non nulle.  Par suite,
l'application obtenue  
$$\Sigma_{\eta}\hfl{}{} {\cal C}_{n}$$  est
une homographie. \end{proof}

 Pour d\'ecrire l'image de $\Sigma$ par $\gamma$ il
suffit de d\'ecrire  les images des sections $\Sigma_{0}$  et $\Sigma_{\infty}.$
Le long de 
$\Sigma_{0}$, le morphisme $\gamma$  est constant et a pour image le point
$s$ d\'efini par la courbe $\beta_{E}+ \check{x}_{0}.$ Le long de la section
$s_{\infty}:\proj\hfl{}{}\Sigma_{\infty}$  de $\Sigma$, la famille ${\cal F}$
isomorphe  au  noyau d'un morphisme surjectif sur
$\proj\times
\proj_2$
$${\cal O}^2\hfl{}{} {\cal O}_{\Delta}$$ 
associ\'e au diviseur $\Delta.$  La restriction de $\gamma$   \`a
$\proj\simeq\Sigma_{\infty}$ associe \`a un point de
$\proj(\Gamma)$ correspondant \`a un diviseur 
$ \eta = \sum_{a\in C}n_{a}a\in $ 
$\proj$,   la courbe du plan projectif dual 
 d\'ecompos\'ee en $n$  droites d\'efinie par $\gamma(s_{\infty}(\eta))=
\sum_{a\in C} n_{a} \check{a}$.

\begin{lemme} 
Dans  l'espace projectif
$ {\cal C}_{n},$ l'image de la section  $\Sigma_{\infty}\subset \Sigma$   par le
morphisme
$\gamma$ est une conique lisse $Q$ qui ne rencontre chacune des droites
$\gamma(\Sigma_{\eta})$, o\`u $\eta \in \proj,$  qu'en un seul point.

\end{lemme}
\begin{proof}
Consid\'erons sur la conique $C$  un fibr\'e inversible $A$ de degr\'e $1,$  de
sorte que $\Theta \simeq A^{\otimes n}.$ Sur  $\proj(H^0(\Theta)^*)\times
C$  on dispose d'une section canonique de ${\cal O}(1)\boxtimes L$  dont le
sch\'ema des z\'eros est une hypersurface plate et finie de degr\'e $n$  au-dessus
de  l'espace projectif $\proj(H^0(\Theta)^*).$  On en d\'eduit un morphisme
$\proj(H^0(\Theta)^*)\hfl{}{} \Div^{n}(C)$
qui est \'evidemment un isomorphisme.  On a d'autre part un isomorphisme 
$S^{n}C \simeq \Div^{n}(C)$ qui identifie le fibr\'e ${\cal O}(1)$  avec 
le fibr\'e d\'eterminant ${\cal D}_{A}=A\boxtimes \dots \boxtimes
A/\hbox{\goth{S}}_{n}.$
 Dans le plongement
$S^{n}C= \Div^{n}(C)
\hookrightarrow {\cal C}_{n}$
l'image r\'eciproque du fibr\'e ${\cal O}(1)$  est  isomorphe \`a ${\cal O}(2).$
L'image de $\proj$  dans $\Div^{n}(C)$  est une droite projective.
Consid\'erons le plongement $j=\gamma \circ s_{\infty}: \proj
\hookrightarrow {\cal C}_{n}.$ L'application lin\'eaire associ\'ee
$H^0({\cal C}_{n},{\cal O}(1))\hfl{}{} H^{0}(\proj,{\cal O}(2))$  est
obligatoirement surjective, sinon  le morphisme $j: \proj \hfl{}{} {\cal C}_{n}$
se  factoriserait \`a travers un morphisme de degr\'e 2, et ne serait pas
injectif.  Le morphisme $j$  se factorise suivant le diagramme
$$\diagram{
\proj &\hfl{i}{}& \proj(H^{0}(\proj,{\cal O}(2))^*)\cr
& \sefl{}{}&\vfl{}{}\cr
&&{\cal C}_{n}\cr
}$$
o\`u $i$ est le plongement de
Veronese de la droite 
$\proj$ dans le plan projectif. Donc l'image de la droite projective $\proj$
est une conique lisse $Q$.
\end{proof}

Montrons que pour $\eta \in \proj,$ la  conique $Q$ ne rencontre  la droite
$\gamma(\Sigma_{\eta})$  qu'en un seul point: sinon, le point d\'efini par la
courbe $\beta_{E}+\check{x}_{0}$  appartiendrait au plan $\Span(Q)$ de la
conique
$Q$ et il existerait des droites $\gamma(\Sigma_{\eta})$  tangentes \`a $Q.$ 
Or le morphisme  induit par la diff\'erentielle $d\gamma$  sur les fibr\'es
normaux
$$N_{\Sigma_{\infty}/\Sigma}={\cal O}(2)\hfl{}{} N_{Q/{\cal C}_{n}}={\cal O}(2)
\oplus {\cal O}(1)^{N-2}$$
o\`u $N= \dim {\cal C}_{n}$  est non nul est donc injectif. Ceci contredit  le
fait qu'il existe parmi les droites $\gamma(\Sigma_{\eta})$  des droites
tangentes \`a $Q.$

\begin{cor} \label{vrai.cone}{\rm (i)}  Le point $s$  n'appartient pas au
plan projectif
$\Span(Q)$  engendr\'e par la conique $Q.$

{\rm (ii)} L'image de $\Sigma$  par le morphisme $\gamma$  dans l'espace
projectif des courbes
${\cal C}_{n}$  est le c\^one quadratique $R$ de sommet le point
$s=\beta_{E}+\check{x}_{0}$  et de base 
 la conique $Q=\gamma(\Sigma_{\infty}).$

\end{prop}

 Consid\'erons l'espace
projectif $\Span(R)$ de dimension 3 engendr\'e par le c\^one $R.$ L'\'enonc\'e
suivant  est la cl\'e de la d\'emontration:
\begin{lemme} \label{non-ramif}
Il existe un ouvert de Zariski non vide ${\cal U}$ de
$\Grass(2,H^0(\Theta))\times (\proj_2\setminus C)$   tel que pour
$(\Gamma,x_{0})\in {\cal U}$
  les propri\'et\'es suivantes soient satisfaites:

{\rm (i)}  La courbe de saut $\beta_{E}$ du fibr\'e de Poncelet est lisse.

 {\rm (ii)}  Le  sous-espace projectif $\overline{T}_{s}$ tangent en $s$ \`a
 l'image du morphisme
$${\cal C}_{n-1} \times \proj_{2}\hfl{}{} {\cal C}_{n}$$ d\'efini par
l'addition des diviseurs de $\proj_2^*$ rencontre  l'espace projectif
$\Span(R)$  au seul point
$s$ d\'efini par le sommet du c\^one $R$.
\end{lemme}

\begin{proof}   D'apr\`es la proposition
\ref{points-singuliers}, si
$\Gamma\subset H^0(\Theta)$ est sans point de base   et  tel que  
$z^{2}H^0(\Theta(-2)) \cap \Gamma=\{0\}$ pour toute droite d'\'equation
$z=0$ la courbe de Poncelet associ\'ee \`a
$(\Gamma,\Theta)$  est lisse.  Donc la condition (i) est
satisfaite. 
  Elle implique que  le morphisme  ${\cal C}_{n-1} \times \proj_{2}\hfl{}{}
{\cal C}_{n}$  est 
non ramifi\'e en $(\beta_{E},\check{x}_{0})$: l'espace tangent projectif
$\overline{T}_{s}$ en
$s$ \`a l'image est alors l'espace projectif des courbes de ${\cal C}_{n}$  qui
contiennent l'intersection
$\beta_{E}\cap \check{x}_{0}.$ V\'erifier la condition (ii) revient \`a
v\'erifier que le plan projectif $\Span(Q)$ engendr\'e par la conique $Q$  ne
rencontre pas $\overline{T}_{s}.$

On d\'esigne pour $z \in \proj_2^*$  par $H_{z}\subset {\cal C}_{n}$ 
l'hyperplan projectif des courbes de degr\'e $n$ de $\proj_2^*$ qui passent
par
$z.$   Tout point
$z
\in
\beta_{E}$  d\'efinit une droite de $\proj_{2}$, dont on
d\'esigne aussi par $z$ l'\'equation, telle que
$\Gamma
\cap z H^0(\Theta(-1)) \not=\{0\}. $
L'intersection $H_{z}\cap  Q$  est r\'eduite au seul point  d\'efini par la
courbe associ\'ee au diviseur
$(s)$ d'une section non nulle  $s \in\Gamma \cap zH^{0}(\Theta(-1)).$  En
effet,  un point 
de $H_{z}\cap  Q$ d\'efinit un diviseur de $C$  associ\'e \`a une section $t
\in \Gamma$ s'\'ecrivant $(t)= \sum_{i}a_{i}$, tel que $a_{1} \in z.$
Les sections $s$  et $t$  s'annulent alors en $a_{1}$, et puisque le
syst\`eme lin\'eaire $\Gamma$  est sans point de base,
les sections $s$ et $t$  sont proportionnelles.
Il en r\'esulte que l'hyperplan $H_{z}$ rencontre le plan $\Span Q$  suivant la
 droite tangente \`a la conique $Q$  d\'etermin\'ee par le diviseur $(s).$ 

Soit $\check{Q}$ la conique duale de $Q$. Consid\'erons le morphisme 
$$\tau:\beta_{E}\hfl{}{} \check{Q}$$
qui associe au point $z$  la tangente $H_{z}\cap \Span(Q).$  Dans
l'identification $\check{Q} \simeq \proj(\Gamma)$, ce  morphisme associe \`a
$z$  l'unique section $s \in \proj(\Gamma)$  qui s'annule sur $z\cap C.$
C'est un morphisme fini  de degr\'e ${1\over 2}n(n-1).$

\begin{lemme}\label{bertini}
Si $x_{0}$  est g\'en\'erique,
la restriction  de  $\tau$ \`a $\beta_{E}\cap \check{x}_0$  est injective.
\end{lemme}

Ce lemme sera d\'emontr\'e ci-dessous.
Supposons $x_{0}$ suffisamment
g\'en\'eral pour  que la propri\'et\'e du lemme \ref{bertini} ci-dessus soit vraie
et que la droite $\check{x}_{0}$ ne soit pas tangente \`a $\beta_{E}$. Si
$\overline{T}_{s}$  rencontre
$\Span(Q)$ en un point
$w,$  le point
$w$  d\'efinit une courbe de degr\'e $n$  qui contient $\beta_{E}\cap
\check{x}_{0}.$  Alors $w $ appartient aux tangentes $  H_{z}
\cap
\Span(Q)$, pour tout $z \in \beta_{E}\cap \check{x}_{0}$. D'apr\`es le
choix de  $x_{0}$ ces tangentes  sont distinctes: le point
$w$ appartient donc \`a $n-1$  tangentes \`a la conique $Q$, ce qui est absurde
puisque $n\geq 5.$ \cqfd
\medskip

{D\'emonstration du lemme \ref{bertini}.}

Consid\'erons l'ensemble $K$
des couples $(s,\ell) \in \proj(\Gamma)\times \proj_2^*$  telles que
$\lg(\ell
\cap
\tau^{-1}(s))\geq 2.$ Par projection sur $\proj(\Gamma)$ on voit que c'est un ferm\'e
de dimension
$1$  et par cons\'equent l'image de $K$ dans $\proj_2^*$  est de dimension
1.  Si on prend
$\check{x}_{0}$ en dehors de ce ferm\'e, la restriction de $\tau$ \`a
$\beta_{E}\cap \check{x}_{0}$  est injective.  \cqfd

\section{Contractions}

On se propose dans cette section de donner un \'enonc\'e
permettant de v\'erifier  que sous certaines conditions, un morphisme est
une contraction sur son image, et donc  g\'en\'eriquement injectif.  Un
morphisme projectif et plat de fibre
$\proj_1$  est appel\'e
$\proj_1-$fibration.

\begin{prop} \label{cle2} Soient $\beta:M\hfl{}{} {\cal C}$  un morphisme
propre de vari\'et\'e alg\'ebriques irr\'eductibles et lisses, et $Y\subset {\cal
C}$  une sous-vari\'et\'e ferm\'ee int\`egre de
${\cal C}.$ On suppose que les conditions suivantes sont satisfaites :

{\rm (i)} Au-dessus  de l'ouvert ${\cal C}\setminus Y,$  le morphisme
$\beta$  est fini.

{\rm (ii)} Le  sous-sch\'ema $D=\beta^{-1}(Y)$  est un diviseur de $M.$

{\rm (iii)}  Le morphisme induit  $\beta|_{D}:D\hfl{}{}Y$ 
 est une $\proj_1$-fibration, et la restriction de ${\cal O}(D)$  \`a la
fibre est de degr\'e $-k<0.$

{\rm (iv)}  Le morphisme  induit par la diff\'erentielle sur les faisceaux
conormaux
$${\cal N}_{Y/{\cal C}}\hfl{}{} \beta_{*}({\cal N}_{D/M})$$
est surjectif.

Alors le morphisme $\beta$  est g\'en\'eriquement injectif.
\end{prop}

Pour d\'emontrer ce r\'esultat, on consid\`ere la factorisation  de  Stein 
 de
$\beta$ et le diagramme commutatif suivant qui en d\'ecoule:
$$\diagram{
D&\hookrightarrow & M\cr
\vfl{\beta|_{D}}{}&&\vfl{}{\widetilde{\beta}}\cr
Y&\hookrightarrow &\widetilde{{\cal C}}\cr
&\sefl{}{}&\vfl{}{\nu}\cr
&&{\cal C}\cr
}$$
Par construction,  $\nu$  est un morphisme fini, et on a
$\widetilde{\beta}_{*}({\cal O}_{M})={\cal O}_{\widetilde{\cal C}}.$ On
obtient en consid\'erant le faisceau conormal ${\cal N}_{Y/{\cal C}}$ de $Y$ 
dans
$\widetilde{\cal C}$ la factorisation  de $^{t}d\beta:$
$$\diagram{{\cal N}_{Y/{\cal C}}&\hfl{^{t}d\nu}{}&
{\cal N}_{Y/\widetilde{{\cal C}}}\cr
&\sefl{}{^{t}d\beta}&\vfl{}{^{t}d\widetilde{\beta}}\cr
&&\beta_{*}(N^*_{D/M})\cr}$$
L'assertion (ii) du lemme
 \ref{f-formelles} ci-dessous montre que sur un ouvert non vide de
$Y,$ 
le morphisme   $^{t}d\widetilde{\beta}$ 
 est g\'en\'eriquement un isomorphisme. D'apr\`es la condition (iv),
la fl\`eche horizontale est donc g\'en\'eriquement surjective. Il en r\'esulte que la
diff\'erentielle
$d\nu$ est g\'en\'eriquement injective le long de $Y$  et par suite  le
   morphisme $\nu$ est g\'en\'eriquement injectif: il suffit  en effet  d'appliquer
le lemme 7.4  du chapitre II  de Hartshorne [8]  au morphisme  d'anneaux
locaux de ${\cal O}_{{\cal C},\xi}\hfl{}{}{\cal O}_{\widetilde{\cal C},\xi
}$ induit par
$\nu,$ o\`u $\xi$   d\'esigne le point g\'en\'erique de $Y$. Le morphisme
$\widetilde{\beta}$  \'etant birationnel, on obtient que le morphisme $\beta$ 
est aussi g\'en\'eriquement injectif,  ce qui ach\`eve la d\'emonstration de la
proposition \ref{cle2}.

Il reste \`a \'etudier le faisceau conormal ${\cal N}_{Y/\widetilde{\cal C}}.$

\begin{lemme} \label{f-formelles} On suppose satisfaites les conditions
{\rm (i),(ii)}, et {\rm (iii)}  de la proposition \ref{cle2}. Soit
$Y^{[\ell]}$ le voisinage infinit\'esimal d'ordre
$\ell$ de
$Y$ dans $\widetilde{\cal C}$
 et $D^{[\ell]}$  le voisinage
infinit\'esimal d'ordre $\ell$  de $D$  dans $M.$

 {\rm (i)} Le morphisme 
$${\cal O}_{Y^{[\ell]}}\hfl{}{} \widetilde{\beta}_{*}({\cal O}_{D^{[\ell]}})$$
est g\'en\'eriquement un isomorphisme.

{\rm (ii)} Le morphisme $$^{t}d\widetilde{\beta}: {\cal N}_{Y/\widetilde{{\cal
C}}}\hfl{}{}\widetilde{\beta}_{*}(N^*_{D/M})$$  est g\'en\'eriquement
un isomorphisme.

\end{lemme}

\begin{proof} 
Ce lemme  est une cons\'equence du th\'eor\`eme des fonctions
formelles (\cf Hartshorne [8]). 
D\'esignons par $I_{Y}$  l'id\'eal de $Y$  dans $\widetilde{{\cal C}}_{n}$, et par
$I_{D}$ l'id\'eal de $D$  dans $M_{n}.$ L'image r\'eciproque du sous-sch\'ema $Y$ 
de $\widetilde{\cal C}$  par le morphisme $\widetilde{\beta}$  est
contenue dans l'image  r\'eciproque de $Y$ par le morphisme $\beta.$
Elle co\"{\i}ncide donc avec le sous-sch\'ema $D.$ Soit
$\eta$  le point g\'en\'erique de
$Y,$  
et $\hat{\cal O}_{\widetilde{{\cal
C}_{n}},\eta}$ le compl\'et\'e de l'anneau local ${\cal O}_{\widetilde{{\cal
C}_{n}},\eta }$ pour la topologie $I_{\eta}-$adique.
Le th\'eor\`eme des fonctions formelles  montre que  le
morphisme canonique 
$$\hat{\cal O}_{\widetilde{{\cal
C}_{n}},\eta}\hfl{}{}\limproj_{\ell}\widetilde{\beta}_{*}({\cal
O}/I_{D}^{\ell})_{\eta}$$
est un isomorphisme. On sait que fibr\'e  conormal de $N^*_{D/M}$  est, le
long des fibres de la projection $D\hfl{}{} Y$, isomorphe \`a ${\cal O}(k).$ 
D'apr\`es le th\'eor\`eme de changement de base, on a 
$R^{1}\beta_{*}(S^{\ell}N^*_{D/M})=0.$  Il en r\'esulte que  pour tout $\ell\geq 0,$
le morphisme canonique 
$$\limproj_{\ell}\widetilde{\beta}_{*}({\cal
O}/I_{D}^{\ell})_{\eta}\hfl{}{} \widetilde{\beta}_{*}({\cal
O}/I_{D}^{\ell})_{\eta}$$
est lui aussi surjectif, et le diagramme commutatif
$$\diagram{\hat{\cal O}_{\widetilde{{\cal
C}},\eta}&\hfl{}{}&\limproj_{\ell}\widetilde{\beta}_{*}({\cal
O}/I_{D}^{\ell})_{\eta} \cr
\vfl{}{}&&\vfl{}{}\cr
({\cal O}/I_{Y}^{\ell})_{\eta} &\hfl{}{}& \widetilde{\beta}_{*}({\cal
O}/I_{D}^{\ell})_{\eta}\cr
}$$
montre que  pour tout entier $\ell>0$ le morphisme canonique ${\cal
O}/I_{Y}^{\ell}
\hfl{}{}
\widetilde{\beta}_{*}({\cal O}/I_{D}^{\ell})$  est surjectif au point g\'en\'erique
de $Y.$  Pour $\ell=1,$ ce morphisme est induit par la projection $D\hfl{}{} Y$
et par suite, c'est un isomorphisme.  Le diagramme commutatif de suites exactes 
$$\diagram{
0&\hfl{}{}& I_Y/I_{Y}^2&\hfl{}{} &{\cal O}/I_{Y}^2&\hfl{}{}&{\cal
O}/I_{Y}&\hfl{}{}0\cr
&&\vfl{}{}&&\vfl{}{}&&\vfl{}{\wr}&\cr
0&\hfl{}{} &\widetilde{\beta}_{*}(I_D/I_{D}^2)&\hfl{}{}&
\widetilde{\beta}_{*}({\cal O}/I_{D}^2)&\hfl{}{}& \widetilde{\beta}_*({\cal
O}/I_{D})&\hfl{}{} 0\cr
}$$
montre que la premi\`ere fl\`eche verticale est surjective au
point g\'en\'erique de $Y.$  

Montrons que le morphisme ${\cal
O}/I_{Y}^{\ell}
\hfl{}{}
\widetilde{\beta}_{*}({\cal O}/I_{D}^{\ell})$  est injectif au point g\'en\'erique
$\eta$  de $Y.$ Consid\'erons  une section locale $u \in {\cal
O}_{\widetilde{C},\eta}.$ On suppose que 
$\widetilde{\beta}^*(u)\in \widetilde{\beta}_{*}(I_{D}^{\ell})_{\eta}$.
 Ainsi,  $\widetilde{\beta}^*(u) $ s'annule \`a l'ordre $\ell $ le long
de
$D$ \`a l'ordre
$\ell.$ Consid\'erons la classe de $\widetilde{\beta}^*(u) $ dans
$\widetilde{\beta}_{*}(I_{D}^{\ell}/I_{D}^{\ell+1}).$  Le morphisme 
$S^{\ell}\widetilde{\beta}_{*}(I_{D}/I_{D}^2)\hfl{}{}
\widetilde{\beta}_*(I_{D}^{\ell}/I_{D}^{\ell+1})$  est d'apr\`es le th\'eor\`eme de
changement de base, surjectif.  Il en r\'esulte qu'il existe
$u_{\ell}\in  I_{Y,\eta}^{\ell}$  tel que $\widetilde{\beta}^*(u-u_{\ell}) $ 
appartienne \`a
$\widetilde{\beta}_{*}(I_{D}^{\ell+1})_{\eta}.$ Par r\'ecurrence, on obtient en
recommen\c{c}ant la construction pr\'ec\'edente pour $i\geq \ell$ des \'el\'ements
$u_{i}\in I_{Y,\eta}^{\ell+i}$  tel que 
$$\widetilde{\beta}_{*}(u- \sum_{\ell\leq i\leq \ell+k}u_{i})\in
\widetilde{\beta}_{*}(I_{D}^{\ell+k+1})_{\eta}.$$
Autrement dit, dans $\limproj_{\ell}\widetilde{\beta}_{*}({\cal
O}/I_{D}^{\ell})_{\eta}$ on a $$\widetilde{\beta}^*(u-\sum_{i\geq \ell}u_{i})=0.$$
Compte-tenu du th\'eor\`eme des fonctions formelles, cette \'egalit\'e  implique
$u=\sum_{i\geq \ell}u_{i}$ dans le compl\'et\'e
$\hat{\cal O}_{\widetilde{{\cal C}},\eta}.$ Mais alors $u$  appartient au
compl\'et\'e
$\widehat{I^{\ell}_{Y,\eta}}.$  Mais  d'apr\`es Atiyah-Macdonald ([1],
proposition 10.15), on a
$$\widehat{I_{Y,\eta}^{\ell}}\cap {\cal O}_{\widetilde{C},\eta}=
I^{\ell}_{Y,\eta}.$$ Ceci prouve que le morphisme ${\cal O}/I_{Y}^{\ell}
\hfl{}{}
\widetilde{\beta}_{*}({\cal O}/I_{D}^{\ell})$ est injectif au point $\eta.$ 
Ceci d\'emontre l'assertion (i). L'assertion (ii)  d\'ecoule du diagramme
commutatif de suite exactes d\'ej\`a \'ecrit, dans lequel les fl\`eches verticales sont
maintenant  des isomorphismes. \end{proof}

\section{Fin de la d\'emonstration}

On suppose que le th\'eor\`eme \ref{main} est vrai \`a l'ordre $n-1.$
Pour appliquer le r\'esultat  de la section ci-dessus, on consid\`ere un ouvert
${\cal C}_{n}^{*}$ de
${\cal C}_{n}$  satisfaisant aux conditions (a)-(d)  d\'ecrites dans
la d\'emonstration du corollaire \ref{cartesien}.

Consid\'erons  la vari\'et\'e   $Y=\psi^{-1}({\cal C}_{n}^*)$  de
 plong\'ee  comme une sous-vari\'et\'e ferm\'ee de
${\cal C}_{n}^*.$
 Nous allons montrer que les conditions de la
proposition
\ref{cle2}  sont satisfaites avec $k=2,$ et quitte  \`a diminuer au besoin
l'ouvert ${\cal C}_{n}^*$, pour le   morphisme 
$\beta^{-1}({\cal C}_{n}^*)\hfl{}{}{\cal C}_{n}^{*}$ induit par $\beta.$ En
effet, les  points du compl\'ementaire $\Omega$ de $Y$ dans ${\cal
C}_{n}^*$  s'ils proviennent de $M_{n},$ repr\'esentent des courbes 
irr\'eductibles, et le morphisme
$\beta$  est alors quasi-fini  d'apr\`es [13] et propre au dessus de cet
ouvert 
$\Omega$, donc fini, ce qui d\'emontre (i). Le
sous-sch\'ema image r\'eciproque de
$Y$ par $\beta $ est ensemblistement  la trace de $D$;  nous allons
montrer qu'il l'est au sens des sch\'emas.  La difficult\'e est  que
l'ouvert ${\cal C}_{n}^*$  ne contient pas
le point correspondant \`a la courbe  correspondant au sommet $s$  du c\^one
$\gamma(\Sigma)$  construit dans la section \S 12.  Nous allons donc dans
un premier temps travailler au-dessus d'un ouvert ${\cal C}_{n}'$  plus
grand que l'ouvert ${\cal C}_{n}^*$  ci-dessus.

Consid\'erons l'ouvert ${\cal C}_{n}'$  de ${\cal C}_{n}$  des courbes qui
satisfont \`a la condition (c), de  sorte qu'au-dessus de cet ouvert, le
morphisme $j:W={\cal C}_{n-1}\times \proj_2\hfl{}{} {\cal C}_{n}$  d\'efini
par l'addition des diviseurs de $\proj_2^*$  est un plongement  sur une
sous-vari\'et\'e ferm\'ee de ${\cal C}'_{n}.$  D\'esignons par $Y'\subset
U_{n-1}\times \proj_2$ , $M_{n}'\subset M_{n}$
et $D'\subset D$ les images r\'eciproques de l'ouvert ${\cal C}_{n}'$  
 de sorte qu'on d\'eduit du diagramme (10.1) le diagramme
commutatif
$$\eqalign{\diagram{D'&&\hookrightarrow&& M_{n}'\cr
\vfl{\pi}{}&&&&\vfl{}{\beta} \cr
Y'&\hfl{\psi}{}& W'&\subset&{\cal C}_{n}'\cr 
}}\eqno(14.1)\quad$$
La diff\'erentielle $d\beta$  induit pour les faisceaux conormaux un
morphisme de  faisceaux alg\'ebriques coh\'erents au-dessus de $Y'$
$$\psi^*(N^*_{W'/{\cal C}'_{n}})\hfl{}{} \pi_{*}(N^*_{D'/M'_{n}})
\eqno(14.2)\quad$$

\begin{lemme} 

{\rm (i)}  On a $R^{1}\pi_{*}(N^*_{D'/M'_{n}})=0$;

{\rm (ii)}
le faisceau $\pi_{*}(N^*_{D'/M'_{n}})$  est localement libre de rang
3, et sa fibre au-dessus d'un point $z\in Y'$ s'identifie \`a l'espace
vectoriel
$H^{0}(N^*_{D'/M'_n}(z))$  des sections  de la restriction
$N^*_{D'/M'_n}(z)$ du fibr\'e conormal $N^*_{D'/M'_n}$ \`a la fibre
$\pi^{-1}(z).$

{\rm (iii)} le morphisme de faisceaux localement libres
$\psi^*(N^*_{W'/{\cal C'}_{n}})\hfl{}{} \pi_{*}(N^*_{D'/M'_{n}}) $ induit
par l'application cotangente \`a $\beta$ est g\'en\'eriquement surjectif. 

\end{lemme}

\begin{proof} Les assertions (i) et (ii)  r\'esultent du fait que le
morphisme $\pi: D\hfl{}{} U_{n-1}\times \proj_2$  est un morphisme
projectif  lisse  dont la fibre s'identifie \`a la droite $\proj_1;$  en
particulier, il est  plat, et de plus, la restriction du fibr\'e conormal
$N^*_{D/M_{n}}$ \`a une fibre est isomorphe \`a ${\cal O}(2).$ Le m\^eme
r\'esultat est \'evidemment vrai pour la fibration  $D'\hfl{}{} Y'$ induite
par $\pi.$ Pour (iii)  on consid\`ere la fibre de
$\pi$ au-dessus du point
$z=(E,x_{0})$, o\`u $E$  est un 
fibr\'e de Poncelet de $U_{n-1} $ associ\'e \`a un syst\`eme lin\'eaire
$(\Gamma,\Theta)$ sur la conique lisse $C$ et
$x_0$ un point de $\proj_2$ suffisamment g\'en\'eraux pour que les
conclusions du lemme \ref{non-ramif}  soient vraies. Puisque la courbe de
Poncelet $\beta_{E}$  est lisse,  ce point appartient \`a l'ouvert $Y'.$
 Consid\'erons la  surface  de Hirzebruch $\Sigma$ associ\'ee 
et son image $f(\Sigma)$, qui rencontre $D$ transversalement  suivant 
la fibre de $\pi$  au-dessus de $z$; consid\'erons la restriction du
morphisme \`a la fibre
$\pi^{-1}(z)$
$$\diagram{N^*_{W'/{\cal C}'_{n}}(z)&\hfl{}{}&
\pi_{*}(N^*_{D'/M_{n}'})(z) \cr
&& \vfl{}{}\cr
&&H^0(N^*_{D'/M'_{n}}(z))\cr}$$
L'application lin\'eaire compos\'ee  est surjective: en   effet, si 3
points
$a_{i}$ sont choisis sur la droite projective $\proj_1=\pi^{-1}(z)$ 
formant un rep\`ere, du fait que $N^*_{D'/M'_{n}}(z)={\cal
O}_{\proj_1}(2)$, le morphisme de restriction
$H^0(N^*_{D'/M'_{n}}(z)) \simeq \prod_{i}N_{D'/M'_{n}}(a_{i})$ est
un isomorphisme. L'application lin\'eaire compos\'ee $$N^*_{W'/{\cal
C}'_{n}}(z)\hfl{}{}
\prod_{i}N^*_{D/M_{n}}(a_{i})$$  est la transpos\'ee de l'application
lin\'eaire 
$\oplus_{i=1,2,3} N_{D'/M'_{n}}(a_{i}) \hfl{}{} N_{W'/{\cal C}_{n}}(z)$
induite par l'application lin\'eaire tangente $d_{a_{i}}\beta$ aux points
$a_{i}.$ Cette application lin\'eaire se factorise suivant le diagramme
$$\diagram{
\oplus_{i=1,2,3} N_{D'/M'_{n}}(a_{i})&\hfl{}{}& T_{s}R \cr
&\sefl{}{}&\vfl{}{}\cr
&&N_{W'/{\cal C}_{n}}(z)\cr
}$$
o\`u $R$  est le c\^one quadratique de sommet $s$  construit dans le
corollaire \ref{vrai.cone}. 
Compte-tenu de la proposition \ref{transverse}, la fl\`eche horizontale
s'identifie \`a l'application lin\'eaire
$\oplus_{i}N_{\Sigma_{0}/\Sigma}(a_{i})\hfl{}{} T_{s}R$ induite par la
diff\'erentielle de $\gamma:$ c'est un isomorphisme. D'apr\`es le lemme
\ref{non-ramif}, la fl\`eche compos\'ee  est injective. Donc la
fl\`eche verticale est injective. Le lemme en d\'ecoule, compte-tenu du
fait que $Y'$ est irr\'eductible, et du fait que les points au-dessus
desquels un morphisme de fibr\'es vectoriels est surjectif forment un
ouvert.
\end{proof}
\bigskip\goodbreak

Fin de la v\'erification.

Consid\'erons le diagramme induit au-dessus de l'ouvert ${\cal C}_{n}^*$
$$\eqalign{\diagram{D^*&\hookrightarrow& M_{n}^*\cr
\vfl{\pi}{}&&\vfl{}{\beta} \cr
Y&\hfl{\psi}{}&{\cal C}_{n}^{*}\cr 
}}\eqno(14.3)\quad$$
dont on sait  d\'ej\`a qu'il est  cart\'esien au sens ensembliste.
La condition (iv) de la  proposition 13.1 est satisfaite sur un ouvert
non vide de
$Y;$ quitte \`a diminuer au besoin  l'ouvert ${\cal C}_{n}^*,$  on peut 
supposer que cette condition (iv) est vraie sur $Y$  tout entier.  Il
en r\'esulte que pour
$a\in D^*$ la diff\'erentielle $d_a\beta$ induit sur les espaces normaux une
application lin\'eaire
$N_{D/M_{n}}(a)\hfl{}{} N_{Y/{\cal C}_{n}}(\pi(a))$  injective.  Par
suite,  l'espace tangent de Zariski  en $a$  au sous-sch\'ema
$\beta^{-1}(Y)$  est r\'eduit \`a l'espace tangent $T_{a}D.$  Ainsi, le sous-sch\'ema $\beta^{-1}(Y)$
est r\'eduit \`a $D^*$ et la condition (ii) de la proposition \ref{cle2}  est
satisfaite. La condition
(iii) \'etant  vraie d'apr\`es le lemme
\ref{bord},  nous sommes exactement dans les conditions d'application
de la proposition \ref{cle2}.  Donc le morphisme $\beta$  est
g\'en\'eriquement injectif.

\bigskip
\centerline{\capit Remerciements}
\medskip
 Le premier 	auteur
remercie  l'ICTP (Trieste) qui l'a invit\'e \`a
pr\'esenter   ce travail \`a l'Ecole d'\'et\'e de g\'eom\'etrie alg\'ebrique en 1999. Il
remercie aussi  le CEFIPRA qui a subventionn\'e son voyage en Inde et  permis
d'exposer les r\'esultats de cet article
 \`a Bombay et Madras  en janvier 2000 dans le cadre du Contrat de 
coop\'eration 1601-2 (Geometry).  
  Une partie
 de cet article a \'et\'e \'ecrite  lors du s\'ejour du second auteur  \`a
Bonn (octobre 1998-mars 1999); il remercie le Max Planck Institute of
Mathematics de son accueil.

\bigskip
\centerline{\capit R\'ef\'erences}

\medskip
\item{\rm [1]}{\capit M.F. Atiyah and I.G.  Macdonald}, {\it Introduction
to commutative algebra,}  Addison-Wesley Publishing Company (1969).

\item{\rm [2]} {\capit W. Barth,} {\ \it Some Properties of Rank-2
Vector Bundles on $\proj_n$}, Inventiones  math., {\bf 226}, 1977, p.
125-150.

\item {\rm [3]} {\capit W. Barth,}{\it\  Moduli of vector bundles on
projective  plane,}{\rm\   Inventiones math. }{\bf\  42 }{\rm\  (1977)}{\ p.
63-91.}

\item{\rm [4]} {\capit H. Bateman,}{\it The quartic curve and its
inscribed configurations,}  American J. Math {\bf 36} (1914) p. 357-386.

\item{\rm [5]} {\capit J.M. Dr\'ezet,} {\it Groupe de Picard des vari\'et\'es
de modules de faisceaux semi-stables sur $\proj_2(\comp)$}, Ann. de
l'Institut Fourier {\bf 38}, (1988) p. 105-168.

\item{\rm [6]} {\capit  G. Ellingsrud, L. G\"ottsche,} {\it Variations of
moduli spaces and Donaldson invariants under change of polarization}, J.
reine angew. Math. Math {\bf 467} (1995)  p. 1-49.

\item{\rm [7]} {\capit G. Ellingsrud, J. Le Potier, S.A. Str\o mme,}
{\it Some Donaldson invariant of $\proj_{2}(\comp)$} Proceedings   du
Symposium Taniguchi (d\'ecembre 1994, Kyoto), Lecture Notes in pure and
applied mathematics, ed. M. Maruyama {\bf 179} (1996) p. 33-38.

\item{\rm [8]} {\capit  R. Hartshorne}, {\it Algebraic geometry},
Springer (1977).

\item{\rm [9]} {\capit M. He}, {\it Espace de modules de syst\`emes
coh\'erents}, Int. J. of Maths, Vol 9, 5 (1998) p.545-598.

\item{\rm [10]}{\capit K. Hulek et J. Le Potier}, {\it Sur l'espace de
modules des faisceaux stables de rang 2, de classes de Chern (0,3) sur
$\proj_2$,}  Annales de l'Institut Fourier, {\bf 39} (1989) p. 251-292.

\item{\rm [11]} {\capit S. Lang}, {\it Algebra}, Addison-Wesley (1984).

\item {[12]} {\capit J. Le Potier,}{\it\  Fibr\'es stables sur le plan
projectif
 et quartiques de L\"uroth,}{\rm\   Expos\'e donn\'e \`a Jussieu }{\rm\ 
(1989).
 }

\item{\rm [13]} {\capit J. Le Potier,} {\it 	Fibr\'e d\'eterminant  et
courbes de saut sur les surfaces alg\'ebriques} ;  Proceedings de la
Conf\'erence de g\'eom\'etrie alg\'ebrique de Bergen (juillet 1989), Complex
projective geometry, London Mathematical Society, Lecture Notes Series
179, p. 213-240.

\item{\rm [14]} {\capit J. Le Potier,} {\it Faisceaux semi-stables de
dimension 1 sur le plan projectif}, Revue Roumaine de Math\'ema\-tiques
Pures et Appliqu\'ees, 
  d\'edicac\'e \`a la m\'emoire de Constantin B\u anic\u a; {\bf 38} (1993) p.
635-678.
\item{\rm [15]}	{\capit J. Le Potier,} {\it Syst\`emes coh\'erents et
structures de niveau}; 
  Ast\'erisque n$^0$ 214, 143 p. (1993). 
	
\item{\rm [16]} {\capit J. Le Potier,} {\it Faisceaux semi-stables et
syst\`emes coh\'erents,}
 Proceedings de  la Conf\'erence de Durham, Cambridge University Press
(1995) p.179-239. 

\item{\rm [17]} {\capit J. Le Potier,} {\it Syst\`emes coh\'erents et
polyn\^omes de Donaldson}, Proceedings   du Symposium Taniguchi (d\'ecembre
1994, Kyoto), Lecture Notes in pure and applied mathematics, ed. M.
Maruyama {\bf 179} (1996) p. 103-128. 

 \item{\rm [18]} {\capit W.-P. Li, Z. Qin,}{\it\  Lower-degree
Donalson polynomial invariants of rational surfaces,}{\rm\   J.
Algebraic geometry }{\bf\  2 }{\rm\  (1993) }{\ p. 413-442. }

\item{\rm [19]} {\capit M. Maruyama,}{\it\ 
Moduli of stable sheaves, II,}{\rm\   J. Math. Kyoto University }{\bf\ 
18 }{\rm\  (1978) }{\ p. 557-614. }

\item {\rm [20]} {\capit M. Maruyama,}{\it\  Singularities of the curve
of jumping  lines of a vector bundle of rank 2 on $\proj_2$}, Algebraic
Geometry, Proc. of  Japan-France Conf. (1982),{\rm\   Lectures Notes in
Math., Springer } {\bf\  1016 }{\rm\  (1983) }{\ p. 370-411}.

\item {\rm [21]} {\capit F. Morley,}{\it\  On the L\"uroth quartic
curve,}{\rm\   Amer. J. Maths }{\bf\  36 }{\rm\  (1918) p. 357-386.}

\item{\rm [22]} {\capit C. T. Simpson,}{\it\  Moduli
of Representations of the Fundamental Group of a Smooth
Variety, I},{\rm\  Publications math\'ematiques de l'IHES, 79
}{\bf\  }{\rm\  (1994) p. 47-129.}

\item{\rm [23]} {\capit S.A. Str\o mme,} {\it Ample divisors on fine
moduli spaces on the projective plane,} Math. Z. {\bf 187} (1984), p.
405-423.

\item{\rm [24]} {\capit M. Toma}, {\it Birational models for varieties
of Poncelet curves,} Manuscripta Math. {\bf 90} (1996) n$^0 1$
p.105-119.

\item {\rm [25]} {\capit G. Trautmann,}{\it\  Poncelet curves and
associated theta  characteristics,}{\rm\   Expositiones Mathematicae
}{\bf\  6 }{\rm\  (1988) } {\ p. 29-64.}

\item {\rm [26]} {\capit A. N. Tyurin}{\it\  The
moduli spaces
 of vector bundles on threefolds, surfaces and curves,} {\rm\  
Forschungsschwerpunkt ``Komplexe Mannigfaltigkeiten'', Schriftenreihe,
Helf 67,  Preprint Erlangen }{\rm\  (1990)}.

\item{\rm [27]} {\capit J. Vall\`es,} {\it Diviseurs inattendus de
droites sauteuses. Fibr\'es de Schwarzenberger,} Th\`ese, Universit\'e Paris 6
(1996).
\bigskip

\noindent Joseph Le Potier\\
Institut de Math\'ematiques de Jussieu\\
Universit\'e Paris 7,
Case Postale 7012\\
2, Place Jussieu, 
75251 Paris-Cedex 05\\
e.mail: jlp@math.jussieu.fr
\medskip

\noindent Alexander Tikhomirov \\
Department of Mathematics,\\
State Pedagogical University of Yaroslavl\\
Respublikanskaya str. 108, Yaroslavl, Russia\\
e.mail: alexandr@tikho.yaroslavl.su
\end{document}